\documentclass[a4paper,11pt]{article}

\usepackage{geometry}										% [showframe]
\usepackage[english]{babel}									% language
\usepackage[utf8]{inputenc}									% input text encoding
\usepackage[T1]{fontenc}									% output encoding

\usepackage{lmodern}										% Latin Modern

\usepackage{xcolor}
\usepackage[normalem]{ulem}
\usepackage{siunitx}
\usepackage{placeins} % for FloatBarrier
\usepackage{bbold} % for \mathbb{1}
\usepackage{algorithm}
\usepackage{algpseudocode}

\usepackage[makeroom]{cancel}
\usepackage{rotating} % for table in landscape orientation

\usepackage{float} % for floating verbatim section
\floatstyle{boxed}
\newfloat{File}{htbp}{}

\usepackage{amsmath}										% improved mathematical typesetting
\usepackage{amstext}										% text in mathematics mode
\usepackage{amsfonts}										% extended set of fonts
\usepackage{amssymb}										% extended symbol collection

\usepackage{authblk}										% [auth-sc,affil-it], to be loaded after titling
\title{Sparse Polynomial Chaos Expansions: Literature Survey and Benchmark}
\author[1]{Nora L\"uthen}
\author[1]{Stefano Marelli}
\author[1]{Bruno Sudret}
\affil[1]{Chair of Risk, Safety and Uncertainty Quantification, ETH Z\"{u}rich, Stefano-Franscini-Platz 5, 8093 Z\"{u}rich, Switzerland}

\date{\today}

\usepackage{caption}										% customized captions for floating environments
\usepackage[caption=false]{subfig}

\usepackage{natbib} % ,sort&compress % number ranges are compressed and hyphenated
\usepackage{hyperref}
\hypersetup{
	pdftitle = {Sparse Polynomial Chaos Expansions: Literature Survey and Benchmark},
	pdfauthor = {Nora L\"uthen, Stefano Marelli, Bruno Sudret},
	pdfsubject = {Manuscript submitted to SIAM/ASA J. Uncertainty Quantification},
	pdfkeywords = {uncertainty quantification, surrogate modelling, sparse regression, sparse polynomial chaos expansions, compressive sensing, experimental design, literature review},
	colorlinks = true,										% false: boxed links; true: colored links
	linkcolor={red!90!black},
	citecolor={green!50!black},
	urlcolor={red!50!black}
}

\newcommand{\ca}{{\mathcal A}}

\newcommand{\cd}{{\mathcal D}}

\newcommand{\cl}{{\mathcal L}}
\newcommand{\cm}{{\mathcal M}}
\newcommand{\cn}{{\mathcal N}}
\newcommand{\co}{{\mathcal O}}

\newcommand{\cs}{{\mathcal S}}

\newcommand{\cu}{{\mathcal U}}

\newcommand{\cx}{{\mathcal X}}

\newcommand{\di}[1]{{\rm d}#1} 				% d droit pour les integrandes
\newcommand{\ve}[1]{\boldsymbol{#1}}			% Vecteur
\newcommand{\enum}{ , \, \dots \,,}
\newcommand{\norme}[2]{\parallel #1\parallel_{#2}}

\newcommand{\Varhat}[1]{{\rm \hat{Var}}\left[ #1 \right]}	% estimateur de la variance

\newcommand{\Esp}[1]{{\mathbb E}\left[ #1 \right]}
\newcommand{\Espe}[2]{{\mathbb E}_{#1}\left[#2\right]}
\newcommand{\Vare}[2]{{\rm Var}_{#1}\left[#2\right]}

\newcommand{\Prob}[1]{{\mathbb P}\left( #1 \right)}	% probabilite de ()

\newcommand{\R}{\mathbb{R}}
\newcommand{\N}{\mathbb{N}}
\newcommand{\alp}{{\ve{\alpha}}}

\renewcommand{\norme}[2]{\left\| #1 \right\|_{#2}}

\newcommand{\NED}{N}
\newcommand{\Nval}{N_\text{val}}
\newcommand{\ival}{^{(i)}_\text{\tiny val}}

\newcommand{\Graycomment}[1]{ \hfill \textcolor{gray}{  $\triangleright$ #1} }

\newcommand{\markOwn}{{*}}
\newcommand{\SPloo}{$\text{SP}_\text{LOO}${}}

\newcommand{\changed}[1]{{#1}}
\newcommand{\changedmath}[1]{{#1}}

\begin{document}

\maketitle

% REQUIRED
\begin{abstract}
	Sparse polynomial chaos expansions (PCE) are a popular surrogate modelling method that takes advantage of the properties of PCE, the sparsity-of-effects principle, and powerful sparse regression solvers to approximate computer models with many input parameters, relying on only few model evaluations.
	Within the last decade, a large number of algorithms for the computation of sparse PCE have been published in the applied math and engineering literature. 
	We present an extensive review of the existing methods and develop a framework for classifying the algorithms. Furthermore, we conduct a unique benchmark on a selection of methods to identify which approaches work best in practical applications. 
	Comparing their accuracy on several benchmark models of varying dimensionality and complexity, we find that the choice of sparse regression solver and sampling scheme for the computation of a sparse PCE surrogate can make a significant difference, of up to several orders of magnitude in the resulting mean-squared error. Different methods seem to be superior in different regimes of model dimensionality and experimental design size. 
	 
\end{abstract}

%% REQUIRED
%\begin{keywords}
%	uncertainty quantification, surrogate modelling, sparse regression, sparse polynomial chaos expansions, experimental design 
%\end{keywords}

%% REQUIRED
%\begin{AMS}
%%  68Q25, 68R10, 68U05
%\end{AMS}

\section{Introduction}
\label{sec:intro}

Computer modelling is used in nearly every field of science and engineering. Often, these computer codes model complex phenomena, have many input parameters, and are \changed{expensive} to evaluate. 
In order to explore the behavior of the model under uncertainty (e.g., uncertainty propagation, parameter calibration from data or sensitivity analysis), many model runs are required. However, if the model is costly, only a few model evaluations can be afforded, which often do not suffice for thorough uncertainty quantification.
In engineering and applied sciences, a popular work-around in this situation is to construct a surrogate model.
A surrogate model is a cheap-to-evaluate proxy to the original model, which typically can be constructed from a relatively small number of model evaluations and approximates the input-output relation of the original model well. 
Since the surrogate model is cheap to evaluate, uncertainty quantification can be performed at a low cost by using the surrogate model instead of the original model.
Therefore, surrogate modelling aims at constructing a metamodel that provides an accurate approximation to the original model while requiring as few model evaluations as possible for its construction.

In this article, we focus on \textit{nonintrusive regression-based sparse polynomial chaos expansions} (PCE), which is a popular surrogate modelling technique, and within the last decade is has received attention from the communities of applied mathematics and engineering.
PCE express the computational model in terms of a basis of polynomials orthonormal with respect to the input random variables \citep{Xiu2002} and work well for globally smooth problems, which are common in many engineering applications. In addition to being a surrogate model, PCE are also often used for uncertainty propagation and sensitivity analysis, since moments and Sobol' sensitivity indices can be computed analytically \citep{SudretRESS2008b}.
\textit{Nonintrusive} PCE treat the model as a black box (unlike intrusive PCE commonly used for solving stochastic PDEs). 
It is often advantageous to compute a \textit{sparse PCE}, which is an expansion for which most coefficients are zero. This can be justified by the \textit{sparsity-of-effects principle} and by \textit{compressibility}:
The sparsity-of-effects principle is a heuristic stating that most models describing physical phenomena are dominated by main effects and interactions of low order \citep{Montgomery:2004}. 
Furthermore, PCE of real-world models are usually either sparse or at least compressible, meaning that the PCE coefficients, sorted by magnitude, decay quickly. 
Additional advantages of sparse expansions are given in Section \ref{sec:sparsePCEdescription}.

Within the last decade, a large number of articles has been published on the topic of regression-based sparse PCE, each containing promising improvements on how to perform sparse PCE but often lacking a thorough comparison to previously published methods.
In this work, we survey the state-of-the-art literature, develop a general framework into which the various approaches can be fit, and carry out a numerical benchmark of a selection of methods to assess which of the many sparse PCE methods perform best on a representative set of realistic benchmark \changed{models}.

The paper is structured as follows. Section~\ref{sec:sparsePCE} contains the description of our framework for classifying the sparse PCE literature as well as the extensive literature review.
Section~\ref{sec:benchmark} contains the benchmark description and the numerical results. Finally, conclusions are drawn in section~\ref{sec:conclusion}. More detailed descriptions of selected sparse solvers and experimental design techniques are given in the Appendices using unified notation.

\section{Framework and literature survey for sparse polynomial chaos expansions}
\label{sec:background}
\label{sec:sparsePCE}

\subsection{Regression-based polynomial chaos expansions}
\label{sec:PCE}

Let $\ve X$ be a $d$-dimensional random vector on a domain $\cd \subset \R^d$ with independent %
components and probability density function $f_{\ve X}(\ve x) = \prod_{i=1}^d f_{X_i}(x_i)$.
Let $L_{f_{\ve X}}^2(\cd)$ be the space of all scalar-valued 
models with finite variance under $f_{\ve X}$, i.e.,
$
L_{f_{\ve X}}^2(\cd) = \{ h: \cd \to \R \ | \ \Vare{\ve X}{h(\ve X)} < + \infty \}
$.
Under certain assumptions on the input distribution $f_{\ve X}$ \citep{Xiu2002,Ernst2012}, there exists a polynomial orthonormal basis $\{\psi_\alp: \alp \in \N^d\}$ for $L_{f_{\ve X}}^2(\cd)$. Since the components of $\ve X$ are assumed to be independent, the basis elements are products of univariate orthonormal polynomials and are characterized by the multi-index $\alp \in \N^d$ of polynomial degrees in each dimension.

We consider a particular model $\cm \in L_{f_{\ve X}}^2(\cd)$
and denote by $Y = \cm(\ve X)$ the corresponding output random \changed{variable}.
$Y$ can be represented exactly through an infinite expansion in $\{\psi_\alp: \alp \in \N^d\}$. 
In practice, however, not all infinitely many coefficients can be computed, and we are interested in a truncated expansion
\begin{equation}
Y = \cm(\ve X) \approx \cm^\text{PCE}(\ve X) = \sum_{\alp \in \ca} c_\alp \psi_\alp(\ve X)
\end{equation}
whose accuracy depends on the choice of the finite set $\ca \subset \N^d$ (i.e., on the basis elements used for the expansion) as well as on the coefficients $c_\alp$.
Several truncation techniques are described in Section~\ref{sec:basis}.

To compute the coefficients, one well-known and practical approach is regression \citep{Isukapalli, Berveiller2006}%
\footnote{Other, earlier approaches for computing the coefficients are stochastic Galerkin and stochastic collocation methods \citep{Ghanembook1991, Xiu:Hesthaven:2005, Shen2020}. For a comparison of their performance to regression-based PCE, see, e.g., \citet{BerveillerThesis,Hosder2007,Doostan2011,Mathelin2012}.}%
. 
The basic regression approach is \textit{ordinary least squares} (OLS). 
Let $\{\ve x^{(k)}\}_{k=1}^{\NED} \subset \cd$ be a sample of the input space called \textit{experimental design} (ED). Let $\ve y = (y^{(1)} \enum y^{(N)})^T$ be the vector of model responses with $y^{(k)} = \cm(\ve x^{(k)})$. Define the matrix of basis function evaluations $\ve \Psi$ with entries $\Psi_{ij} = \psi_{j}(\ve x^{(i)})$, where the basis functions are enumerated in an arbitrary way. 
Denoting the number of basis functions with $P$, we see that the \textit{regression matrix} $\ve\Psi$ is an $N \times P$-matrix. 
Then, the OLS regression problem can be written as 
\begin{equation}
\hat{\ve c} = \arg\min_{\ve c \in \R^P} \norme{\ve \Psi \ve c - \ve y}{2}.
\label{eq:OLS}
\end{equation}
For a unique and robust solution, a heuristic number of model evaluations is $N \approx 2P,3P$ \citep{Hosder2007,FajraouiMarelli2017}, which can be infeasible for high-dimensional or high-degree PCE approximations.

\subsection{Sparse PCE}
\label{sec:sparsePCEdescription}

Sparse coefficient vectors are determined through \textit{sparse regression}, which, in addition to a good regression fit, requires that the solution be sparse.
This constraint on sparsity is realized e.g.\ by adding as a regularization term the $\ell^0$-``norm'' or the $\ell^1$-norm of the coefficient vector to the OLS formulation of \eqref{eq:OLS} (see Appendix \ref{app:solvers} for more details).
Many sparse regression methods used in PCE were originally developed in the context of compressive sensing \citep{Donoho2006, Candes2006}.
For an introduction to the concepts and ideas of compressive sensing, see, e.g., \citet{Candes2008a, Bruckstein2009, Kougioumtzoglou2020}.

Unlike OLS, compressive sensing methods allow one to use fewer design points $\NED$ than basis functions $P$ and still recover the true sparse solution, or find a good sparse approximation to it.
This and its robustness to noise, both of which are induced by the sparsity constraint, are the main reasons why sparse PCE are preferred to full PCE in practical settings when the number of model evaluations necessary for OLS-based PCE would be infeasible to compute.
Note that while the use of sparse PCE for engineering models can be justified by compressibility and the sparsity-of-effects heuristic (section \ref{sec:intro}), the main goal in sparse PCE is to compute a good surrogate model from a few model evaluations and not to find the sparsest possible expansion. 
The assumption of sparsity is used as a tool for finding robust solutions to underdetermined systems of linear equations.

The first publications on sparse regression-based PCE proposed greedy forward-backward selection algorithms \citep{BlatmanCras2008, BlatmanPEM2010} and introduced the LARS algorithm for sparse PCE \citep{BlatmanJCP2011}. On the mathematical side, \citet{Doostan2011} analyzes convergence properties for sparse Legendre PCE when the design points are sampled from the uniform distribution.
Another early work is \citet{Mathelin2012} demonstrating that sparse PCE are less costly and more accurate than PCE based on Smolyak sparse grids. 
Since then, a large number of articles has been published on the topic of sparse PCE suggesting new methods for specific aspects of the sparse regression procedure. In the following, we present a framework into which the existing literature can be fit. The framework provides an overview of the available choices and enables a structured comparison of their impact on the performance of the resulting sparse PCE. Naturally, some new combinations of methods arise that have not yet been considered in the literature.

\subsection{Framework: Classifying the literature on sparse PCE}
\label{sec:framework}
\begin{figure}
	\centering
	\includegraphics[width=.5\textwidth]{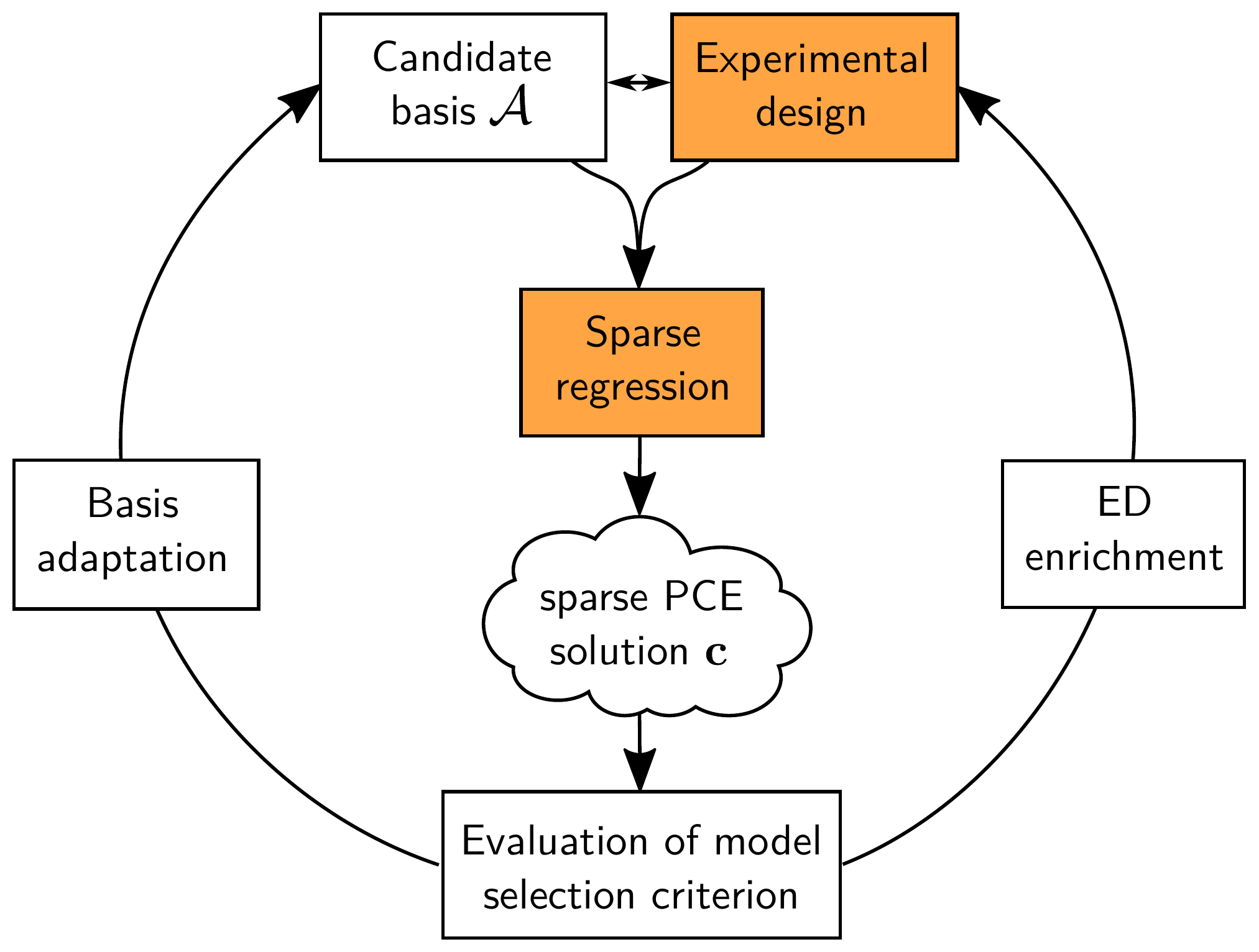}
	\caption{Framework for computing sparse PCE. For each component of the framework, a number of methods has been proposed in the literature. In the first part (section~\ref{sec:background}) of this paper, we review the literature for each of the components. Details on selected methods are given in Appendices \ref{app:ED} and \ref{app:solvers}. In the second part (section~\ref{sec:benchmark}), we conduct a benchmark of selected methods for the components marked in orange, performing a single iteration of the framework. Iterative basis adaptation and experimental design enrichment are not considered in this work and are left for future benchmarks.}
	\label{fig:framework}
\end{figure}

Here, we present the framework we developed in order to gain an overview of the extensive literature proposing new methods for computing sparse PCE.
Figure~\ref{fig:framework} shows a sketch of this framework. To compute a sparse PCE, the first step is to choose a set $\ca$ of candidate polynomials for the expansion (Section \ref{sec:basis}) as well as an experimental design (Section \ref{sec:ED}). 
The experimental design defines the locations of the model evaluations.
Once the model evaluations are obtained, the sparse solution can be computed by applying a sparse regression solver (see Section \ref{sec:solvers}). This solver often depends on a number of hyperparameters that have to be selected carefully in order to get good results. 
Then, a suitable model selection criterion is evaluated (Section \ref{sec:modelselection}). If the obtained solution is satisfactory, the process can be stopped. 
Otherwise, the basis can be adapted (usually augmented; see Section~\ref{sec:basis}), and/or the experimental design can be enriched (see Section~\ref{sec:EDenrichment}).
This process is repeated until the value of the model selection criterion either is satisfactory or cannot be reduced.

In addition to the components shown in Figure~\ref{fig:framework}, there are methods (we call them \textit{enhancements}) that aim at generally improving the solution to the sparse regression problem by, e.g., adapting the input space or preconditioning the regression matrix. They are discussed in Section~\ref{sec:enhancements}.

\subsection{Choice of basis and basis adaptation}
\label{sec:basis}
\label{sec:iterativealgs}

The approximation quality of a truncated PCE hinges on the polynomial functions available for building the surrogate model, which are characterized by the associated set of multi-indices $\ca \subset \N^d$.
We call the finite set of polynomials $\{\psi_\alp: \alp \in \ca\}$ included in the current truncated PCE model \textit{basis} and call its members \textit{candidate polynomials} or \textit{candidate basis functions}. 
A sparse PCE algorithm will find nonzero coefficients only for a subset $\ca^\text{active} \subset \ca$ of the basis functions, which we call \textit{active basis functions}.
On the one hand, $\ca$ should include enough candidate polynomials to facilitate a good approximation. On the other hand, unnecessary basis functions decrease the ratio $N/P$ of model evaluations to unknown coefficients and
deteriorate the regression matrix. 
Therefore, it is beneficial to carefully select the polynomials to be included in the expansion.

The choice of basis is often motivated by the \textit{sparsity-of-effects principle}, a heuristic guideline stating that most real-world models are well approximated by terms of low degree and low interaction order.
The following are popular ways to construct a basis:
\begin{itemize}
	\item \textit{Total-degree} A total-degree basis of degree $p$ is defined by  
	$ \ca^p = \{\alp: \norme{\alp}{1} \leq p \} $. 
	\item \textit{Hyperbolic truncation} 
	Let $p$ be fixed. Define the q-norm-truncated basis 
	\begin{equation}
	\ca^{p,q} = \{ \alp: \norme{\alp}{q} \leq p \}
	\end{equation}
	with $q \in (0,1]$ \citep{BlatmanJCP2011} \changed{and the quasi-norm $\norme{\ve x}{q} = \left(\sum_{i=1}^d|x_i|^q\right)^{\frac{1}{q}}$}. For $q = 1$, \changed{$\ca^{p,1}$} is the total-degree basis of order $p$. For smaller $q$, this truncation scheme excludes terms with high interaction order while keeping univariate polynomials up to degree $p$.
	\item \textit{Interaction order} 
	The interaction order of the basis can be restricted by defining
	\begin{equation}
	\ca^{p, r} = \{\alp \in \ca^p: \#\{i:\alpha_i \neq 0, i = 1 \enum d\} \leq r\}
	\end{equation}
	\citep{BlatmanCras2008,UQdoc_13_104}.
	This is useful for reducing the number of basis functions especially in high dimensions and when it is known (e.g., for physical reasons) that only a certain number of variables might interact.
\end{itemize}

Instead of using a fixed basis $\ca$, it can be beneficial to employ an iterative scheme which starts from a small set of basis functions (low-dimensional, low-order) and, after computing a sparse solution, repeatedly adapts the basis by including a set of the most promising candidate polynomials and possibly removing others. This is called \textit{basis adaptivity} (not to be confused with Gaussian adaptation \citep{Tipireddy2014}; see also Section \ref{sec:enhancements}).

A simple instance of basis adaptivity is \textit{degree adaptivity} \citep{BlatmanPEM2010}, which is based on total-degree bases. The procedure starts with a basis of low total degree and iteratively increases the total degree of the basis.
Finally, a model selection criterion is used to select the best basis and associated sparse solution.
Similarly, q-norm and interaction order, or a combination of all three, can also be used to design a basis adaptation scheme \citep{BlatmanPEM2010,BlatmanJCP2011}. 
This basis adaptivity is solution-agnostic in the sense that it does not use any information from the solutions computed in previous runs for the augmentation of the basis.
Another solution-agnostic method is the dimension- and order-incrementing algorithm of \citet{Alemazkoor2017}. 
Two methods that adapt the basis based on the active terms of the previous sparse solution are \changed{forward neighbor basis adaptivity}
\changed{\citep{Sargsyan2014, Jakeman2015}} and \changed{anisotropic degree basis adaptivity}
\citep{Hampton2018}.
These approaches keep the size of the basis small by strictly controlling which functions are added to the basis, often starting with a constant surrogate model and adding dimensions only when necessary.
A discussion and benchmark of basis-adaptive methods for sparse PCE is available in \citet{LuethenIJUQ2021}.

\subsection{Experimental design}
\label{sec:ED}

Generally, experimental design techniques aim to select points in order to achieve certain goals related to exploration of the space (space-filling design), or to achieve certain properties of the regression matrix such as (in expectation) orthonormal columns, small determinant, small condition number, etc.
For regression-based PCE, there are several main classes of experimental design techniques:
\begin{itemize}
	\item \textit{Sampling based on the input distribution.} The samples are drawn from the input distribution. Techniques like Latin hypercube sampling (LHS) can be used to improve the space-filling properties.
	\item \textit{Sampling from a different distribution} \changed{(also called \textit{induced sampling} \citep{Guo2020}).} A different distribution and associated basis are constructed that have better properties than the input distribution and its basis.
	\item \textit{Choosing points according to an optimality criterion from a candidate set.} Certain properties of the regression matrix are optimized by choosing the design points from a suitable candidate set.
\end{itemize}
Note that none of these sampling methods considers the evaluations $\ve y$ of the computational model. 
For model-aware sampling techniques (active and supervised learning), see Section \ref{sec:EDenrichment}.

The following sections contain an overview of available sampling methods for sparse PCE. Selected experimental design techniques are described in more detail in Appendix \ref{app:ED}.
Note that due to the large amount of literature on experimental design for sparse PCE, our review cannot be exhaustive. There are many more approaches available, including the deterministic Weil points \citep{Zhou2014}, sparse grids \citep{Perko2014}, randomized or subsampled quadrature points \citep{Berveiller2006, Tang2014, Guo2017}, etc.

\subsubsection{Sampling based on the input distribution}
The most basic sampling method is Monte Carlo (MC) sampling, where the points are sampled independently from the input distribution \citep{Doostan2011,Hampton2015}.
LHS \citep{McKay1979} aims at distributing the design points in a more space-filling way than MC sampling, using a stratification of the input quantile space in each dimension. 
LHS is known to filter main effects; i.e., it reduces the variance of linear regression estimators when the quantity of interest is dominated by terms of interaction order one \citep{Shields2016}. 
LHS with sample decorrelation can further reduce the variance \citep{Owen1994}.
LHS can also be used together with a criterion such as maximin distance (maximize the minimal distance between the design points in quantile space) \citep{Pronzato2012}, where several LHS designs are generated and the one that optimizes the criterion is returned.

A generalization of LHS that combines it with stratified sampling is \textit{Latinized partially stratified sampling} \citep{Shields2016}, which filters both main effects and low-order interaction terms and has been shown to consistently outperform LHS in high-dimensional cases.

Other space-filling/low-discrepancy methods are Sobol' sequences \citep{Sobol1967} and Halton sequences \citep{Halton1960}, which are deterministic but appear to be quasi-random and space-filling in low dimensions. 

\subsubsection{Sampling from a different distribution}
Several methods consider the \textit{coherence} parameter of a basis, defined by
\begin{equation}
\mu(\ca, \{\psi_\alp\}) = \sup_{\ve x \in \cd} \max_{\alp \in \ca} |\psi_\alp(\ve x)|^2,
\label{eq:coherence}
\end{equation}
which can be used to bound the number of samples needed for accurate recovery by $\ell^1$-minimization \citep{Candes2011, Hampton2015}.
For Hermite and Legendre polynomial bases, the coherence parameter of a total-degree basis grows exponentially with the total degree $p$ \citep{Rauhut2012,Yan2012,Hampton2015}.

To construct a \textit{coherence-optimal design}, a new probability distribution and its associated orthonormal basis are constructed that achieve minimal coherence \citep{Hampton2015, Hampton2015b}.
The new basis can be derived from the original PCE basis by multiplying each member by a weight function.
Coherence-optimal samples can be drawn by Markov Chain Monte Carlo (MCMC) \citep{Hampton2015} or by rejection sampling (see section~\ref{app:coh-opt-rejection}).
A related sampling scheme is obtained by constructing a new probability distribution and associated orthonormal basis which have improved but not optimal coherence; however, the distribution is constructed to belong to some classical family and is therefore straightforward to sample. This is called \textit{asymptotic sampling} \citep{Hampton2015, Hampton2015b} and results in a Chebyshev distribution for uniform input, and in a uniform distribution (within a ball of degree-dependent radius) for Gaussian input.
Numerical experiments confirm the expected performance gain of coherence-optimal sampling over both MC and asymptotic sampling, and of asymptotic sampling over MC in the case of low dimension $d$ and high total degree $p$. For high-dimensional problems with low degree, MC often performs better than asymptotic sampling \citep{Hampton2015}.

The so-called \textit{Christoffel sparse approximation (CSA)} \changed{\citep{Jakeman2017,Narayan2017,Cohen2017}} is a related approach which constructs a new orthonormal basis that minimizes the quantity  
\begin{equation}
\tilde\mu(\ca, \{\psi_\alp\}) = \sup_{\ve x \in \cd} \left(\frac{1}{|\ca|} \sum_{\alp \in \ca} |\psi_\alp(\ve x)|^2\right)^\frac12.
\label{eq:CSA_coherence}
\end{equation}
As for coherence-optimal sampling, the new basis can be derived from the original basis by multiplying each member by a weight function, which results in a weighted regression problem.
The corresponding probability distribution is chosen to be the so-called \textit{weighted pluripotential equilibrium measure}, which for bounded distributions is the Chebyshev distribution.
For one-dimensional Gaussian input, this measure is a symmetric Beta distribution with degree-dependent bounds.

Note that all three sampling methods described in this section introduce weights and therefore modify the objective function into a weighted regression problem $ \ve W \ve \Psi \ve c \approx \ve W \ve y$. Since the objective function belongs to the scope of the solver, these methods cannot be considered as pure sampling methods in the sense of being completely independent of the solver.

\subsubsection{Choosing points according to an optimality criterion from a candidate set}
The following methods choose points from a candidate set in order to optimize properties of the regression matrix. Candidate points can be sampled, e.g., using MC, LHS \citep{FajraouiMarelli2017}, coherence-optimal sampling \citep{Diaz2018, Alemazkoor2018}, or Christoffel sampling \citep{Shin2016b}. Note that some of these methods introduce weights, resulting in a weighted regression problem. The candidate set can have a large influence on the resulting design.

\begin{itemize}
	\item \textit{D-optimal sampling} aims at maximizing the determinant 
	$
	D(\ve \Psi) = \det(\frac1N \ve \Psi^T \ve \Psi)^\frac1P
	$
	of the information matrix \citep{Kiefer1959}. Note that $D(\ve \Psi) = 0$ if $\NED < P$.
	Maximizing the D-value is connected to minimizing the variance of the coefficients of the PCE estimate \citep{Zein2013}.
	Algorithms for D-optimal designs include greedy augmentation \citep{Dykstra1971}, exchange techniques \citep{Fedorov2013, Cook1980, Nguyen1992, Zein2013},
	maxvol \citep{Mikhalev2018}, gradient descent \citep{Zankin2018}, and rank-revealing QR decomposition (RRQR)/subset selection \citep{Diaz2018, Gu1996}.
	The advantage of the last method is that it can also be applied for wide matrices $\ve\Psi \in \R^{N \times P} $ where $N < P$.	
	
	\item \textit{S-optimal sampling} (also called ``quasi-optimal'' in \citet{Shin2016a}) selects samples from a large pool of candidate points so that the PCE coefficients computed using the selected set are as close as possible to the coefficients computed from the whole set of candidate points \citep{Shin2016a}.
	The S-value is defined by%
	\footnote{This definition assumes that the columns of the matrix $\ve\Psi_\text{cand}$, containing the evaluations of all candidate points, are mutually orthogonal.}
	\begin{equation}
	S(\ve\Psi) = \left( \frac{\sqrt{\det \ve\Psi^T \ve\Psi}}{\prod_{i=1}^P \norme{\Psi_i}{2}} \right)^{\frac1P}
	\end{equation}
	where $\Psi_i$ denotes the $i$th column of the regression matrix. Its maximization has the effect of maximizing the column orthogonality of the regression matrix while at the same time maximizing the determinant of the information matrix \citep{Shin2016a}. 
	Note that $S(\ve\Psi) = 0$ if $N < P$.
	An S-optimal experimental design can be computed using a greedy exchange algorithm \citep{Shin2016a, Shin2016b, FajraouiMarelli2017}.
	
	\item \textit{Near-optimal sampling} simultaneously minimizes the two matrix properties \textit{mutual coherence} and \textit{average cross-correlation} \citep{Alemazkoor2018}, both of which quantify the correlation between normalized columns of the regression matrix 
	(see section~\ref{sec:near-opt} for the definitions of these properties).
	A near-optimal design can be built by a greedy algorithm \citep{Alemazkoor2018}.
	Note that for near-optimal sampling, it is not necessary that $N \geq P$\changed{, since this method does not rely on the determinant of the information matrix}.
	
\end{itemize}

\subsubsection{Sequential enrichment of the experimental design}
\label{sec:EDenrichment}
Instead of sampling the whole experimental design at once, it has been proposed to use \textit{sequential enrichment}. 
Starting with a small experimental design, additional points are chosen based on the last computed sparse PCE solution or on an augmented basis.
In the context of machine learning, sequential sampling is also known as \textit{active learning} \citep{Settles2012}.
Sequential enrichment has been proposed in the context of S-optimal sampling \citep{FajraouiMarelli2017}, D-optimal sampling \citep{Diaz2018}, and coherence-optimal sampling \citep{Hampton2018}.
\citet{Zhou2019b} suggest an enrichment strategy based on approximations to the expected quadratic loss function, i.e., the mean-squared error. \citet{Ji2008} and \citet{Seeger2008} propose choosing points that minimize the differential entropy of the posterior distribution of the coefficients (using a Bayesian regression setting).
In all cases, numerical examples show that the sequential strategy generally leads to solutions with a smaller validation error compared to nonsequential strategies.
Due to the complexity of the topic and the already large extent of our benchmark, this strategy, albeit promising, is not explored further in this paper.

\subsection{Solution of the minimization problem}
\label{sec:solvers}

There are many formulations of the regression problem that lead to a sparse solution, such as $\ell^0$-minimization, $\ell^1$-minimization (basis pursuit denoising (BPDN), LASSO), $\ell^1-\ell^2$ minimization, Bayesian methods, etc.\ (see also Appendix \ref{app:solvers}).
Based on these formulations, a vast number of sparse solvers has been proposed in the compressed sensing literature; see, e.g., \citet{carronWebsite} and the surveys of \citet{Qaisar2013, Zhang2015c, Arjoune2017}.
We focus here on solvers that have been proposed in the context of sparse PCE.
Of course, it is straightforward to use any other sparse solver to compute a sparse PCE.

The following solvers have been proposed in the sparse PCE literature:
\begin{itemize}
	\item \textit{Convex optimization solvers.} 
	$\ell^1$-minimization in its various formulations is a (constrained) convex optimization problem. 
	Least angle regression (LARS) \citep{Efron2004, BlatmanJCP2011,UQdoc_13_104} is an iterative method that adds regressors one by one according to their correlation with the current residual, and updates the coefficients following a least angle strategy. With the LARS-LASSO modification, which allows for backwards elimination of regressors, LARS is able to generate the whole LASSO path \citep{Efron2004}. Unmodified LARS can also be classified as a greedy method.
	SPGL1 \citep{Vandenberg2008, SPGL1} solves the BPDN formulation by solving a succession of LASSO instances using the spectral projected gradient (SPG) method.
	Other solvers belonging to this class are e.g.\ the solvers implemented in $\ell_1$magic \citep{L1Magic} and SparseLab \citep{SparseLab}.
	
	\item \textit{Greedy methods} are variants of stepwise regression where the regressors are added to the model one by one according to some selection criterion, aiming at finding a heuristic solution to the intractable $\ell^0$-minimization formulation. 
	Orthogonal matching pursuit (OMP) \citep{Tropp2007,Doostan2011,UQdoc_13_104} is a classical forward selection algorithm in which orthonormalized regressors are added to the model one by one according to their correlation with the residual, and the coefficients are computed by least-squares. \citet{Baptista2019} suggests extensions to OMP such as parallelization, randomization and a modified regressor selection procedure.
	Subspace pursuit (SP) \citep{Dai2009, Diaz2018, SPcode} is an iterative algorithm that repeatedly uses least squares on a subset of regressors. 
	LARS \citep{Efron2004, BlatmanJCP2011} without the LASSO modification (allowing for removal of regressors) can also be classified as a greedy method.
	Another greedy method is ranking-based sparse PCE \citep{Tarakanov2019} which employs batch updating, coordinatewise gradient descent of the elastic net formulation, and a correlation- and stability-based ranking procedure for the regressors.
	Many more greedy stepwise regression techniques have been proposed, utilizing various selection criteria, solvers, and stopping criteria. An overview of methods following this scheme is given in section~\ref{app:greedystepwise}.
	
	\item \textit{Bayesian compressive sensing (BCS)} (a.k.a.\ sparse Bayesian learning) is a class of methods that use a Bayesian setting to find a sparse solution. They impose a sparsity-inducing prior on the coefficients, whose parameters are again considered to be random variables with a hyperprior \changed{\citep{Tipping2001,Ji2008,Sargsyan2014,Tsilifis2020}}. The solution is typically the maximum a posteriori estimate of the coefficients and can be computed e.g.\ by differentiation \citep{Tipping2001}, expectation-maximization \citep{Figueiredo2003, Wipf2004b}, expectation-propagation \citep{Seeger2008}\changed{, variational inference \citep{Tsilifis2020, Bhattacharyya2020}}, or a fast approximate algorithm \citep{Faul2002, Tipping2003}.
	An extension called FastLaplace with an additional layer of hyperparameters has been proven to attain even sparser solutions \citep{Babacan2010, FastLaplace}. 
	A greedy algorithm using the Bayesian setting to select the regressors is the greedy Bayesian Kashyap information criterion (KIC)-based algorithm \citep{Shao2017}.
	
	\item \textit{Iteratively reweighted methods.}
	Iteratively reweighted $\ell^1$-minimization uses the coefficients computed in a previous iteration to 
	construct a weighted $\ell^1$-minimization problem \citep{Candes2008b, Yang2013}.
	\citet{Cheng2018b} suggest an iterative reweighted method with D-MORPH regression \citep{Li2010} as its computational core, which is a technique that follows a certain path, defined by a quadratic objective function, on the manifold of solutions to the underdetermined system. 
		
\end{itemize}

Each of the solvers mentioned above features one or more hyperparameters whose values must be calibrated. This is usually done by cross-validation. Popular choices are leave-one-out (LOO) cross-validation (accelerated for least-squares solutions) \citep{BlatmanPEM2010, BlatmanJCP2011}, LOO cross-validation with a modification factor for small sample sizes \citep{Chapelle2002, BlatmanJCP2011}, and $k$-fold cross-validation 
\citep{Doostan2011,Jakeman2015,Huan2018}.

Selected sparse regression solvers are described in more detail in Appendix \ref{app:solvers}.

\subsection{Model selection criterion}
\label{sec:modelselection}
\label{sec:crossvalidation}

To decide whether to continue iterating in the framework or stop the process, we need to assess how well the current sparse solution performs.
Our main quantity of interest is the \textit{generalization error}, which quantifies the mean-square accuracy of the surrogate.
It is given by
\begin{align}
E_\text{gen} = \Espe{\ve X}{\left( \cm(\ve X) - \cm^\text{PCE}(\ve X) \right)^2} 
\label{eq:generalizationerror}
\end{align}
where $\cm$ is the computational model, $\ve X$ is the random input vector, and $\cm^\text{PCE}$ is the sparse PC surrogate.

The generalization error can be approximated by the \textit{validation error}, which is the MC estimate of \eqref{eq:generalizationerror} on a validation set $\{(\ve x\ival, y\ival): \ve x\ival \sim_\text{i.i.d.} f_{\ve X}, y\ival = \cm(\ve x\ival), i=1\enum N_\text{val}\}$.
To make the validation error independent of the scaling of the model, it is convenient to use the \textit{relative mean-squared error} defined by 
\begin{align}
\text{RelMSE} = \frac{\sum_{i=1}^{\Nval} (y\ival - \cm^\text{PCE}(\ve x\ival))^2}{\sum_{i=1}^{\Nval} (y\ival - \bar{y})^2}
\label{eq:RelMSE}
\end{align}
where $\bar{y} = \frac{1}{\Nval} \sum_{i=1}^{\Nval} y\ival$.

The best surrogate model $\cm^\text{PCE}$ (defined by $\ca$ and $\ve c$) is the one that has the smallest generalization error. 
In practical applications, the generalization error typically cannot be computed, and a large validation set is not available due to computational constraints. 
Instead, we define a \textit{model selection criterion} that acts as a proxy for the generalization error. 
A typical stopping criterion in the PCE framework of Figure \ref{fig:framework} is the observation that the model selection criterion no longer improves.
The following model selection criteria have been proposed in the sparse PCE literature:
\begin{itemize}
	\item \textit{$k$-fold cross-validation (CV)} \citep{Hastie:2001, Jakeman2015, Hampton2018}, which approximates the validation error by building a surrogate several times on different subsets of the data, and evaluating the error on the remaining data points. 
		
	\item \textit{Leave-one-out (LOO) cross-validation} \citep{Hastie:2001,BlatmanPEM2010, BlatmanJCP2011}, which is $\NED$-fold cross-validation (where $\NED$ is the size of the experimental design). 
	For PCE approximations computed by OLS, there exists an efficient formula to evaluate the LOO error \cite[Appendix D]{BlatmanJCP2011}.
	
	\item \textit{Modified LOO} \citep{BlatmanJCP2011}, which uses a correction factor for the LOO which was derived for the empirical error for OLS with small sample size \citep{Chapelle2002}. The correction factor depends on the experimental design and the active basis functions.
	
	\item \textit{Kashyap information criterion (KIC)} \citep{Shao2017,Zhou2019c}, an approximation to the Bayesian model evidence, which is the likelihood of observations given the model.
	
	\item \textit{Sparsity} \citep{Alemazkoor2017}, which uses the idea that a larger basis should lead to a sparser solution when the necessary basis functions enter the model, unless the ratio of basis functions to model evaluations becomes too large.
	
\end{itemize}

A model selection criterion is also often used to determine the hyperparameter of the sparse solver (see Section \ref{sec:solvers}).

If cross-validation is used to select the solver hyperparameter, this estimate of the validation error is often too optimistic due to model selection bias. Instead of reusing this estimate for model selection, it is better to perform an outer loop of $k$-fold or LOO cross-validation, a procedure called \textit{double cross-validation} or cross-model validation \citep{Baumann2014, LiuWiart2020b}.

\subsection{Further enhancements of sparse PCE}
\label{sec:enhancements}
There are many enhancements to the simple scheme for sparse PCE presented in Figure \ref{fig:framework}.
The following methods have been suggested to improve the accuracy of the solution and reduce the number of model evaluations needed:
\begin{itemize}
	\item \citet{Alemazkoor2018b} construct a preconditioning matrix for a given regression matrix which reduces the mutual coherence while avoiding deterioration of the signal-to-noise ratio.
	\item \citet{Huan2018} suggest a technique called stop-sampling, which guides the decision of whether to obtain more samples (sequential ED enrichment) by observing the decrease of the CV error.
	\item In the case when \changed{prior} information about the magnitude of the coefficients is available, \citet{Peng2014} use this information to construct a weighted regression problem which allows a more accurate solution with fewer points (similar to iteratively reweighted $\ell^1$-minimization).
	\item \citet{LiuWiart2020a} use resampled PCE, which is a technique for improving the PCE solution by aggregating the results of several solver runs on different subsets of the data. Only the terms that are chosen most often by the solvers are retained in the final solution. 
		
	\item Several methods exist to reduce the dimension of the input space before computing the sparse PCE.
	Unsupervised methods are principal component analysis (PCA) and kernel PCA \citep{Lataniotis2019}.
	``Basis adaptation'' methods (referring to a basis of the input random space) determine a suitable rotation of the input space, often assumed to be independent standard Gaussian, into new coordinates which permit a sparser representation in fewer coordinates (see \citet{Tipireddy2014,Yang2018,Tsilifis2019} and others).
	A related technique is nonlinear PCE-driven partial least squares (PLS) \citep{Papaioannou2019, Zhou2020}, which reduces the input dimension by identifying directions in the input space that are able to explain the output well in terms of a sum of one-dimensional PCEs.

\end{itemize}

\section{Numerical results}
\label{sec:benchmark}

\subsection{Benchmark design}
While the number of methods for computing sparse PCE is large, to the best of the authors' knowledge there is no comprehensive benchmark study on this topic. Most publications only compare the newly developed method to one or two baseline methods.
An overview of publications containing comparisons of sparse PCE methods is presented in Appendix~\ref{app:benchmarkstudies}.

Since the number of possible combinations of sampling schemes, sparse regression solvers, basis adaptation schemes, model evaluation criteria etc.\ is huge (see Section \ref{sec:framework} and thereafter), we restrict our benchmark as follows to some of the most promising and best-known methods: 
\begin{itemize}
	\item We consider the sampling schemes MC, LHS, coherence-optimal, and D-optimal.
	LHS is used together with a maximin criterion to improve the space-filling property (using the MATLAB function \texttt{lhsdesign}).
	D-optimal designs are constructed from a coherence-optimal candidate set%
	\footnote{We have also conducted all benchmark experiments with D-optimal designs constructed from LHS candidate sets, but we do not display these results, because we found that in most cases, D-opt(LHS) sampling performs (significantly) worse than most other sampling schemes, and often worse than its candidate set LHS. This matches with the results of \citet{FajraouiMarelli2017} who observed this in a sequential enrichment setting and with the LARS solver.}
	using the subset selection/RRQR algorithm, which allows for the construction of D-optimal experimental designs with size $N$ smaller than the number of regressors $P$ \citep{Diaz2018}. There is no such algorithm for S-optimal sampling, which is why we do not consider the latter in this benchmark.
	We do not consider Sobol' sequences, since they have been shown to be outperformed by LHS in sparse PCE applications \citep{FajraouiMarelli2017}.
	Near-optimal sampling can realistically be used only for rather small bases ($P \in \co(100)$), since its algorithm scales as \changed{$\co(MP^2)$, with $M = 10P$ as suggested by \citep{Diaz2018} (see below)}. We use it with a coherence-optimal candidate set for two models with small basis.
	\item We consider the sparse regression solvers LARS, OMP, subspace pursuit (SP), FastLaplace (which we call here BCS), and SPGL1.
	Each of these solvers involves at least one hyperparameter, whose range is chosen according to reasonable guesses. For LARS and OMP, the hyperparameter is the number $K$ of selected regressors and its range is $[1, \min\{P,N-1\}]$. For SP, $K$ must fulfill $2K \leq \min\{P, N\}$. For BCS and SPGL1, the hyperparameter $\sigma$ is chosen from the range $\sigma^2 \in  N \changedmath{\cdot} \Varhat{\ve y} \changedmath{\cdot} [10^{-16}, 10^{-1}]$ which resembles a suitable range of possible relative MSE values. The hyperparameter values of LARS and OMP are determined by modified LOO cross-validation, while the hyperparameters of SP, BCS, and SPGL1 are determined by $k$-fold cross-validation (Section \ref{sec:solvers}).
	\changed{In addition, we consider a variant of SP which uses LOO cross-validation instead of $k$-fold cross-validation, which we name \emph{\SPloo}.}
	\item We only consider the nonadaptive setting, in which both the basis and the size of the experimental design are fixed before the sparse PCE is computed.
	\item For each model, we define a reasonable range of \changed{5--7} experimental design sizes. Each experiment is repeated 30--50 times to account for statistical uncertainty. The experimental designs are generated anew for each repetition and each ED size. All solvers are tested on the same ED realizations.
	\item The \changed{coherence-optimal} candidate sets \changed{from which the D-optimal designs are selected} have size $M = 10P$ as in \citep{Diaz2018}. For computational reasons, \changed{they} 
	are not sampled completely anew for each replication, but are rather drawn uniformly at random without replacement from a larger set of size $2M = 20P$ as in \citep{Diaz2018}. 
	\item Since we are interested in sparse PCE for the purpose of surrogate modelling, our main quantity of interest is the relative mean-squared error%
	\footnote{Note that some authors such as \citet{Doostan2011,Hampton2015,Diaz2018, Alemazkoor2018} choose to normalize instead by $\sum_{\ve x \in \cx_\text{val}} \cm(\ve x)^2$  or use the unnormalized mean-squared error \citep{Shin2016b}. To assess the recovery of sparse vectors just as in compressed sensing, some consider the error in the coefficient vector instead of the error in the model approximation \citep{Alemazkoor2018}.}% 
	\footnote{Since a typical application of PCE is the computation of moments and Sobol' indices, the error in these quantities is another possible performance measure. However, globally accurate prediction as considered in this paper is more challenging that the prediction of moments and Sobol' indices, which are accurate if the largest-in-magnitude coefficients are estimated accurately. If a globally accurate surrogate model can be constructed, typically also the moments and Sobol' indices are accurate.}
	(RelMSE) as defined in \eqref{eq:RelMSE}.
	We investigate the RelMSE for several models, sparse solvers, and experimental design techniques.
	Typically, the practical interest lies in small experimental designs. 
	\item Since the experimental design is random, the resulting validation error is a random variable. We visualize the data with boxplots. When comparing the performance of different methods, we consider the median performance and the spread of the resulting validation error. However, often there can be considerable overlap of validation errors between methods.
\end{itemize}

\subsection{Software}
For the implementation of the benchmark, we use the general-purpose uncertainty quantification software UQLab \citep{MarelliUQLab2014}. 
UQLab supports the integration of other software packages.%
\footnote{A description of how to use custom sparse solvers and sampling schemes in the UQLab framework can be found in the supplementary material.}
We utilize the following code:
\begin{itemize}
	\item UQLab for MC sampling and LHS \citep{MarelliUQLab2014}.
	\item \texttt{DOPT\_PCE} for D-optimal sampling (subset selection/RRQR) and subspace pursuit \citep{Diaz2018, SPcode}.
	\item An in-house developed rejection-based implementation of coherence-optimal sampling.
	\item An in-house implementation of near-optimal sampling based on the description by \citet{Alemazkoor2018}.
	\item UQLab for the solvers LARS and OMP \citep{MarelliUQLab2014}.
	\item \texttt{spgl1-1.9} for SPGL1 \citep{Vandenberg2008, SPGL1}.
	\item \texttt{FastLaplace} for the hierarchical implementation ``FastLaplace'' of BCS \citep{Babacan2010, FastLaplace}.
\end{itemize}

\subsection{Benchmark: Considered models}
\label{sec:benchmark_models}
\changed{Our benchmark is performed on a selection of 11 computational models of varying complexity and input dimensionality, which are typical benchmark models in the context of sensitivity and reliability analysis.
An overview of these models is given in Table~\ref{table:models}. For details on the models, we refer the reader to the respective publications supplied in the last column of the table. 
While of course not representative of all possible classes of engineering models, we believe that this sample provides a good testing ground for the comparative performance among different approaches for computing sparse PCE. }

\begin{table}[htbp]
	\footnotesize
	\centering
	\caption{\changed{Overview of the 11 computational models used in our benchmark. Italic font denotes finite element (FE) models , all other models are analytical. 
			For each model, a static total-degree basis with hyperbolic truncation defined by $p$ and $q$ is used. The values are chosen to fulfill $P \approx \frac{10}{3} N_\text{max}$, where $N_\text{max}$ is the largest tested experimental design size. 
			The values for $p$ in parentheses for the Ishigami and borehole models refer to the smaller basis used in Section~\ref{sec:results_nearopt}.
			The column ``Reference'' provides the relevant literature in which the models and their probabilistic inputs are described in detail. }}
	\label{table:models}
	\renewcommand{\arraystretch}{1.2}
	\begin{tabular}{l c >{\centering}p{.22\textwidth} >{\centering}p{.11\textwidth} >{\centering}p{.05\textwidth} >{\centering\arraybackslash}p{.11\textwidth}} 
		\hline
		Model & Dimension & Input distributions & Basis & $N_\text{max}$ & Reference \\
		\hline
		Ishigami function & 3 & uniform & $p = 14$ ($12$), \newline $q=1$& $200$&  \citep{BlatmanJCP2011}  \\
		Undamped oscillator & 6 & Gaussian & $p = 5$, \newline $q=1$& $150$ &\citep{Echard2013}  \\
		Borehole function & 8 & Gaussian, lognormal, uniform & $p = 5$ ($4$), \newline $q=1$& $300$ & \citep{Harper1983}  \\ %\citep{Morris1993}
		Damped oscillator & 8 & lognormal & $p = 5$, \newline $q=1$ & $400$ &\citep{DubourgThesis} \\
		Wingweight function & 10 & uniform  & $p = 4$, \newline $q=1$& $300$ &\citep{Forrester2008} \\
		\textit{Truss model} & 10 & lognormal, Gumbel  & $p = 4$, \newline $q=1$& $300$ & \citep{BlatmanJCP2011} \\
		%		\hline
		Morris function & 20 & uniform & $p = 8$, \newline $q=0.5$& $400$ & \citep{BlatmanRESS2010}  \\ 
		%\citep{Morris1991}(Blatman\&Sudret version) \\
		\textit{Structural frame model} & 21 & lognormal, Gaussian; {dependent input variables} & $p = 8$, \newline $q=0.5$& $400$ & \citep{BlatmanPEM2010}  \\
		\textit{$2$-dim diffusion model} & 53 & Gaussian & $p = 4$, \newline $q=0.5$& $500$ & \citep{KonakliRESS2016}  \\
		\textit{$1$-dim diffusion model} & 62 & Gaussian & $p = 4$, \newline $q=0.5$& $500$ & \citep{FajraouiMarelli2017}  \\
		\changed{$100$D} function & 100 & uniform & $p = 4$, \newline $q=0.5$& $1400$ & UQLab \newline example\footnotemark  \\
		\hline	
	\end{tabular}
\end{table}

\changed{In addition to analyzing aggregated performance on all 11 models, we investigate the behavior of the methods in detail on a subset of four \textit{spotlight models}, each of which possesses characteristic properties that might influence the approximation quality of sparse PCE methods:
the \textit{Ishigami function} is low-dimensional and highly compressible in the PCE basis but requires a high-degree basis to be approximated accurately. The \textit{borehole function} is smooth and nonlinear and therefore is an example for a well-behaved engineering model. A \textit{two-dimensional diffusion model}, a stochastic heat diffusion PDE in two physical dimensions, is high-dimensional, not analytical, and the magnitude of its expansion coefficients decays only slowly. Finally, the \textit{$100$D function} is high-dimensional, analytical, and compressible.}

The Ishigami model is the well-known three-dimensional, highly nonlinear, smooth analytical function
\begin{equation}
f(X_1, X_2, X_3) = \sin(X_1) + a \sin^2(X_2) + b X_3^4 \sin(X_1)
\end{equation}
taking uniform input $\ve X \sim \cu([-\pi, \pi]^3)$. A typical choice is $a = 7, b = 0.1$. For this function, any sparse solver should be able to find a sparse solution.

The borehole function simulates the water flow through a borehole between two aquifers \citep{Harper1983}. It is an eight-dimensional nonlinear function which, despite having an analytical form, is not trivial to approximate. It is defined by
\begin{equation}
B(r_w, L, K_w, T_u, T_l, H_u, H_l, r) = \frac{2\pi T_u (H_u - H_l)}{\ln\left(r / r_w\right) \left( 1 + \frac{2 L T_u}{\ln\left(r / r_w\right) r_w^2 K_w} + \frac{T_u}{T_l}\right)}.
\end{equation}
Its input random variables and their distributions are provided in Table \ref{table:borehole}.

\begin{table}
	\caption{Borehole function: Input random variables and their distributions}
	\label{table:borehole}
	\centering
	\footnotesize
	\begin{tabular}{l l l}
		\hline
		Variable & Distribution & Description \\
		\hline
		$r_w$ & $\cn(0.10, 0.0161812)$ & borehole radius \\
		$L$ & $\cu([1120, 1680])$ & borehole length \\
		$K_w$ & $\cu([9855, 12045])$ & borehole hydraulic conductivity  \\
		$T_u$ & $\cu([63070, 115600])$ & transmissivity of upper aquifer  \\
		$T_l$ & $\cu([63.1, 116])$ & transmissivity of lower aquifer  \\
		$H_u$ & $\cu([990, 1110])$ & potentiometric head of upper aquifer  \\
		$H_l$ & $\cu([700, 820])$ &  potentiometric head of lower aquifer \\
		$r$ & $ \text{Lognormal}([7.71, 1.0056])$ & radius of influence  \\
		\hline
	\end{tabular}
\end{table}

The two-dimensional heat diffusion model \citep{KonakliRESS2016} is defined by the partial differential equation (PDE)
\begin{equation}
- \nabla \cdot (\kappa(\ve x) \nabla T(\ve x)) = Q \mathbb{1}_A(\ve x) \qquad \text{ in }\Omega = [-0.5, 0.5]^2
\end{equation}
with boundary conditions $T = 0$ on the top boundary and $\nabla T \cdot \ve n = 0$ on the left, lower, and right boundaries of the square domain $\Omega$, where $\ve n$ denotes the outer unit normal (see \citep{KonakliRESS2016} for an illustration of the setup). 
Here, the source is in $A = [0.2, 0.3]^2$ with strength $Q = 500$.
The output quantity of interest is the average temperature in $B = [-0.3, -0.2]^2$. 
The diffusion coefficient $\kappa(\ve x)$ is modelled by a lognormal random field with mean $\mu_\kappa = 1$ and standard deviation $\sigma_\kappa = 0.3$. The autocorrelation function of the underlying Gaussian random field is an isotropic squared-exponential with length scale $l = 0.2$.
The random field $\kappa(\ve x)$ is discretized using the EOLE method \citep{DerKiureghian1993} with $d = 53$ terms, which comprises $99\%$ of its variance. 
The solution to an individual heat diffusion problem is computed using an in-house finite \changed{element} code \citep{KonakliRESS2016}.
The input comprises $d=53$ independent standard normal random variables.

Finally, the so-called \changed{100D} function is an analytical model of the form 
\begin{align}
	f(\ve X) &= 3 - \frac{5}{d} \sum_{i = 1}^{d} i X_i +
	\frac{1}{d}\sum_{i = 1}^{d} i X_i^3  + 
	\frac{1}{3d}\sum_{i = 1}^{d} i  \ln\left(X_i^2 + X_i^4\right) \\
	& \quad + X_1 X_2^2  + X_2 X_4 - X_3 X_5 + X_{51} + X_{50}X_{54}^2
	\nonumber
\end{align}
taking uniform inputs $X_i \sim \cu([1,2]), i \neq 20$, and $X_{20} \sim \cu([1,3])$.
We use $d = 100$.
This function was designed for sensitivity analysis: the first-order sensitivity indices of the input variables are generally nonlinearly increasing with their index, with certain variables having especially high sensitivity. The model also contains four interaction terms. 
It is an example from UQLab.\textsuperscript{\ref{fn:uqlab}}
\footnotetext{\url{https://www.uqlab.com/sensitivity-high-dimension}\label{fn:uqlab}}

For each of the models, we use a fixed basis for the benchmark. In general, the best total degree $p$ and the hyperbolic truncation $q$ are a priori unknown. We heuristically choose $q = 1$ for low-dimensional models ($d \leq 10$) and $q = 0.5$ for high-dimensional models (\changedmath{$d \geq 20$}). The degree $p$ is chosen so that the number of basis functions $P$ is approximately $\frac{10}{3}N_\text{max}$, where $N_\text{max}$ is the largest number of experimental design points for the specific benchmark. This choice is based on the reasoning that for an experimental design of size $N$, sparse solvers like LARS often select an active basis of size $ \approx \frac{N}{3}$, and that the candidate basis might be 10 times larger than the final active basis to be sufficiently rich.
We focus on rather small experimental designs, since our goal is not to investigate the convergence of the methods as $N \to \infty$ (which has been demonstrated elsewhere), but to decide which methods are most efficient for small $N$.
This results in 
\changed{the choice of values for $p$, $q$, and $N_\text{max}$ displayed in Table~\ref{table:models}.}

\subsection{Results: Comparison of solvers}
\label{sec:results_only_solvers}
First, we use a fixed sampling scheme (LHS) to compare the performance of the \changed{six} solvers LARS, OMP, subspace pursuit (SP) using $k$-fold cross-validation, \changed{Subspace Pursuit using LOO cross-validation (\SPloo{})}, FastLaplace (BCS), and SPGL1 on \changed{all 11} benchmark models described above.

\begin{figure}[htbp]
	\centering
	\subfloat[][Ishigami function]
	{\includegraphics[width=.49\textwidth]{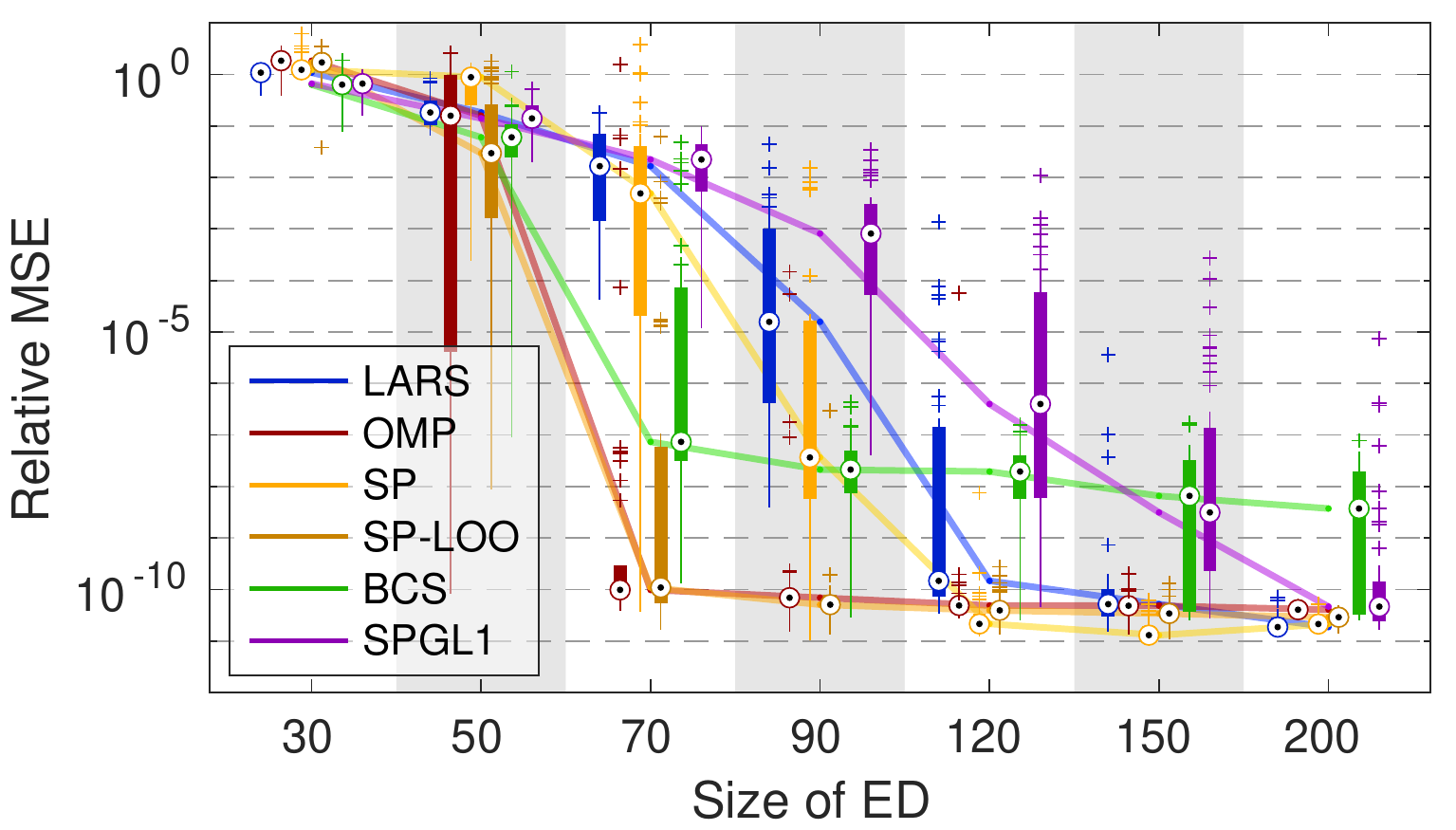}\label{fig:onlysolvers_ishigami}}
	\hfill
	\subfloat[][Borehole model]
	{\includegraphics[width=.49\textwidth]{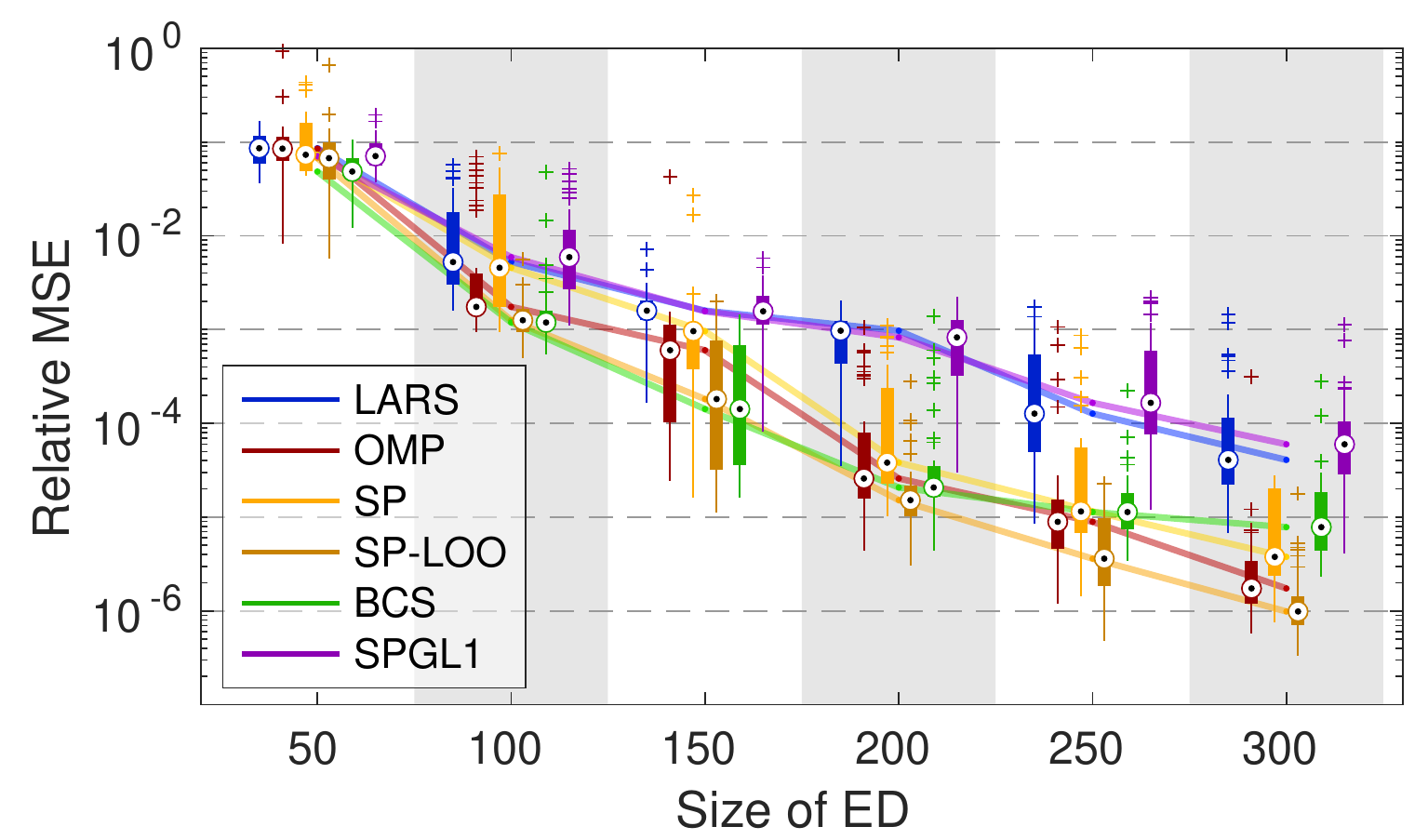}\label{fig:onlysolvers_borehole}}
	\\
	\subfloat[][Two-dimensional diffusion]
	{\includegraphics[width=.49\textwidth]{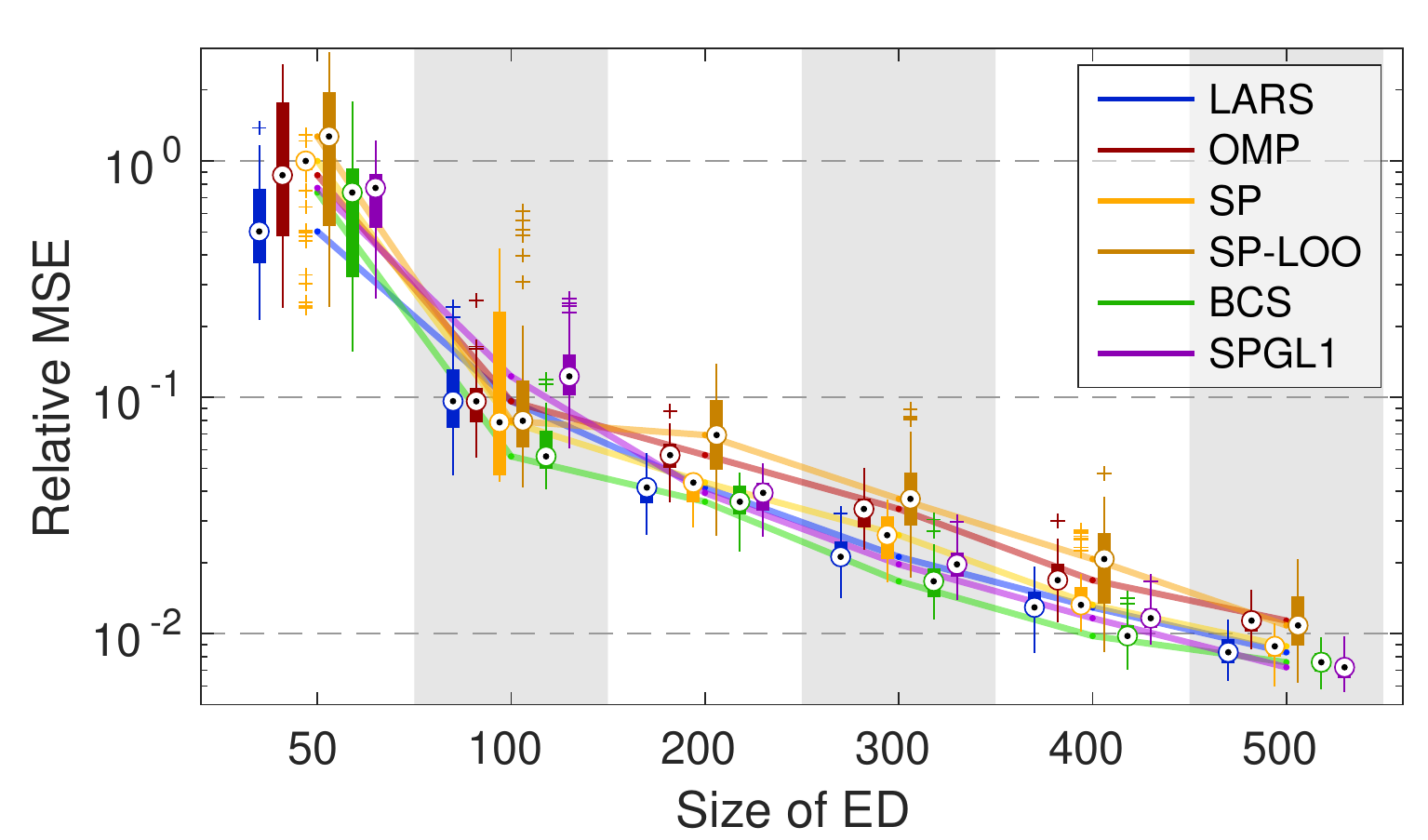}\label{fig:onlysolvers_diffusion2D}}
	\hfill
	\subfloat[][\changed{100D} function]
	{\includegraphics[width=.49\textwidth]{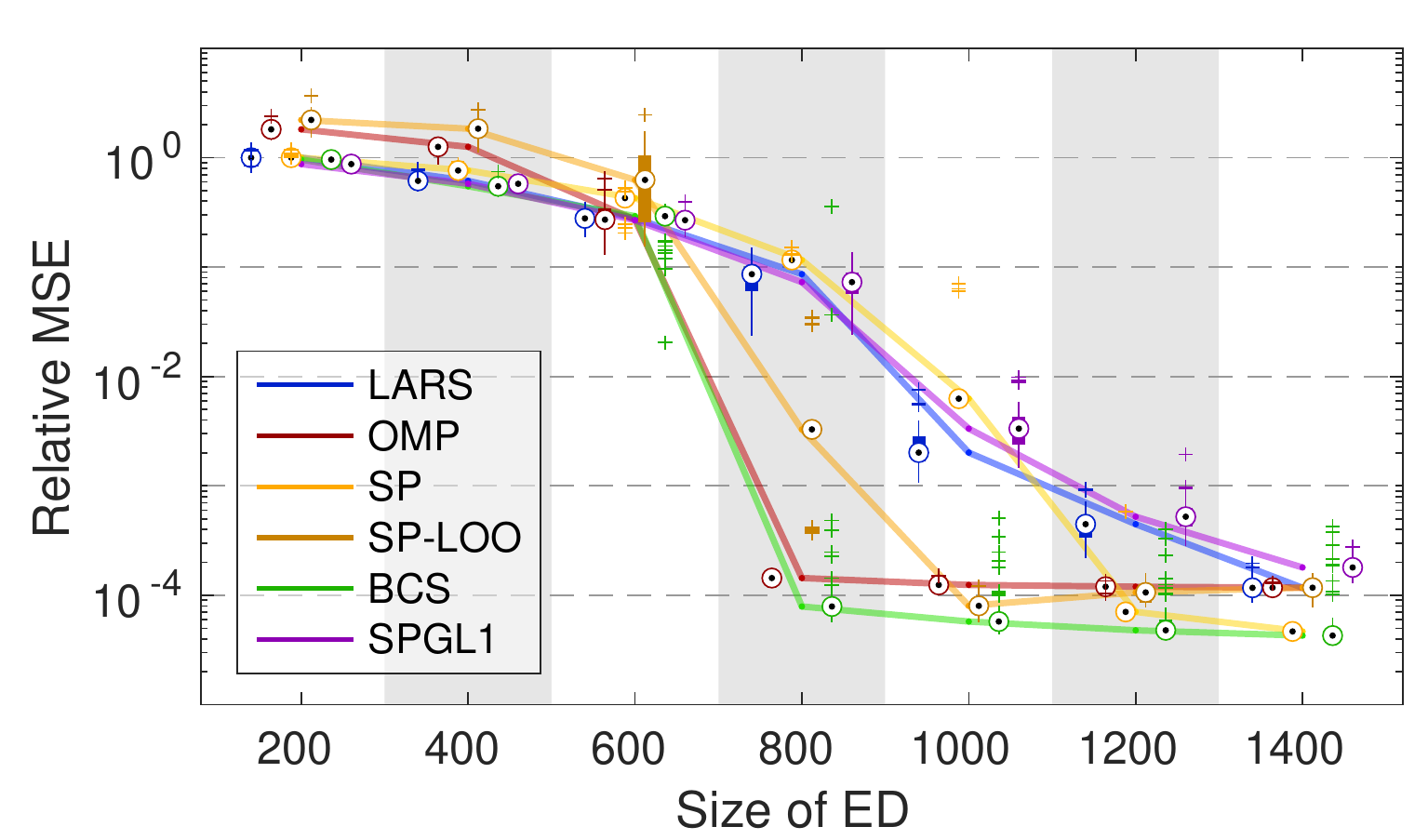}\label{fig:onlysolvers_highdimfct}}
	\caption{\changed{Boxplots of relative MSE against ED size from the solver benchmark for the four spotlight models. Results for six sparse solvers and LHS design. For the corresponding plots for the seven remaining models, see Figure~\ref{fig:results_more_models} in Appendix~\ref{app:additional_results}.}}
	\label{fig:results_onlysolvers}	
\end{figure}

In \changed{Figure~\ref{fig:results_onlysolvers}} we display boxplots (50 replications) of relative MSE against experimental design size for all \changed{six} solvers \changed{for the four spotlight models}.
\changed{For the remaining seven models, the corresponding boxplots (30 replications) of relative MSE against experimental design size are provided in Figure~\ref{fig:results_more_models} in Appendix~\ref{app:additional_results}.}
In the plots, the lines as well as the dot inside the white circle denote the median of the relative MSE.
We make the following observations:
\begin{itemize}
	\item For the smallest experimental designs, all solvers perform similarly \changed{poorly: there is not enough information in the ED to construct an accurate surrogate model}. For larger experimental designs, there can be considerable differences between the solvers' generalization errors of up to several orders of magnitude\changed{, which demonstrates that the solvers do not use the available information in identical ways}.
	\item BCS, \changed{\SPloo{},} and OMP are often among the best solvers. BCS performs especially well for smaller experimental design sizes.
	In the case of the Ishigami model, it seems to plateau earlier than the other solvers. It also tends to find sparser solutions than the other solvers (not shown in plots)\changed{, which might explain these observations: sparse solutions are advantageous when only limited data is available, but at the same time they carry the risk of ignoring important terms}. For \changed{large ED sizes and the two highly compressible models (Ishigami and \changed{100D} function), BCS} has a larger spread than the other solvers\changed{, possibly because the sparsity-enforcing procedure does not always include all of the important terms}. \changed{In contrast, the greedy solver OMP returns rather dense solutions, which seem to generalize well.}
	\changed{\SPloo{} performs well in general for low-dimensional models.}
	\item SP does not perform well for small ED sizes, but for large ED sizes it \changed{sometimes} outperforms the other solvers. Together with BCS, it tends to find sparser solutions than the other solvers \changed{(not shown in plots)}.
	\item LARS and SPGL1 \changed{often achieve a similar generalization error, which is often larger, sometimes significantly, than that of the other solvers.} SPGL1 tends to return rather dense solutions (not shown in plots)\changed{, which might not generalize as well as other solutions}.
	\item \changed{Some models characterized by relatively poor compressibility (e.g. diffusion and frame) show comparable performance among all solvers. This is expected, as in such cases the sparsity assumption is a rather weak proxy for solution quality.}
	\item \changed{The exceptions to the general observations outlined above are the damped oscillator and the Morris function, for which the solver performance is reversed, with LARS and SPGL1 among the best solvers, and OMP and \SPloo{} among the worst.
		In these cases, however, none of the methods achieves satisfactory accuracy within the available computational budget.}
\end{itemize}

\begin{figure}[htbp]
		\centering
		\subfloat[][low-dim models ($d \leq 10$), small ED sizes]
		{\includegraphics[width=.56\textwidth]{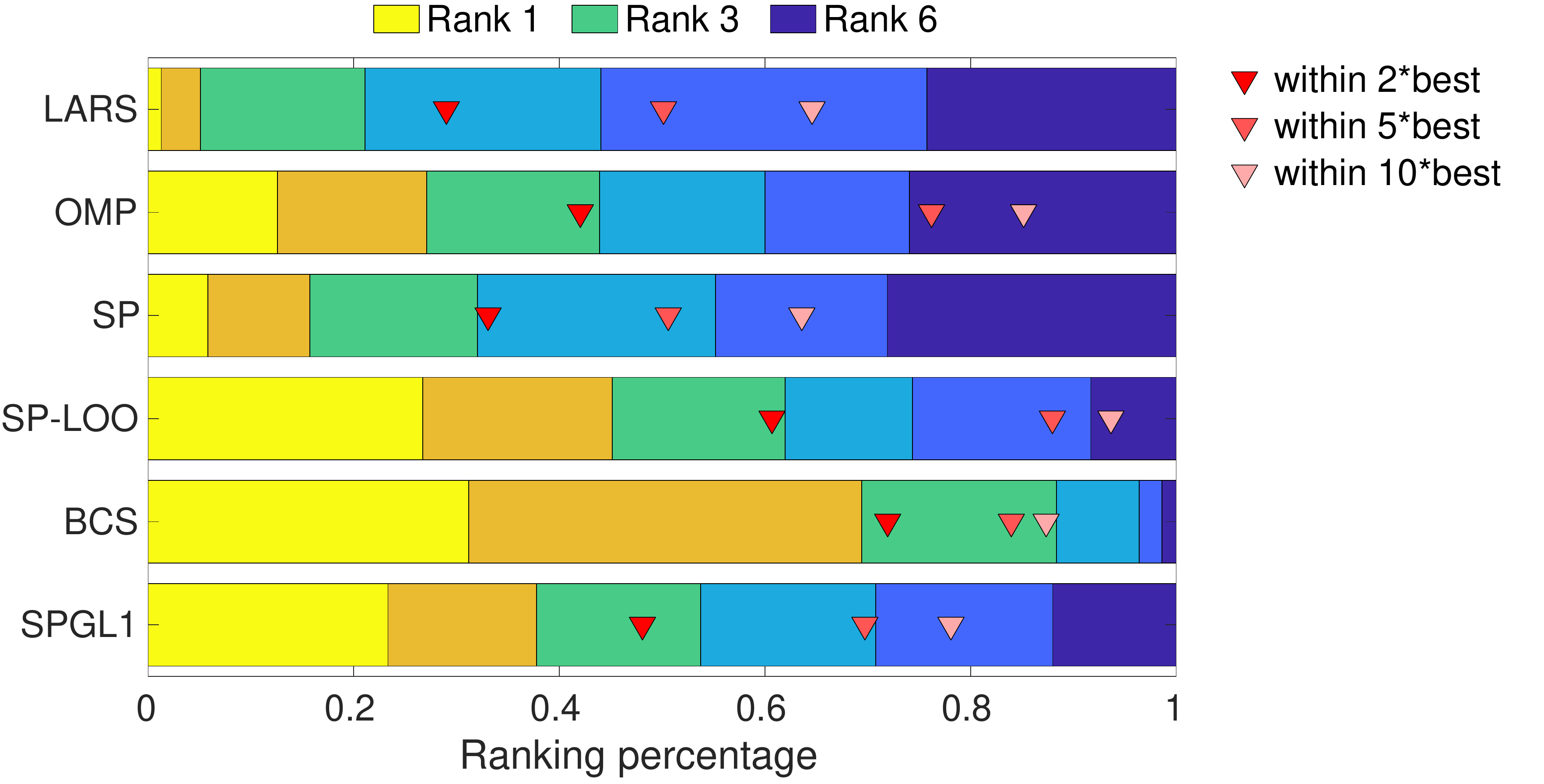}\label{fig:results_more_models_aggregated1}}
		\hfill
		\subfloat[][low-dim models ($d \leq 10$), large ED sizes]
		{\includegraphics[width=.435\textwidth, trim = 0 0 8.2cm 0, clip]{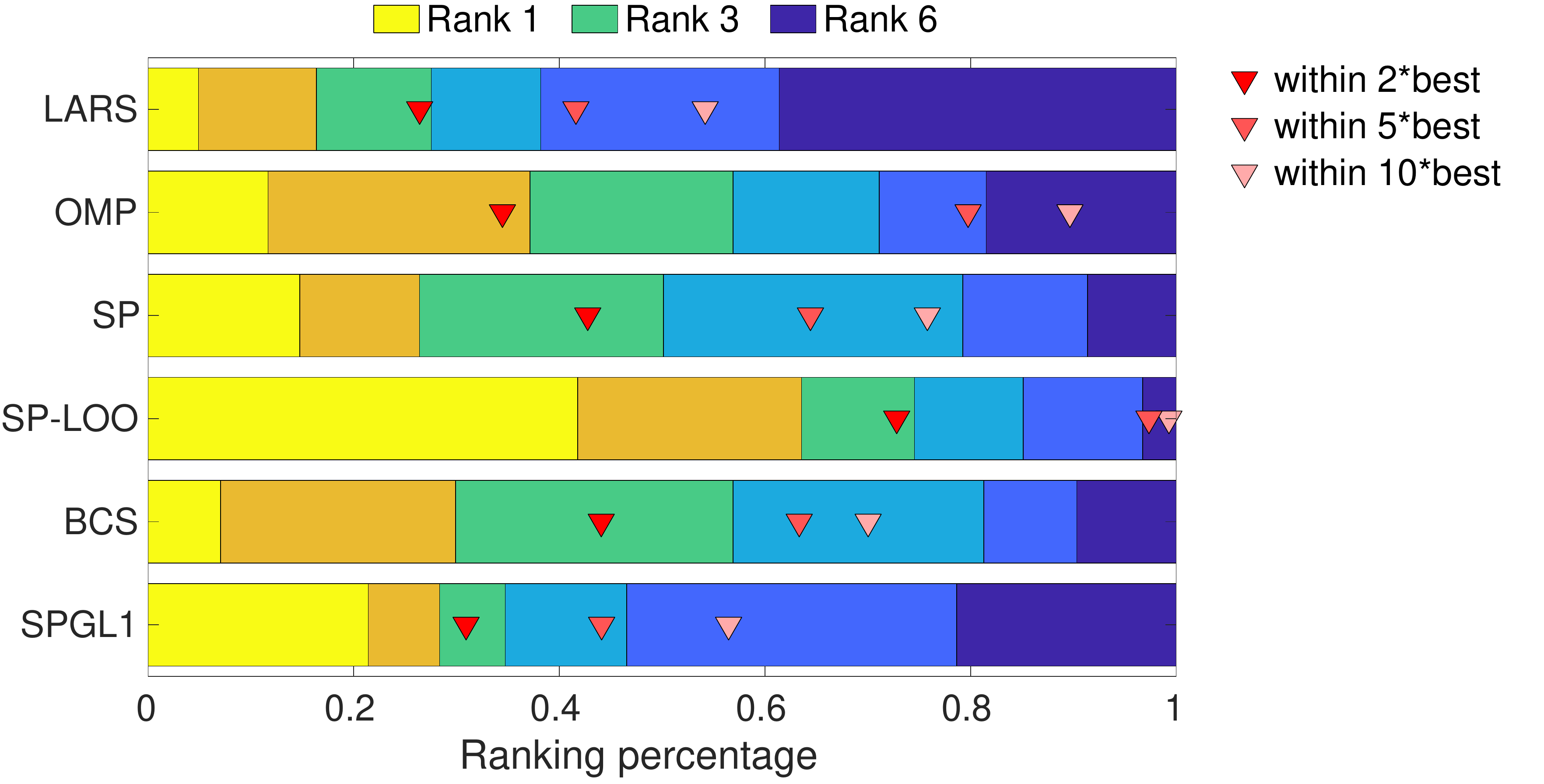}\label{fig:results_only_solvers_lowdim_largeED}}
		\\
		\subfloat[][high-dim models ($d \geq 20$), small ED sizes]
		{\includegraphics[width=.56\textwidth]{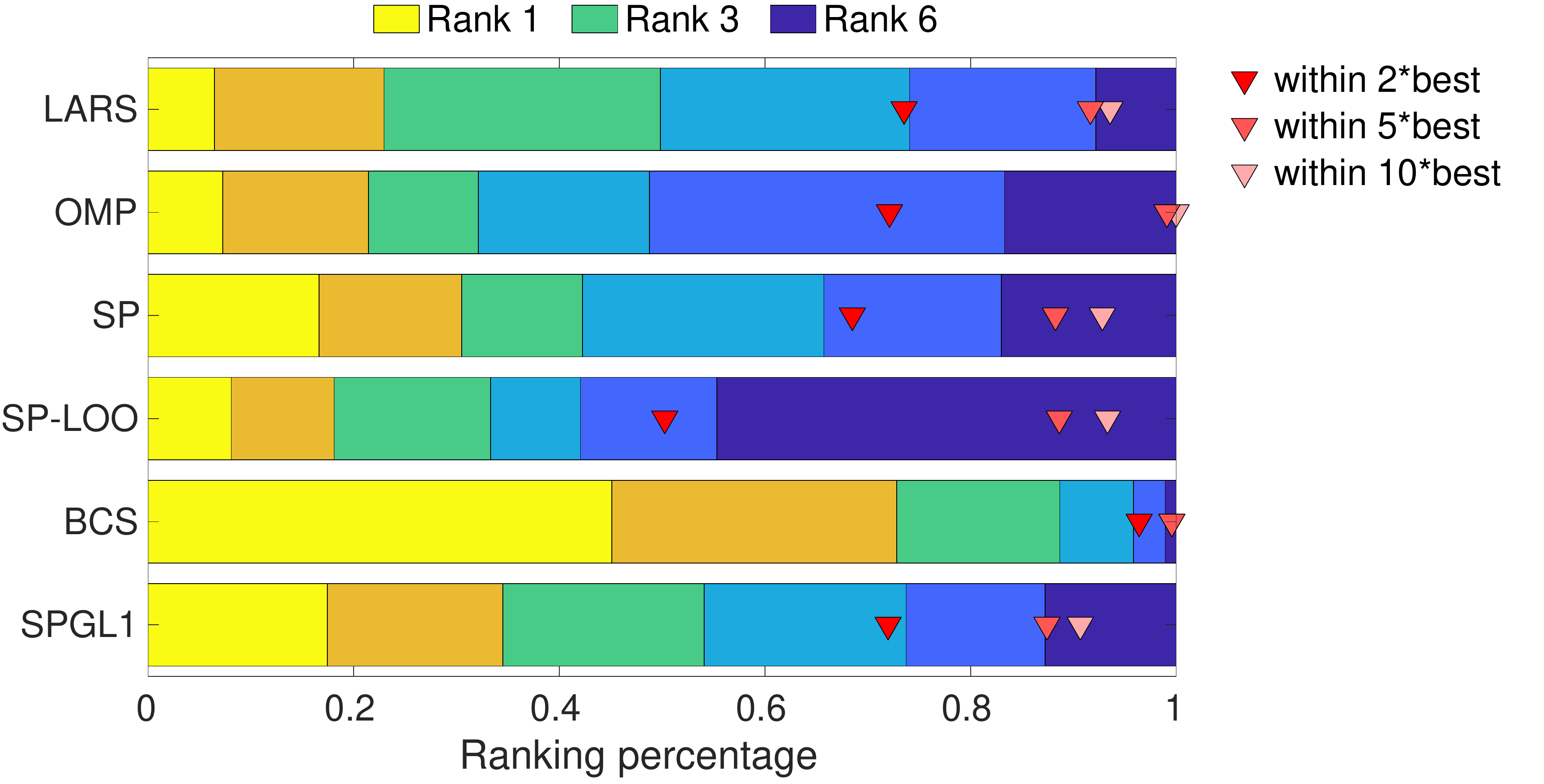}}
		\hfill
		\subfloat[][high-dim models ($d \geq 20$), large ED sizes]
		{\includegraphics[width=.435\textwidth, trim = 0 0 8.2cm 0, clip]{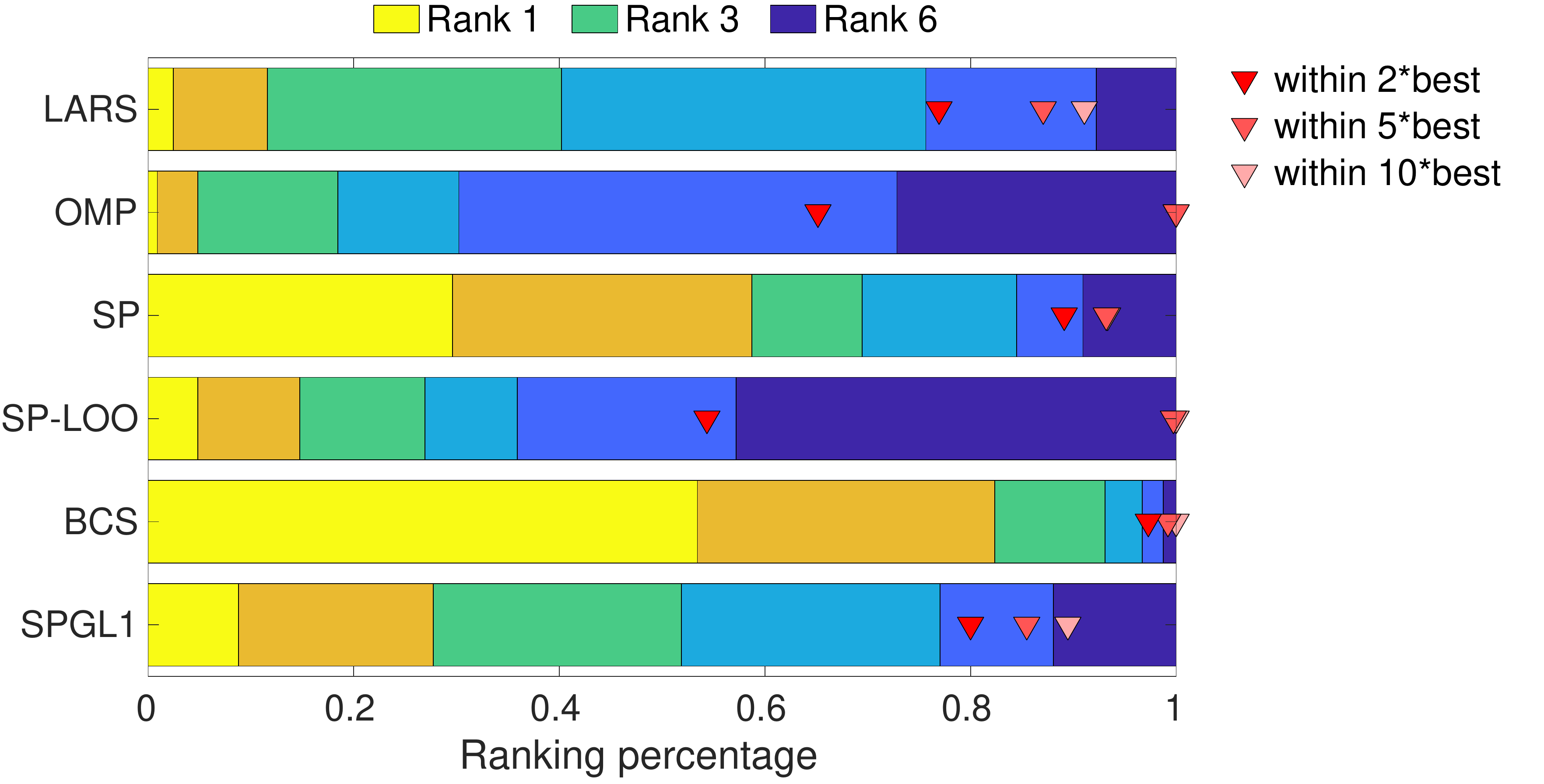}}
		\caption{\changed{Aggregated results for the solver benchmark (Section~\ref{sec:results_only_solvers}), separately for low-dimensional (a),(b) and high-dimensional models (c),(d). Left: small ED sizes. Right: large ED sizes. For each model and experimental design, the ranking of the six solvers is determined. The stacked bar plots visualize how often each solver achieved the respective rank. The triangle markers in hues of red additionally demonstrate in how many runs the obtained relative MSE was within a factor of $\{2,5,10\}$ of the smallest relative MSE achieved on this experimental design.} 
		}
		\label{fig:results_more_models_aggregated}	
	\end{figure}

	To objectively assess the performance of the methods, we now aggregate the results across models. Since all solvers are tested on the same set of experimental designs, we can determine the ranking of solvers for each experimental design (50 replications $\times$ $6-7$ ED sizes for the four spotlight models, and 30 replications $\times$ $5-7$ ED sizes for each of the seven additional models, resulting in 2620 PCEs) and count how often each solver achieved each rank. This is displayed in Figure~\ref{fig:results_more_models_aggregated} in the form of stacked bar plots, where the counts are given as percentages. The counts have been normalized by the number of replications and ED sizes used for each model, so that each of the models contributes equally to the final percentages.
	
	This ranking alone, however, does not provide a complete picture; e.g., a solver ranked last can be off by orders of magnitude or barely worse than the best-performing one. Therefore, we added an additional set of triangle markers detailing for how many of the EDs
	the respective solver returned a result that was within two, five, or 10 times the smallest relative MSE attained by any of the six solvers on the same ED. For example, the red triangle in the top row of Figure~\ref{fig:results_more_models_aggregated1} indicates that in ca.\ 30\% of the runs, the solution returned by LARS had a validation error that was at most twice as large as that of the best solution on the same ED.
	
	We have grouped the analysis results separately for the 6 low-dimensional ($d \leq 10$) and the 5 high-dimensional ($d \geq 20$) models, because we observed that dimension had a significant impact on the rankings, and this information is readily available even for black-box models. We also analyze small and large ED sizes separately, where the first half of considered ED sizes (3--4) are regarded small, and the second half large.
	
	We make the following observations:
	\begin{itemize}
		\item \textit{Low-dimensional models, small ED sizes}:
		BCS is the best solver most often (31\% of runs) and also most often within two times the smallest error (72\% of runs). However, it comes within one order of magnitude of the smallest relative MSE (a property which we here call \textit{robustness}) in only 87\% of runs, while \SPloo{} achieves this in 94\% of runs. In the two former metrics, \SPloo{} comes second. 
		LARS and SP perform worst, followed by OMP and SPGL1. Here and in the following, we observe that OMP and \SPloo{} often result in quite robust solutions.
		
		\item \textit{Low-dimensional models, large ED sizes}:
		\SPloo{} outperforms the other solvers in all respects. It provides the smallest relative MSE of all solvers in 42\% of runs, is in 73\% of runs within two times of the smallest relative MSE, and is even in 99\% of runs within one order of magnitude of the smallest relative MSE.
		For the other solvers, the statistics confirm the observations highlighted by the spotlight models: LARS and SPGL1 overall do not perform well. SP and BCS come close to the best solution quite often, but are less robust, whereas OMP is robust but often not as close to the best solution.
		
		\item \textit{High-dimensional models}:
		For the small as well as the large ED sizes, we see that BCS performs exceptionally well. It is the best solver in 45\% (53\%) of runs and comes within two times the smallest relative MSE in even 96\% (97\%) of runs. Both BCS and OMP attain in \textit{all} cases a relative MSE within 10 times the smallest relative MSE. SP performs better for large rather than small ED sizes. 
		While \SPloo{} performed best for the low-dimensional models, here it shows poor performance.
		Note that all solvers come within 10 times the smallest relative MSE in more than 89\% of all runs, showing that the choice of solver has a smaller impact for high-dimensional than for low-dimensional models.
		
	\end{itemize}

\subsection{Results: Comparison of sampling schemes together with solvers}
\label{sec:results_sampling}

	We pair the five solvers LARS, OMP, SP, \SPloo{}, and BCS%
	\footnote{Since SPGL1 did not perform well in the previous section, and is quite slow, we do not include it further in this benchmark.}
	with the sampling schemes MC, LHS, coherence-optimal, and D-optimal based on a coherence-optimal candidate set. We use the abbreviations coh-opt and D-opt(coh-opt) for the latter two.
	Since \SPloo{} performed poorly for high-dimensional models in the solver benchmark of the previous section, we do not consider it here for the high-dimensional models.
	
	We run the benchmark for the low-dimensional models Ishigami, undamped oscillator, borehole, damped oscillator, and wingweight function, and for the high-dimensional models Morris function, two-dimensional diffusion, one-dimensional diffusion, and \changed{100D} function. The truss model has Gumbel input, for which we (as of now) cannot construct a coherence-optimal sample. The same holds for the structural frame model with its dependent input. 
	
\begin{figure}[htbp]
		\centering
		\subfloat[][$d \leq 10$, small ED sizes]
		{\includegraphics[width=.55\textwidth]{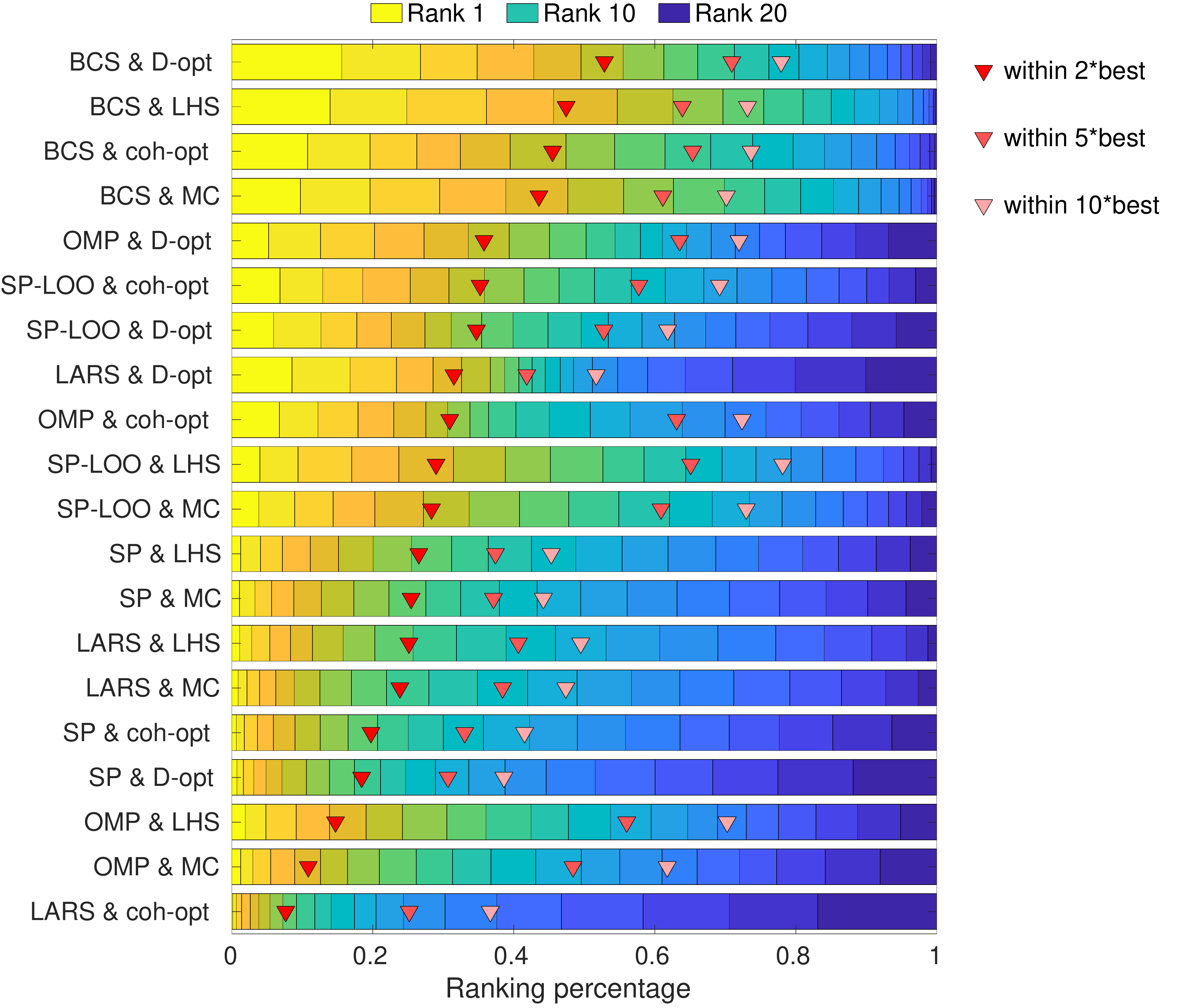}
			\label{fig:results_sampling_aggregated_small}}
		\subfloat[][$d \leq 10$, large ED sizes]
		{\includegraphics[width=.43\textwidth,trim=0 0 8.5cm 0, clip]{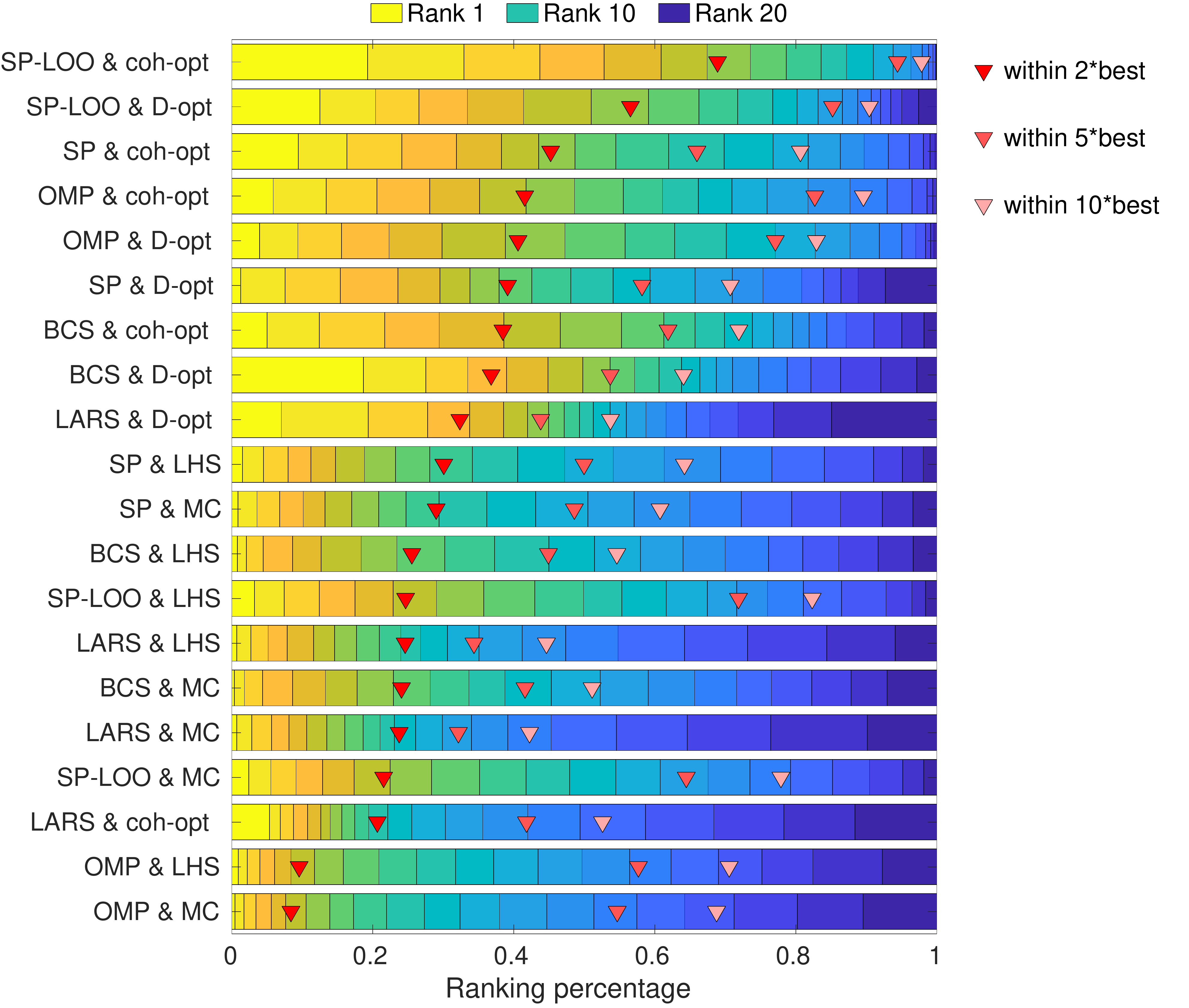}
			\label{fig:results_sampling_aggregated_large}}
		\\
		\subfloat[][$d \geq 20$, small ED sizes]
		{\includegraphics[width=.55\textwidth]{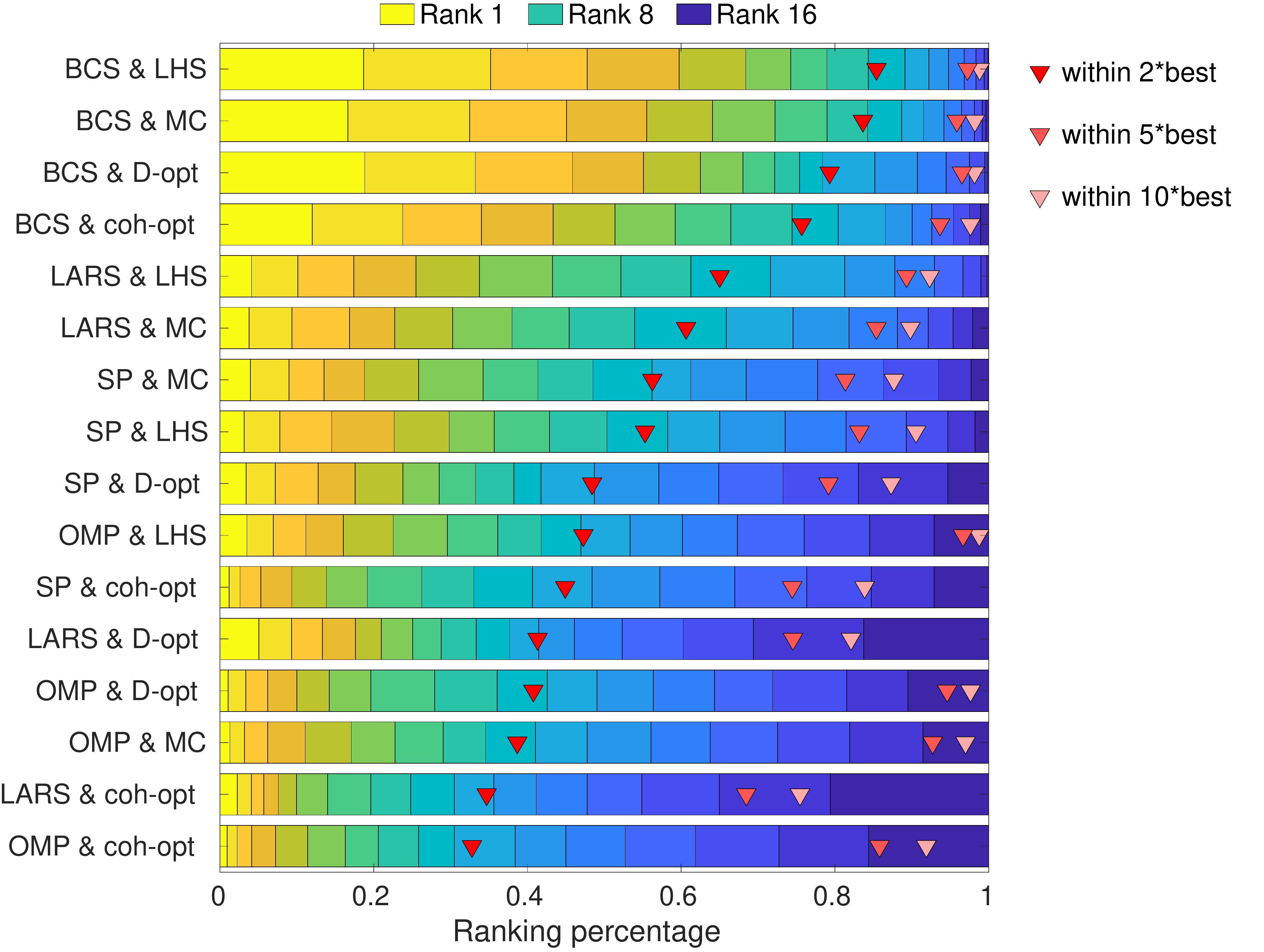}
			\label{fig:results_sampling_aggregated_small_highdim}}
		\subfloat[][$d \geq 20$, large ED sizes]
		{\includegraphics[width=.43\textwidth,trim=0 0 8.3cm 0, clip]{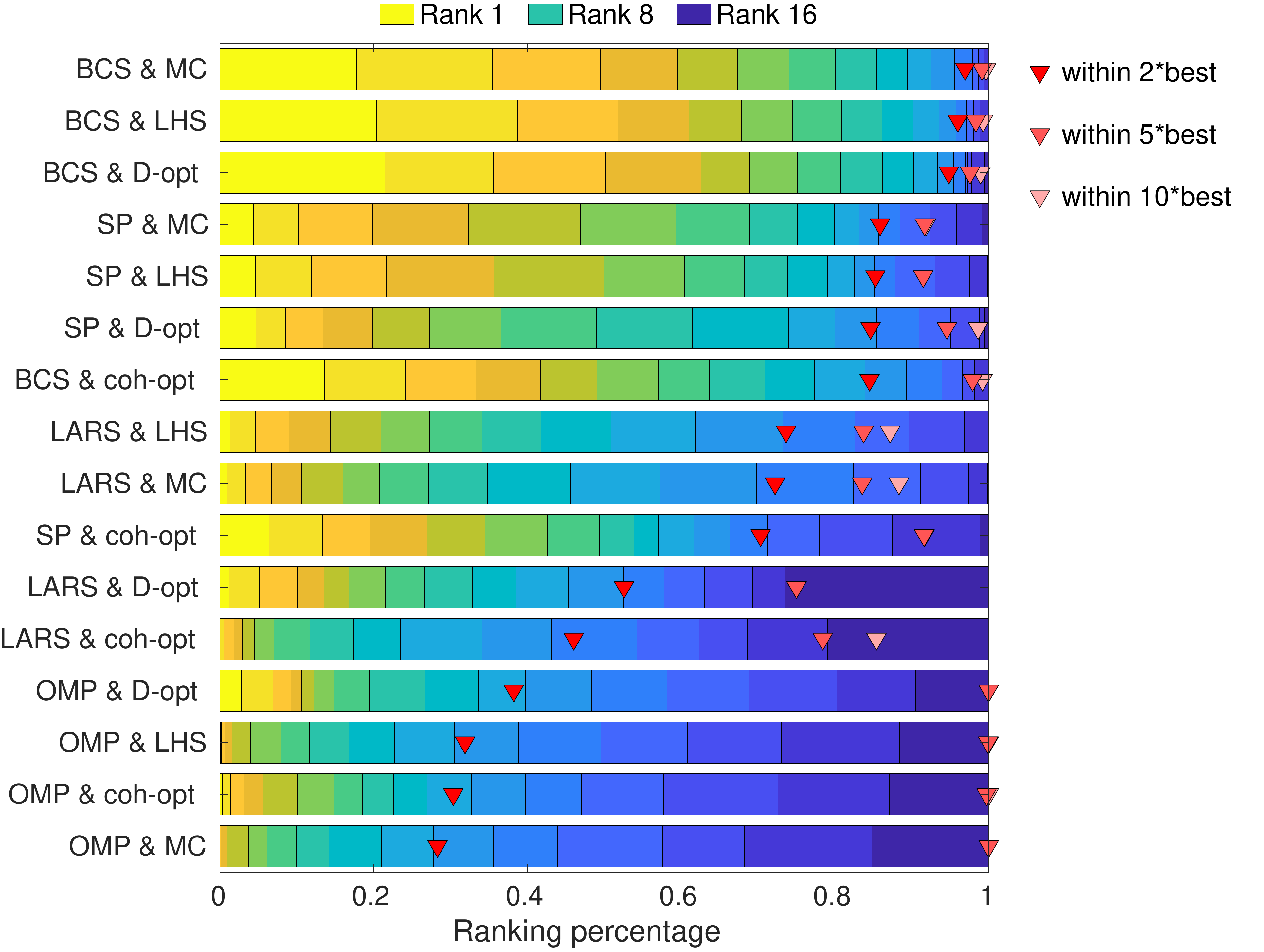}
			\label{fig:results_sampling_aggregated_large_highdim}}
		\caption{\changed{Aggregated results for the five low-dimensional models Ishigami, undamped oscillator, borehole, damped oscillator, and wingweight (top), and for the four high-dimensional models Morris function, structural frame, two-dimensional diffusion, and \changed{100D} function (bottom). Separately for small (a),(c) and large (b),(d) experimental designs. For the low-dimensional (high-dimensional) case, we investigate five (four) solvers and four sampling schemes, resulting in 20 (16) combinations. For each model and repetition, the ranking of all the combinations is determined (note that as opposed to Figure~\ref{fig:results_more_models_aggregated}, here the comparison is done on different EDs, which are matched randomly. Results are bootstrapped four times by random permutations to increase robustness). The stacked bar plots visualize how often each combination achieved the respective rank. The triangle markers in hues of red additionally demonstrate in how many runs the obtained relative MSE was within a factor of $\{2,5,10\}$ of the smallest relative MSE achieved in this comparison.
				The combinations are sorted by the percentage of runs in which they achieved a relative MSE within two times the smallest relative MSE of the respective random pairing (red triangle marker).
				Plots of relative MSE against ED size for the spotlight models can be found in Appendix~\ref{app:additional_results}, Figures~\ref{fig:results_ishigami_additional}--\ref{fig:results_highdimfct_additional}.}}
		\label{fig:results_sampling_aggregated}	
	\end{figure}

	Detailed boxplots of the relative MSE against ED size for the spotlight models, showing how each solver performed when paired with the sampling schemes, can be found in Appendix~\ref{app:additional_results}, Figures~\ref{fig:results_ishigami_additional}--\ref{fig:results_highdimfct_additional}.
	For the sake of readability, in this section we only show results that are aggregated over models, separately for the low- and high-dimensional cases. 
	For every model and repetition index, we determine the relative ranking of the 20 (16) combinations (5 (4) solvers $\times$ 4 sampling strategies). We also determine which of the combinations came within a factor of $\{2,5,10\}$ of the smallest relative MSE among this set (robustness). 
	Then, we count how often each combination achieved each rank, and how often each combination achieved a relative MSE within a factor of the smallest relative MSE. 
	The results are displayed in Figure~\ref{fig:results_sampling_aggregated} in the form of stacked bar plots for the ranks, with triangle markers denoting the percentage of robust runs. 
	The combinations are sorted by the percentage of runs in which they achieved a relative MSE within two times the smallest relative MSE, because we find that this metric is a good compromise between performance and robustness.
	We analyze small and large ED sizes separately, where the first half of considered ED sizes (3--4) are regarded as small and the second half as large.

	Note that, as opposed to the aggregated results in Figure~\ref{fig:results_more_models_aggregated}, where the solvers are compared on the same experimental designs, here the comparison is done on \textit{different} experimental designs, which are matched randomly.%
	\footnote{\changed{The relative MSE of each combination is interpreted as a random variable $E_\text{solver}^\text{sampling}$, where the randomness is induced by the randomness in the experimental design. E.g., $E_\text{BCS}^\text{LHS}$ is the random variable of relative MSE attained by BCS applied to an LHS design of specified size. The reference error (``smallest relative MSE'') $E^*$ is a random variable as well, defined as the minimum over \textit{one realization} of each combination of methods: $E^* = \min_{s \in \text{sampling}, t \in \text{solvers}} E_{s}^{t}$. The plots in Figure~\ref{fig:results_sampling_aggregated} are therefore read as follows: e.g., in the low-dimensional, small ED size case (Figure~\ref{fig:results_sampling_aggregated_small}): 
			$ \Prob{E_\text{BCS}^\text{D-opt} = E^*} = 0.16, \quad \Prob{E_\text{BCS}^\text{D-opt} \leq 10 E^*} = 0.78.$%
		}} 
	Our results in Figure~\ref{fig:results_sampling_aggregated} are bootstrapped four times using random permutations of the replication index, corresponding in total to 250 replications, to minimize the influence of this randomness (which is in any case not large, as can be seen from permutation tests).
	
	\changed{From Figure~\ref{fig:results_sampling_aggregated} and Figures~\ref{fig:results_ishigami_additional}--\ref{fig:results_highdimfct_additional}, we} make the following observations:
	\begin{itemize}
		\item \changed{There can be considerable differences in the performance of different combinations of solvers and sampling schemes. The differences are larger for low-dimensional models, visible in the spread of triangle markers in Figures~\ref{fig:results_sampling_aggregated_small} and \ref{fig:results_sampling_aggregated_large}.
			For high-dimensional models, many combinations find an error that is close to the smallest error, which can be seen from the red triangle markers at high percentages in Figures~\ref{fig:results_sampling_aggregated_small_highdim} and \ref{fig:results_sampling_aggregated_large_highdim}, and from the clustered boxplots in Figures~\ref{fig:results_diffusion2D_additional} and \ref{fig:results_highdimfct_additional}, (e) and (f).}
		
		\item MC and LHS perform comparably, and for the high-dimensional models almost identically. For the low-dimensional models, LHS sampling has, in most cases, median error and variability that are the same as or smaller than MC. 
		This is consistent with the literature \citep{Shields2016}.
		\changed{These observations are confirmed by the plots in Figure~\ref{fig:results_sampling_aggregated}: for almost every solver, the combination with LHS is slightly better than the corresponding one with MC in each of the metrics. 
		}
		
		\item \changed{The advanced sampling schemes coh-opt and D-opt(coh-opt) show a clear advantage over MC and LHS sampling for low-dimensional models and large experimental designs (consistent with theoretical considerations and numerical experiments \citep{Hampton2015}). For low-dimensional models and small experimental designs, they show mixed performance; for high-dimensional models, they perform the same as or worse than LHS and MC.}
		
		\changed{It is known that coh-opt sampling leads to a greater improvement over MC sampling for low-dimensional, high-degree expansions than for high-dimensional, low-degree expansions \citep{Hampton2015,Alemazkoor2018}.
			Note also that all numerical experiments in the literature testing coh-opt sampling were performed in $d \leq 30$ dimensions, using only models with uniform input, or manufactured sparse PCE, i.e., polynomial models with an exactly sparse representation \citep{Hampton2015,Hampton2015b,Alemazkoor2018,Diaz2018}.}

		\changed{Both} coh-opt and D-opt are sampling methods that aim to improve properties of the regression matrix. They are adapted to the candidate basis. If the candidate basis is large and contains many regressors that are not needed in the final sparse expansion, this adaptation might even deteriorate the solution.
		
		\item \changed{BCS is one of the best-performing solvers, almost regardless of sampling scheme. The exceptions are low-dimensional models with large experimental designs, where \SPloo{} with coh-opt sampling outperforms all other solvers. This might be related to BCS plateauing earlier than other solvers (see Figure~\ref{fig:onlysolvers_ishigami} and \ref{fig:onlysolvers_borehole}). It seems BCS is preferable whenever the information content is low (small ED sizes or high-dimensional models).}
		
		\item \changed{OMP and \SPloo{} are generally quite robust (within one order of magnitude of the best solution). However, OMP often does not come close to the best solution, especially when paired with LHS or MC. BCS is more robust for high-dimensional models than for low-dimensional models. LARS and SP show mixed performance, with LARS being one of the least robust solvers.}
		
		\item \changed{Aggregating the results for each sampling scheme separately (not shown), we observe that the behavior of the solvers is very similar in terms of ranking and robustness to the behavior observed on LHS (Figure~\ref{fig:results_more_models_aggregated}), suggesting that the ranking of solvers is mostly independent of the sampling scheme.
	}
		
	\end{itemize}

	Note that the results in Figures~\ref{fig:results_ishigami_additional}, \ref{fig:results_borehole_additional}, and \ref{fig:results_highdimfct_additional} exhibit plateauing for larger sample sizes. This indicates that the maximal accuracy achievable with this set of basis functions has been reached. \changed{Using a larger basis might lead to more accurate solutions, if it contains an important regressor that was previously missing. However, note that a larger basis can also lead to less accurate solutions: when the experimental design size is held fixed while a larger basis is used, } 
	the ratio of experimental design points to basis functions \changed{is}
	smaller, and the properties of the regression matrix might deteriorate.

	\subsection{Results: Comparison of sampling schemes together with solvers, using a smaller candidate basis}
	\label{sec:results_nearopt}
	
\begin{figure}[htbp]
		\centering
		\subfloat[][OMP]{\includegraphics[width=.49\textwidth, height=0.25\textheight, keepaspectratio]{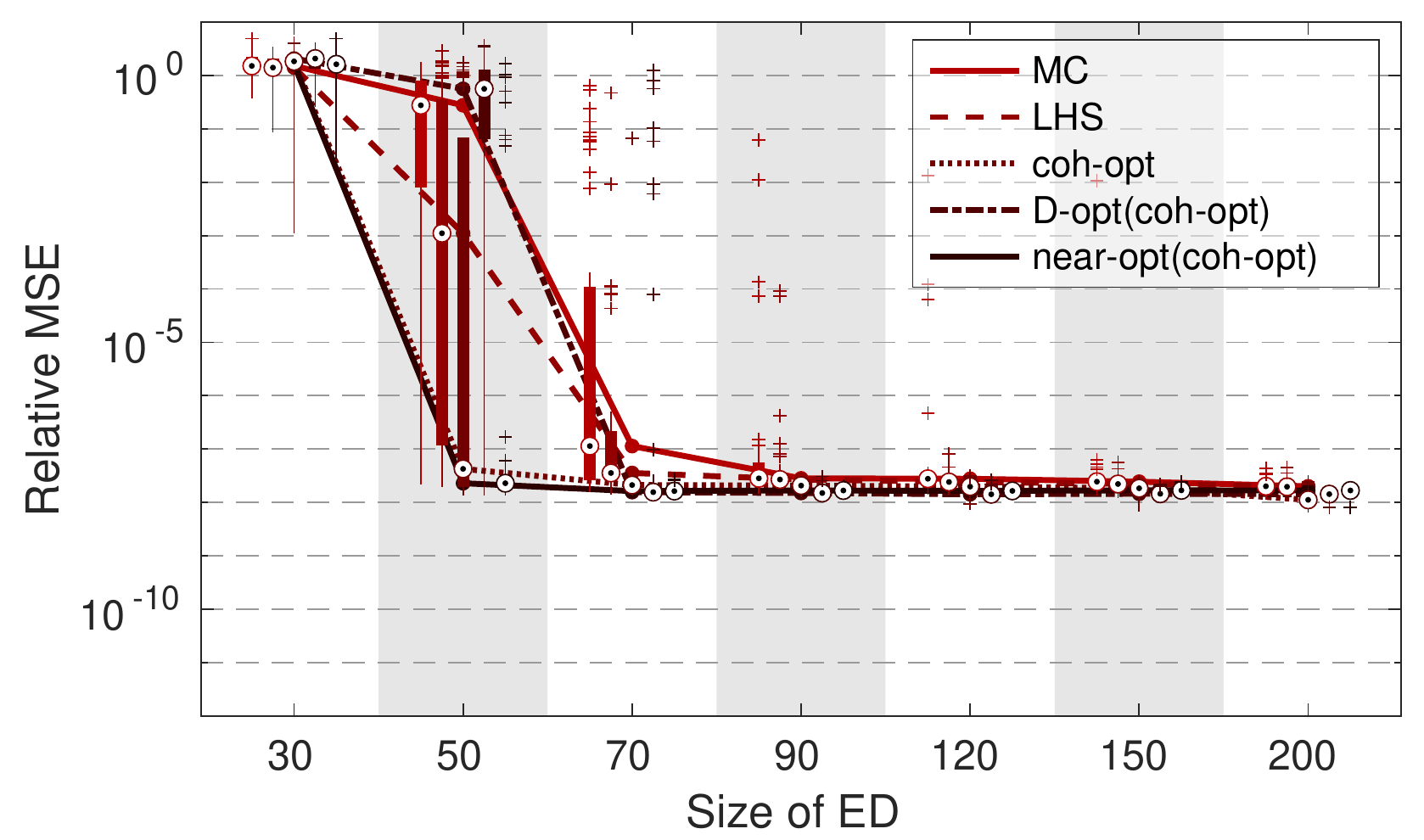}}
		\hfill
		\subfloat[][\SPloo{}]{\includegraphics[width=.49\textwidth, height=0.25\textheight, keepaspectratio]{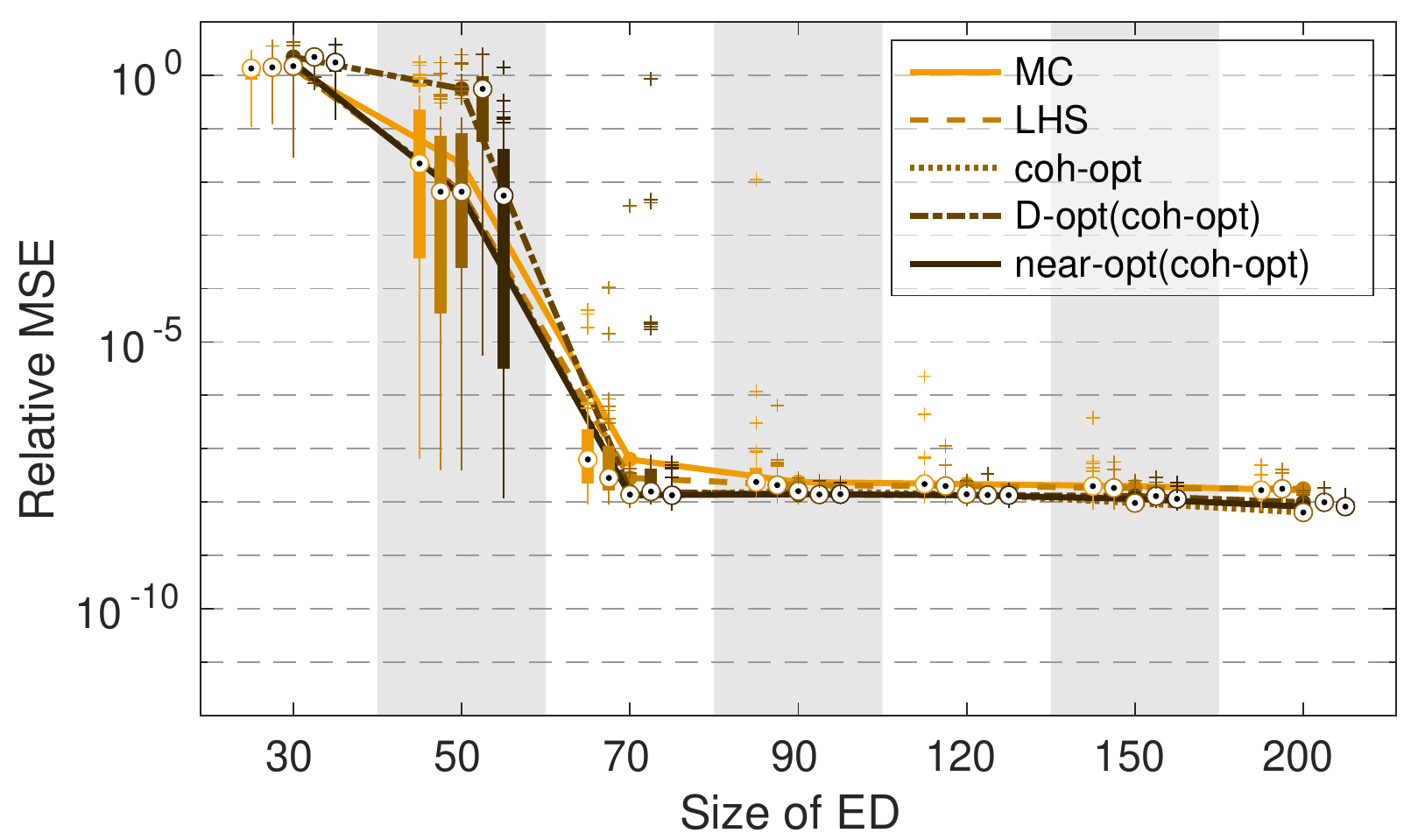}}
		\caption{\changed{Results for the Ishigami model with a smaller basis ($d = 3, p = 12, q=1$). Results for two sparse solvers and five experimental design schemes. 50 replications. For the remaining plots, see Figure~\ref{fig:results_smallbasis_additional} in Appendix~\ref{app:additional_results}.}}
		\label{fig:results_ishigami_smallbasis}	
	\end{figure}

\begin{figure}[htbp]
		\centering
		\subfloat[][OMP]{\includegraphics[width=.49\textwidth, height=0.25\textheight, keepaspectratio]{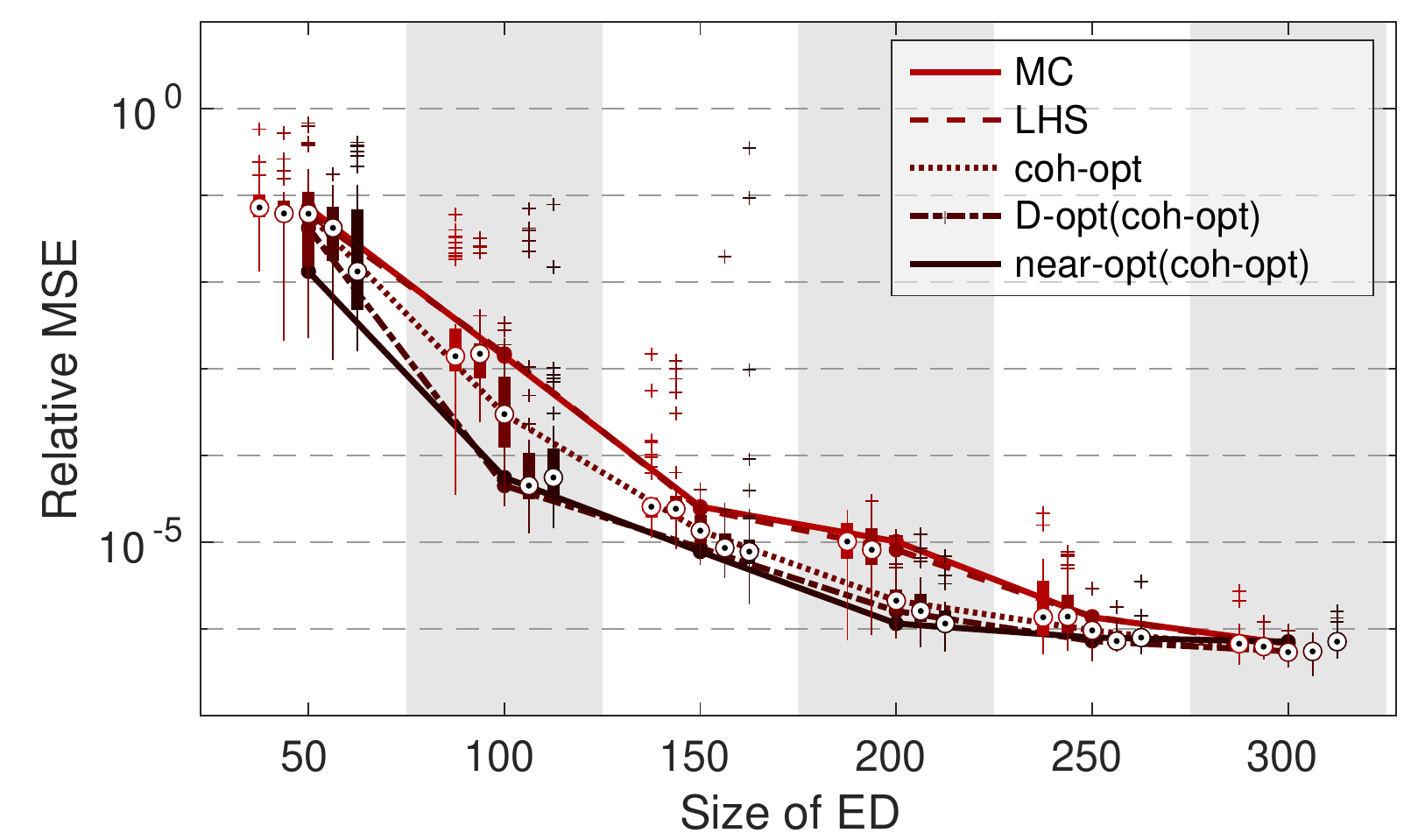}}
		\hfill
		\subfloat[][\SPloo{}]{\includegraphics[width=.49\textwidth, height=0.25\textheight, keepaspectratio]{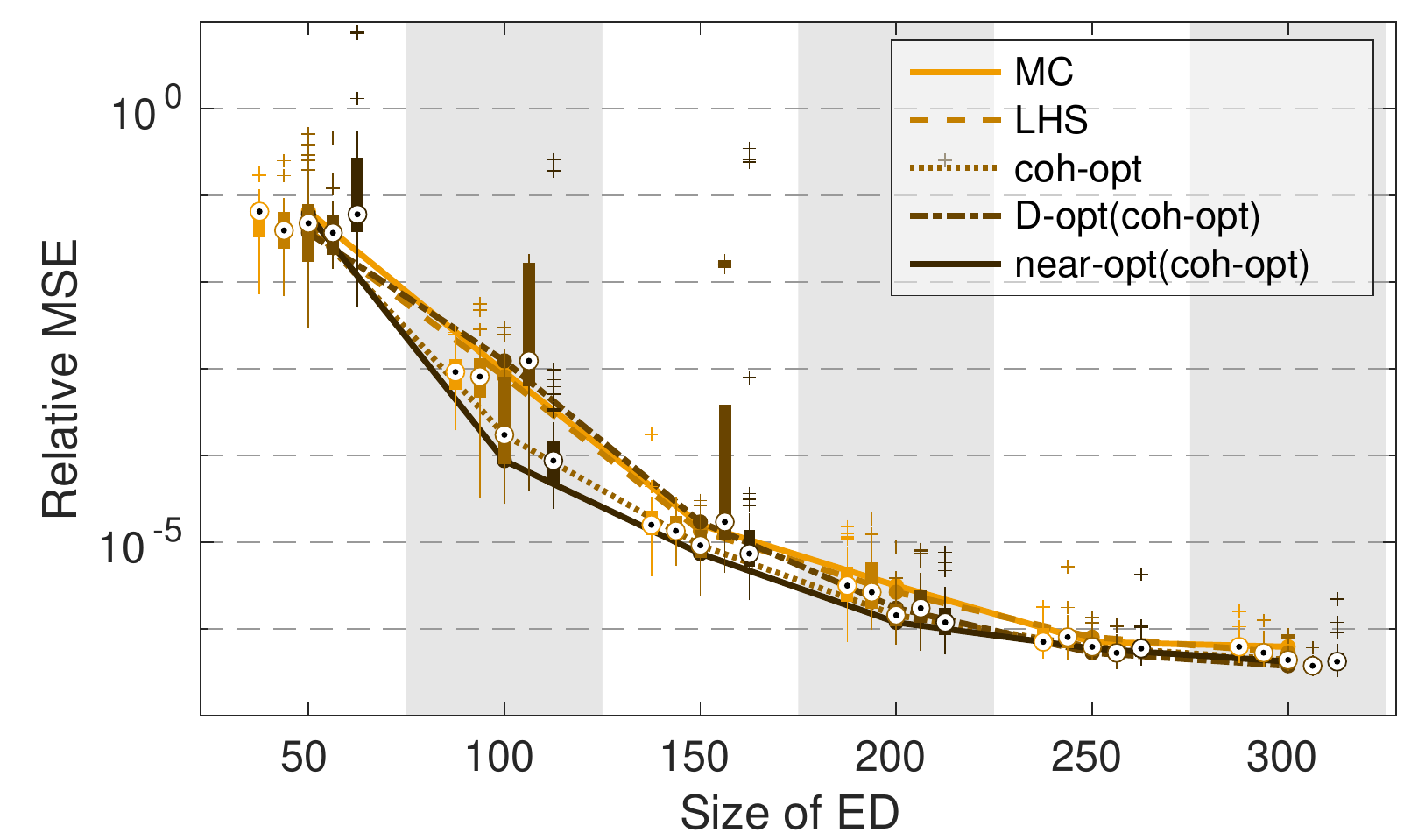}}
		\caption{\changed{Results for the borehole model with a smaller basis ($d = 8, p = 4, q = 1$). Results for two sparse solvers and five experimental design schemes. 50 replications. For the remaining plots, see Figure~\ref{fig:results_smallbasis_additional} in Appendix~\ref{app:additional_results}.}}
		\label{fig:results_borehole_smallbasis}	
	\end{figure}

	We repeat the experiments from the previous section for the Ishigami and borehole models, using a smaller candidate basis for which near-optimal sampling is feasible. The tested solvers are LARS, OMP, SP, \changed{\SPloo{},} and BCS. We use the sampling schemes MC, LHS, coh-opt, D-opt(coh-opt) and near-opt(coh-opt). 
	\changed{Boxplots of relative MSE against ED size} are shown in Figures~\ref{fig:results_ishigami_smallbasis} and \ref{fig:results_borehole_smallbasis}.
	\changed{For the sake of conciseness, we only show the combinations involving OMP and \SPloo{}. The remaining plots are provided in Appendix~\ref{app:additional_results}, Figure~\ref{fig:results_smallbasis_additional}.}
	
	We observe the following:
	\begin{itemize}
		\item Since the basis is smaller, the relative MSE reaches a plateau already for smaller experimental design sizes. 
		
		\item Most qualitative observations regarding solver and sampling performance are the same as in the previous section, where a larger basis was used.
		\item Near-optimal sampling often achieves the same or a slightly smaller error than \changed{coh-opt} sampling, which is consistent with \citep{Alemazkoor2018}. In many cases, near-optimal sampling achieves the smallest median error. For the Ishigami \changed{model}, near-opt\changed{imal sampling} additionally exhibits small variability\changed{, while for} the borehole model, \changed{it} has a rather large spread, i.e., several outliers.
	\end{itemize}

	\FloatBarrier

\section{\changed{Discussion and conclusions}}
\label{sec:conclusion}

In this paper, we investigated sparse PCE methods with the goal of computing accurate surrogate models based on a few model evaluations. 

We presented a literature survey and a framework describing the general computation procedure for sparse PCE. 
We have seen that the existing literature on sparse PCE can be fit into this framework and that methods developed for different components of the framework can naturally be combined. 

In order to give recommendations to practitioners who want to use sparse PCE surrogates for their applications, we performed a numerical benchmark based on \changed{11} example functions which are intended to resemble real-world engineering problems presenting different challenges.
We tested several popular sparse solvers and sampling schemes on a fixed set of basis functions, using a range of experimental design sizes and \changed{30--}50 replications, and
made the following observations:
\begin{itemize}
	\item The choice of sampling scheme and sparse regression solver can make a difference of up to several orders of magnitude for the relative MSE. Mostly, the rankings of solvers and sampling schemes seem to be independent of one another: an experimental design that works best for one solver will also perform well with other solvers, and the ranking of solvers looks similar independent of which sampling scheme is used. \changed{Both solvers and} sampling schemes make a greater difference for low-dimensional models.
	\item \changed{For low-dimensional models and small ED sizes, the solver BCS performs best, regardless of sampling scheme (with D-opt(coh-opt) being slightly preferable), while the solver \SPloo{} (a variant of SP) appears to be especially robust.}
	\item \changed{For low-dimensional models and large ED sizes, \SPloo{} together with coherence-optimal sampling outperforms all other combinations.} 
	\item \changed{For low-dimensional models, and when the basis is small enough to make it feasible, near-optimal sampling outperforms all other sampling schemes, regardless of the solver}.
	\item \changed{For high-dimensional models, BCS is by far the best solver. All solvers perform better when paired with LHS; in other words, no advanced sampling scheme appears competitive compared to LHS for such problems, whatever the solver used.
	} 
\end{itemize}
\changed{The benchmark results demonstrate that in costly, real-world applications it is worth choosing the sparse PCE training strategy carefully, since the methods can make a substantial difference in the quality of the resulting surrogate.
	While a more accurate surrogate model is generally desirable, in industrial applications it might have a higher impact for purposes such as optimization, rather than, e.g., sensitivity analysis.
}

\changed{Our conclusions are based on a number of benchmark models, which we consider representative of engineering models in terms of dimensionality and complexity. Naturally, however, no selection of models can cover the whole space of engineering models.}
\changed{Further work would be required to understand the connection among model properties, basis choice, experimental design size, and sparse PCE techniques like solvers and sampling schemes.}

All results were obtained using a fixed basis based on a heuristic choice (see Section \ref{sec:benchmark_models}).
Generally, when the optimal degree of the basis is unknown, degree adaptivity (based on a cross-validation error) can be a useful strategy.
\changed{Due to time and space constraints, this} was not investigated in the present work.
Adaptivity critically depends on the availability of an accurate error estimator.
Some of the best solvers in this study (i.e., OMP and BCS) %do not provide such estimators, and 
tend to underestimate the generalization error (not shown in the plots), which might be a drawback in the setting of adaptive degree selection and might change the effect and ranking of solvers and sampling methods. 
\changed{For a detailed discussion and benchmark of basis-adaptive schemes, we refer the reader to \citep{LuethenIJUQ2021}.}

\changed{As evident from the extensive literature on the topic, sparse PCE is an already well-established technique, as well as an active field of methodological research. Recent innovations include Bayesian techniques for sparse PCE and the identification of suitable rotated coordinates for the expansion.
	Such innovative ideas are expected to lead to further improvements in the computation of sparse PCE, which will in turn benefit applications as well as all advanced methods that use sparse PCE as one of their building blocks (see, e.g., \citep{SchoebiIJUQ2015,Chatterjee2019,ZhuIJUQ2020,MarelliSPCE2021}).}

\changed{PCE is a popular metamodelling tool in the engineering community, and many different methods are available. Up to now, the choice of which of the many PCE methods to apply was mostly left to chance or the personal experience of the practitioner. 
	In our benchmark, we explored a significant set of methods that have received attention in the past few years. We hope that this work can serve as a basis for further benchmarking efforts, in order to identify which of the many available methods are most suitable for real-world problems. These might include sequential enrichment of experimental design, Gaussian adaptation of the input space, stepwise regression algorithms and many other ideas for sparse solvers, as well as methods for extremely high-dimensional problems.
}

\changed{Our benchmark code is available on request. The solvers BCS and \SPloo{} will be made available in the 1.4.0 release of UQLab. For a description of how to add custom sampling schemes and sparse solvers for PCE to UQLab, we refer the reader to the supplementary material accompanying this paper.}
\changed{
	To facilitate easier benchmarking of PCE techniques on a large number of examples in a standardized setup, we are actively engaging in designing and developing a benchmarking platform for surrogate modelling methods similar to the UCI machine learning repository%
	\footnote{\url{https://archive.ics.uci.edu}}
	or the structural reliability platform RPrepo%
	\footnote{\url{https://rprepo.readthedocs.io/en/latest/reliability_problems.html}}
	where data sets, models, and methods can be made available for testing and benchmarking.}

\section*{Acknowledgments}
We thank Dr.\ Emiliano Torre for helpful discussions.

%\bibliographystyle{apalike}%siam(plain) % alpha, apalike show authors and year
%\bibliography{../my_bib,BiblioRSUQ,BiblioSudret}

\clearpage

\appendix

\section{Details on experimental design sampling techniques}
\label{app:ED}

It depends on certain properties of the regression matrix $\ve\Psi$ whether or not sparse regression techniques are able to find the true sparse solution of a linear system of equations (assuming that it exists). 
In the context of polynomial chaos, the entries of the regression matrix are the basis polynomials evaluated at the design points. The basis polynomials are determined by the distribution of the input random variables and the choice of the index set $\ca$, while the design points $\{\ve x^{(j)}\}_{j=1}^N$ can be chosen freely from the input space to optimize properties of the resulting regression matrix.

In the following, we present sampling schemes that have been proposed in the literature for the computation of sparse PCE. Some of the schemes come with theoretical results about their performance for sparse PCE, others have heuristic justification or have guarantees for least-squares regression. They can be broadly grouped into three categories: 
\begin{itemize}
	\item Sampling according to the input distribution
	\begin{itemize}
		\item MC \citep{Doostan2011, Hampton2015}
		\item LHS \citep{McKay1979}
	\end{itemize}
	\item Sampling from a modified distribution \changed{(induced sampling)}
	\begin{itemize}
		\item asymptotic \citep{Hampton2015}
		\item coherence-optimal \citep{Hampton2015}
		\item Christoffel sparse approximation \citep{Jakeman2017}
	\end{itemize}
	\item Optimizing matrix properties
	\begin{itemize}
		\item D-optimal \citep{Diaz2018}
		\item S-optimal \citep{Shin2016a, FajraouiMarelli2017}
		\item near-optimal \citep{Alemazkoor2018}
	\end{itemize}
\end{itemize}
Some of the sampling schemes are nontrivial or costly to evaluate\changed{, or even not available for all input distributions}. However, the bottleneck in surrogate modelling for practical applications is typically the repeated evaluation of the model, which justifies the use of a complex sampling scheme if it allows better approximation with fewer samples.

\subsection{Sampling according to the input distribution}
This class of sampling methods consists of all methods that are oblivious to the 
\changed{choice of truncation set $\ca$}
and whose main objective is to distribute design points evenly in the quantile space. Heuristically, the more uniformly the points are distributed in the quantile space, the more information about the model is captured in the model evaluations, since no region of the input domain is forgotten. 

LHS is one technique for achieving a space-filling design. For each component of the input random vector $\ve X$, the corresponding quantile space is divided into $\NED$ intervals. In each interval, one point is sampled uniformly at random. Then, the points for each dimension are combined randomly into vectors and finally transformed into the input space using an isoprobabilistic transform. LHS can be shown to reduce the variance of linear regression estimates when the main effects are dominant, i.e., when the most important terms have interaction order one \citep{Shields2016}.
LHS can be combined with heuristic criteria such as the maximin distance strategy, where among several random LHS designs the one with the largest minimal pairwise distance between points is chosen, to further improve on the space-filling property.

Stratified sampling is a related sampling technique in which the input space is divided into disjoint regions, called strata, from which points are sampled and weighted according to the probability mass of their stratum. 
Stratified sampling reduces the variance of statistical estimators \citep{McKay1979}.
There exists a range of methods between stratified sampling and LHS, called partially stratified sampling, which are able to reduce the variance of statistical estimators when interaction terms are dominant \citep{Shields2016}. The authors propose an additional method called
Latinized partially stratified sampling (LPSS) which combines LHS and stratified sampling with the aim tif minimizing the variance of the resulting estimator. It is especially beneficial when there is prior knowledge about which variable groups interact, and it has been used for several problems with input dimension $d = 100$.

MC sampling, i.e., sampling from the input distribution, is a special case of the coherence-based theory detailed in Section \ref{app:coherence} below, and bounds on the coherence and the associated number of points needed for recovery can be derived \citep{Doostan2011, Rauhut2012, Yan2012, Hampton2015}. 

\subsection{Sampling from a different distribution}

The ability of sparse regression to recover the true sparse solution (if it exists; otherwise it recovers the best sparse approximation) largely depends on the regression matrix. In the case of PCE, the entries of this matrix are the evaluations of the basis polynomials at the experimental design points. The points can be chosen in a way that improves the recovery properties of the matrix.

Several approaches exist in which the $\ell^1$-minimization problem is modified into a weighted problem and samples are drawn not from the input distribution, but from a suitable modified distribution. 
The idea of these approaches is as follows. Define a weight function $w(x): \Omega \to \R$ in a suitable way, which will be explained later. For an ED $\{\ve x^{(k)}\}_{k=1}^N$, define the diagonal matrix $\ve W = \text{diag}(w(\ve x^{(1)}) \enum w(\ve x^{(\NED)}))$. 
Then the following modified system is solved:
\begin{equation}
\min_{\ve c} \norme{\ve c}{1} \text{ s.t. } \norme{\ve W \ve \Psi \ve c - \ve W \ve y}{2} \leq \epsilon.
\end{equation}
Depending on $w(\ve x)$, this modification can improve or deteriorate the solution $\ve c$. 
Of course, the weight function is chosen to improve the solution.
The matrix $\ve W \ve \Psi$ can also be interpreted as the evaluation of a modified basis $\{ \tilde\psi_\alp(\ve x) = w(\ve x)\psi_\alp(\ve x)\}_\changedmath{\alp \in \ca}$.
To ensure orthonormality of the columns of $\ve W \ve \Psi$, the design points are drawn from a suitably modified input distribution $f_{\tilde{\ve X}}$.

\subsubsection{Coherence, isotropy, and weighted orthonormal systems}
\label{app:coherence}
In this section, we define concepts that are the basis for guarantees on accuracy and stability for different sampling distributions. We mainly follow the exposition in \citep{Hampton2015}.

In the setting of PCE, the \textit{coherence} of an orthonormal system $\{\psi_\alp\}_{\alp \in \ca}$ is defined by
\begin{equation}
\mu(\ca, \{\psi_\alp\}) = \sup_{\ve x \in \cd} \max_{\alp \in \ca} |\psi_\alp(\ve x)|^2.
\label{eq:coherence_app}
\end{equation}
For distributions for which this quantity would be $\infty$, such as a Gaussian distribution, see the remark below.

A second important concept is \textit{isotropy} \citep{Candes2011}: a random matrix, whose rows are chosen randomly following some distribution $\ve a \sim F_{\ve a}$, is isotropic if it holds that $\Esp{\ve{a}^T \ve a} = \mathbb{1}$. In the case of PCE, $F_{\ve a}$ is induced by propagating the input distribution $F_{\ve X}$ through the basis functions. By construction, the regression matrix of standard PCE is isotropic if the ED is sampled from the input distribution.
Under the assumption that the regression matrix $\ve\Psi$ is isotropic, 
the number of samples needed for perfect recovery of sparse solutions in the noiseless case is proportional to $\mu(\ca, \{\psi_\alp\}) s \log(P)$ with high probability \citep{Candes2011}, where $s$ is the sparsity of the solution vector and $P = |\ca|$ is the number of basis functions. A similar result holds in the noisy case.

Thus, an orthonormal system $\{\psi_\alp\}_{\alp \in \ca}$ with low coherence $\mu(\ca, \{\psi_\alp\})$ requires fewer samples for perfect recovery. The goal of coherence-optimal sampling is to find a weighted system 
$\{ \tilde\psi_\alp(\ve x) = w(\ve x)\psi_\alp(\ve x)\}_\changedmath{\alp \in \ca}$ 
that achieves $\mu(\ca, \{\tilde\psi_\alp\}) < \mu(\ca, \{\psi_\alp\})$
and is orthonormal with respect to some distribution $\tilde f_{\ve X}$.

The ideas of isotropy and coherence were applied to PCE by \citet{Hampton2015}, who construct an isotropic regression matrix with improved coherence as follows.
Let $B:\cd \to \R$ be the tight upper bound for the polynomial basis,
\begin{equation}
B(\ve x) = \max_{\alp \in \ca} |\psi_\alp(\ve x)|.
\end{equation}
Let $G:\cd \to \R$ be a loose upper bound with $G(\ve x) \geq B(\ve x) \ \forall \ve x \in \cd$. ($G$ is useful because using a simple expression for the upper bound can in some cases result in $\tilde{f}_{\ve X}$ being a well-known distribution that can be sampled from easily.)
Define a new probability distribution $\tilde{f}_{\ve X}(\ve x)$ by
\begin{equation}
\tilde{f}_{\ve X}(\ve x) = c^2 G(\ve x)^2 f_{\ve X}(\ve x),
\label{eq:cohimprovedpdf}
\end{equation}
where $c = \left( \int_\Omega f_{\ve X}(\ve x) G(\ve x)^2 \di{\ve x} \right)^{- \frac12}$ is the normalizing constant.
Then, with the weight function
\begin{equation}
w(\ve x) = \frac{1}{c G(\ve x)},
\label{eq:weight_cohopt}
\end{equation}
the set of functions $\{\tilde\psi_\alp(\ve x) = w(\ve x) \psi_\alp(\ve x)\}_\changedmath{\alp \in \ca}$ 
is an orthonormal system with respect to the distribution $\tilde{f}_{\ve X}$. This follows directly
from the orthonormality of $\{\psi_\alp\}_{\alp \in \ca}$ with respect to $f_{\ve X}$. 
Furthermore, if $G = B$, the coherence $\mu(\ca, \{\tilde\psi_\alp\})$ is minimal.
%In other words, the PC matrix created from the functions $\{w(\xi) \psi_\alp(\xi)\}_{\alp \in \ca}$ by sampling from the distribution $\tilde{f}_X$ is isotropic.

\paragraph{Remark} Some polynomial bases (e.g. Hermite polynomials) do not have a finite upper bound. It is still possible to obtain similar results by considering a smaller domain $\cs \subset \cd$ on which the upper bound is finite and the isotropy is still approximately fulfilled. The modified probability distribution is then $\tilde{f}_{\ve X}(\ve x) = c^2 G(\ve x)^2 f_{\ve X}(\ve x) \mathbb{1}_{\cs}(\ve x)$.

\subsubsection{Sampling using a loose upper bound ("asymptotic sampling")}
In the case of Legendre and Hermite polynomials, and using a certain loose upper bound \changedmath{$G(\ve{x}) \geq \max_{\alp \in \ca} |\psi_\alp(\ve{x})|$}, analytical expressions for distributions with improved coherence can be obtained \citep{Hampton2015}.

In the case of Legendre polynomials on $[-1,1]^d$, a loose upper bound on the polynomials is given by 
$G(\ve x) \propto \prod_{i=1}^d(1-x_i^2)^{-\frac14}$, 
which leads to the Chebyshev distribution
$\tilde{f}_{\ve X}(\ve x) = \prod_{i=1}^d \frac{1}{\pi \sqrt{1-x_i^2}}$ 
and to the weight function
$w(\ve x) = \prod_{i=1}^d(1-x_i^2)^{\frac14}$.

In the case of Hermite polynomials for standard Gaussian variables, a loose upper bound on the polynomials is given by $G(\ve x) \propto \exp(\frac14 \norme{\ve x}{2}^2)$,
and the subset $\cs$ is chosen to be the $d$-dimensional ball with radius $\sqrt{2}\sqrt{2p+1}$. This leads to a uniform distribution $\tilde{f}_{\ve X}$ on $\cs$ and to the weight function $w(\ve x) = \exp(- \frac14 \norme{\ve x}{2}^2)$.

Additionally, asymptotic distributions for Laguerre polynomials (corresponding to the Gamma distribution) and for Jacobi polynomials (Beta distribution) have been implemented in the software package COH-OPT \citep{COHOPT}.

For Legendre polynomials, asymptotic sampling has a smaller coherence than standard sampling in the case $d < p$ (asymptotically). In the case $d > p$, which is more common in applications, standard sampling has (asymptotically) a smaller coherence.
According to theory, the sampling scheme with smaller coherence should exhibit better recovery rates.
This is confirmed numerically \cite[section 5.1]{Hampton2015}. For Hermite polynomials, the same observation is made.

\subsubsection{Coherence-optimal sampling}
The choice $G=B$ leads to the minimum possible  coherence $\mu(\ca, \{\tilde\psi_\alp\})$
\cite[Theorem 4.5]{Hampton2015}.
$B$ is simple to evaluate for a single point $\ve x \in \cd$, but its functional form is in general not known. Therefore, \citet{Hampton2015} suggest sampling $\tilde{f}_{\ve X} \propto B^2 f_{\ve X}$ using Markov chain Monte Carlo (MCMC) sampling with proposal distribution equal to the input distribution in the case $d \geq p$ and equal to the asymptotic distribution in the case $d < p$. 
The resulting (unnormalized) weights are $w(\ve x) = \frac{1}{B(\ve x)}$.
As expected from theory, numerical examples indicate that coherence-optimal sampling achieves better recovery and a smaller error in various norms than both standard and asymptotic sampling \citep{Hampton2015}.
Coherence-optimal sampling can be shown to have good properties also when used as a sampling scheme for least-squares regression \citep{Hampton2015b}.

\label{app:coh-opt-rejection}

MATLAB code for MCMC-based coherence-optimal sampling is available \citep{Hampton2015b, Hampton2015}.
However, MCMC-based coherence-optimal sampling can be very slow for high-dimensional input.
An alternative is rejection-based coherence-optimal sampling.
Here, samples $\ve x_\text{cand}$ are generated from a proposal distribution $f_\text{prop}$, which has the property that there is a $\gamma \in \R$ such that $\gamma f_\text{prop}(\ve x) \geq \tilde f_{\ve X}(\ve x)$ for all $\ve x \in \cd$. 
Uniform random numbers $u \sim_\text{i.i.d.} \cu([0,1])$ are generated. 
A proposed point $\ve x_\text{cand}$ is accepted if $u \leq \frac{\tilde f(\ve x_\text{cand})}{\gamma f_\text{prop}(\ve x_\text{cand})}$.
This is the implementation used in this benchmark. 
\changed{We use a product proposal density whose marginals are determined by the input marginals, the dimension of the problem, and the degree of the expansion: we choose uniform proposal marginals for uniform input marginals. For Gaussian input marginals, we use Gaussian proposal marginals if $d \geq p$; otherwise, we use the corresponding asymptotic distribution. As usual, lognormal input is mapped to Gaussian random variables before sampling \citep{BlatmanJCP2011}. }

Note that for Gaussian input, coherence-optimal and asymptotic sampling have a significantly larger spread than input sampling, as can be seen from Figure~\ref{fig:candidateSamplingVisGaussian}. Their support is the ball of radius $r = \sqrt{2}\sqrt{2p + 2}$ (as implemented in \citep{COHOPT}).
This can potentially cause problems in engineering applications, for which simulations may be less accurate when the input parameters are far from typical operating conditions.

\subsubsection{Christoffel sparse approximation}
A similar weighted sampling scheme is \changed{Christoffel sparse approximation \citep{Narayan2017, Jakeman2017, Cohen2017}}.
Those authors propose to use the weight function 
\begin{equation}
w(\ve x) = \left( \frac{1}{|\ca|} \sum_{\alp \in \ca} |\psi_\alp(\ve x)|^2 \right)^{-\frac12}
\end{equation}
which leads to a modified basis that has pointwise minimal average squared basis magnitude (compare to \eqref{eq:weight_cohopt} with $G = B$). This quantity \eqref{eq:CSA_coherence} is a measure similar to coherence \eqref{eq:coherence} and is used by \citet{Hampton2015b} and \citet{Cohen2017} together with the induced probability measure to obtain convergence results for \changed{weighted} least-squares regression.
\citet{Narayan2017} and \citet{Jakeman2017} choose as probability distribution the so-called \textit{weighted pluripotential equilibrium measure} (possibly degree-dependent), which asymptotically coincides with $\tilde f(\ve x) = c^2 w(\ve x)^2 f(\ve x)$ when the total degree of the truncated basis $p \to \infty$. However, the modified basis is not orthonormal with respect to this measure, which leads to weaker theoretical recovery results. Theoretical results are available only for the univariate case. In numerical examples, the method performs well for low-dimensional high-degree cases and often very similarly to asymptotic sampling. In high dimensions, it performs worse than input sampling (i.e., MC). It has not been compared to coherence-optimal sampling.

\subsection{Choosing points according to an optimality criterion from a candidate set}
The following methods aim to improve the properties of the regression matrix by choosing the ``best'' design points from a large candidate set. The methods differ in the criterion defining what are the ``best'' points.
Most presented algorithms are greedy or heuristic and are actually only able to find a suboptimal design (local optimum).
The choice of the candidate set obviously influences the quality of the resulting design. 
In the literature, candidate sets were sampled from MC \citep{Diaz2018}, LHS \citep{FajraouiMarelli2017}, coherence-optimal sampling \citep{Diaz2018, Alemazkoor2018}, or Christoffel sparse approximation \citep{Shin2016b}. In the case of coherence-optimal sampling or Christoffel sparse approximation, the resulting optimized sample inherits the weights. It also often preserves the spread of the candidate set, as can be seen in Figure~\ref{fig:candidateSamplingVisGaussian}.

\subsubsection{D-optimal sampling}
D-optimal design of experiments \citep{Kiefer1959, Dykstra1971} aims at maximizing the determinant of the so-called \textit{information matrix} $\displaystyle\frac1N \ve\Psi^T \ve\Psi \in \R^{P \times P}$. The D-value is defined as
\begin{equation}
D(\ve\Psi) = \det( \ve\Psi^T \ve\Psi).
\end{equation}
Sometimes the determinant of the inverse information matrix is minimized \citep{Nguyen1992}, or the $P$th root is taken for normalization purposes \citep{Diaz2018}.
The maximization of this determinant is connected to the minimization of the variance of the PCE coefficient estimate \citep{Nguyen1992,Zein2013}.
Note that $D(\ve\Psi) = 0$ if $N < P$.

There exists a large selection of methods for constructing D-optimal experimental designs. For an overview of methods for constructing designs following alphabetic optimality criteria (such as A-, D-, or E-optimality), see \cite[Section 4.5]{Hadigol2018}.

Here, we only discuss D-optimal sampling based on rank-revealing QR decomposition (RRQR) \citep{Diaz2018}, since this is the technique used in our benchmark.
We decided to use RRQR-based D-optimal sampling because it can be used even in the case $N < P$ when other D-optimal methods fail due to singularity of the information matrix. 
Note that RRQR is not guaranteed to find a design with maximal D-value but only a local optimum \cite[Section 3.4]{Diaz2018}.

Let $\ve\Psi_\text{cand} \in \R^{M \times P}$ be the regression matrix evaluated at a set of $M$ candidate points.
The goal is to select $N \leq M$ points from this candidate set with the property that the D-value of the resulting regression matrix $\ve\Psi \in \R^{N \times P}$ is as large as possible.
Since in the case of sparse PCE often $N < P$, which leads to $D(\ve\Psi) = 0$, another strategy is necessary.
The RRQR decomposition, also known as pivoted QR, aims at permuting the columns of the original matrix in a way that ensures the R-matrix of the associated QR decomposition is as well-behaved\footnote{$\ve R =\begin{pmatrix} A_k & B_k \\ 0 & C_k \end{pmatrix} $ where $A_k$ is well-conditioned and $\norme{C_k}{2}$ is small} as possible. 
This is useful for inexpensively determining the numerical rank of a matrix \citep{Hong1992, Gu1996}.
RRQR has a strong connection to SVD and to the selection of submatrices of maximal determinant \citep{Hong1992}. 
\citet{Gu1996} propose a pivoted QR decomposition where pivots are chosen to maximize the determinant of the resulting quadratic submatrix of R. The exchange of rows is based on a formula relating the determinant of the quadratic submatrix of R before and after the row exchange by a simple factor \cite[Lemma 3.1]{Gu1996}. 
This algorithm can be used together with SVD to perform subset selection to construct an initial experimental design from a large set of candidate samples \citep{Seshadri2017, Diaz2018}. Here, first an SVD of the matrix $\ve\Psi_\text{cand}^T$ is computed. Then RRQR is applied to the transpose of the matrix consisting of the first $N$ right singular vectors. The resulting permutation matrix is used to determine the points to be chosen from the candidate set.

\subsubsection{Quasi-optimal sampling based on the S-value}
\label{sec:Svalue}
Here, the idea is to select samples from a pool of candidate points so that the PCE coefficients obtained using the selected set are as close as possible to the coefficients that would be obtained if the whole set of candidate points was used \citep{Shin2016a}.
Under the assumption that the columns of the matrix $\ve\Psi_\text{cand}$ are mutually orthogonal, the S-value is defined by%
\begin{equation}
S(\ve\Psi) = \left( \frac{\sqrt{\det \ve\Psi^T \ve\Psi}}{\prod_{i=1}^P \norme{\Psi_i}{2}} \right)^{\frac1P},
\end{equation}
where $\Psi_i$ denotes the $i$th column of the regression matrix $\ve\Psi$. Its maximization has the (heuristic) effect of maximizing the column orthogonality of the regression matrix while at the same time maximizing the determinant of the information matrix \citep{Shin2016a}. 
It holds that $S(\ve\Psi) \in [0,1]$ due to Hadamard's inequality. If $N < P$, $S(\ve\Psi) = 0$. If $N \geq P$, $S(\ve\Psi) = 1$ if and only if the columns of $\ve\Psi$ are mutually orthogonal.
There exists an update formula for the S-value when the regression matrix is augmented by one row, which thus avoids the repeated calculation of determinants.

\citet{Shin2016a} suggest a greedy algorithm that in every iteration augments the current matrix by an additional row which maximizes the S-value of the resulting matrix among all candidate rows. When the current number of rows in the matrix $\ve \Psi$ is smaller than the number of columns, the procedure can be adapted to avoid $S(\ve\Psi) = 0$.
We do not include it in our benchmark because it is not well suited for situations where there are more basis polynomials than design points, which is the case in sparse PCE without experimental design enrichment. 
However, in a sequential enrichment context \citep{FajraouiMarelli2017} and for least-squares regression \citep{Shin2016b}, this algorithm performs well.

\subsubsection{Near-optimal sampling}
\label{sec:near-opt}
The coherence parameter \eqref{eq:coherence_app} gives a bound on the recovery rate, but it is not the only criterion that has been studied with respect to recovery accuracy. 
Two other matrix properties related to recovery accuracy are mutual coherence and average cross-correlation.
Both of them consider the correlation between normalized columns of the regression matrix, i.e., their scalar product. 
They are scalar measures of how ``orthonormal'' the columns of a rectangular matrix $\ve\Psi \in \R^{N \times P}$ with $N < P$ are. 
The heuristic idea is that columns should point in as different directions as possible, so that the multiplication with sparse coefficient vectors, which results in a linear combination of a subset of the columns, is ``as unique as possible''. This facilitates the recovery of the true sparse solution (assuming that it exists).

The \textit{mutual coherence} is defined by
\begin{equation}
\mu(\ve\Psi) = \max_{i\neq j} \frac{|\Psi_i^T \Psi_j|}{\norme{\Psi_i}{2}\norme{\Psi_j}{2}},
\end{equation}
where $\Psi_i$ denotes the $i$th column of the regression matrix $\ve\Psi \in \R^{N \times P}$.
The mutual coherence is the worst-case cross-correlation between any two columns of the matrix. It is zero for orthonormal matrices and positive for $N < P$. 

The \textit{average \changed{(squared)} cross-correlation} is defined by
\begin{equation}
\gamma(\ve\Psi) = \frac{1}{\changedmath{P(P-1)}} \norme{\mathbb{1}_\changedmath{P} - \tilde{\ve\Psi}^T\tilde{\ve\Psi}}{F}^2 = \frac{1}{\changedmath{P(P-1)}} \sum_{i\neq j} \frac{|\Psi_i^T \Psi_j|^{\changedmath{2}}}{\norme{\Psi_i}{2}^{\changedmath{2}}\norme{\Psi_j}{2}^{\changedmath{2}}}
\end{equation}
where $\tilde{\ve\Psi}$ is the column-normalized version of $\ve\Psi$, and $\Psi_i$ denotes the $i$th column of the regression matrix. The norm is the Frobenius-norm, taking the sum of squares of all matrix entries, and the factor $\changedmath{P(P-1)}$ is the number of column pairs.

\citet{Alemazkoor2018} suggest simultaneously optimizing mutual coherence and average cross-correlation by using the greedy procedure described in Algorithm \ref{alg:nearoptimalsampling} below: In each iteration, the current regression matrix is augmented by one row. This row corresponds to that point $\ve x_j$ from the large pool of candidate points which minimizes the (normalized) distance of $(\mu_j', \gamma_j') \in \R^2$ to the ``utopia point'' $(\min(\ve\mu'), \min(\ve\gamma'))$ among all candidate points.

\begin{algorithm}
	\caption{Near-optimal sampling \citep{Alemazkoor2018}.}
	\label{alg:nearoptimalsampling}
	\begin{algorithmic}[1] 
		\State Sample a large number $M$ of candidate points from the coherence-optimal distribution and compute candidate rows arranged in a matrix $\ve\Psi_\text{cand}$
		\State Initialize $\ve \Psi_\text{opt(1)}$ to be a random row from $\ve\Psi_\text{cand}$ 
		\For{$i=2\ldots N$}
		\For{$j=1\ldots M$}
		\State $\ve\Psi_\text{temp} = $ row-concatenate$(\ve \Psi_\text{opt(i-1)}, \Psi^{(j)}_\text{cand})$
		\State $\mu_j' = \mu(\ve\Psi_\text{temp})$
		and $\gamma_j' = \gamma(\ve\Psi_\text{temp})$ 
		\EndFor
		\State $\ve{\mu}' = (\mu_1' \enum \mu_{M}')$ and $\ve{\gamma}' = (\gamma_1' \enum \gamma_{M}')$
		\State $j^* = \arg\min_j \left( \frac{\mu_j' - \min(\ve{\mu}')}{\max(\ve\mu') - \min(\ve\mu')} \right)^2 + \left( \frac{\gamma_j' - \min(\ve{\gamma}')}{\max(\ve\gamma') - \min(\ve\gamma')} \right)^2$
		\State $\ve\Psi_\text{opt(i)} = $ row-concatenate$(\ve\Psi_\text{opt(i-1)}, \Psi^{(j^*)}_\text{cand})$
		\EndFor 
	\end{algorithmic}
\end{algorithm}

The algorithm is called near-optimal because it is a greedy algorithm, finding only a local optimum, and because optimized mutual coherence and average cross-correlation are only hinting at, but not guaranteeing, good recovery accuracy \citep{Alemazkoor2018}.
Its computational complexity is \changed{$\co(NMP^2)$}, where $N$ is the size of the final experimental design, $M$ is the number of candidate samples (chosen to be, e.g., proportional to $P$ \citep{Diaz2018}), and $P$ is the number of regressors. 
This makes the algorithm prohibitively expensive in the case of large bases ($P$ in the order of thousands), which is why we do not use it for some of the benchmark examples.

\subsection{Illustration of sampling schemes}
\begin{figure}[h!]
	\centering
	
	\subfloat[][MC]{
		\includegraphics[width=.315\textwidth]{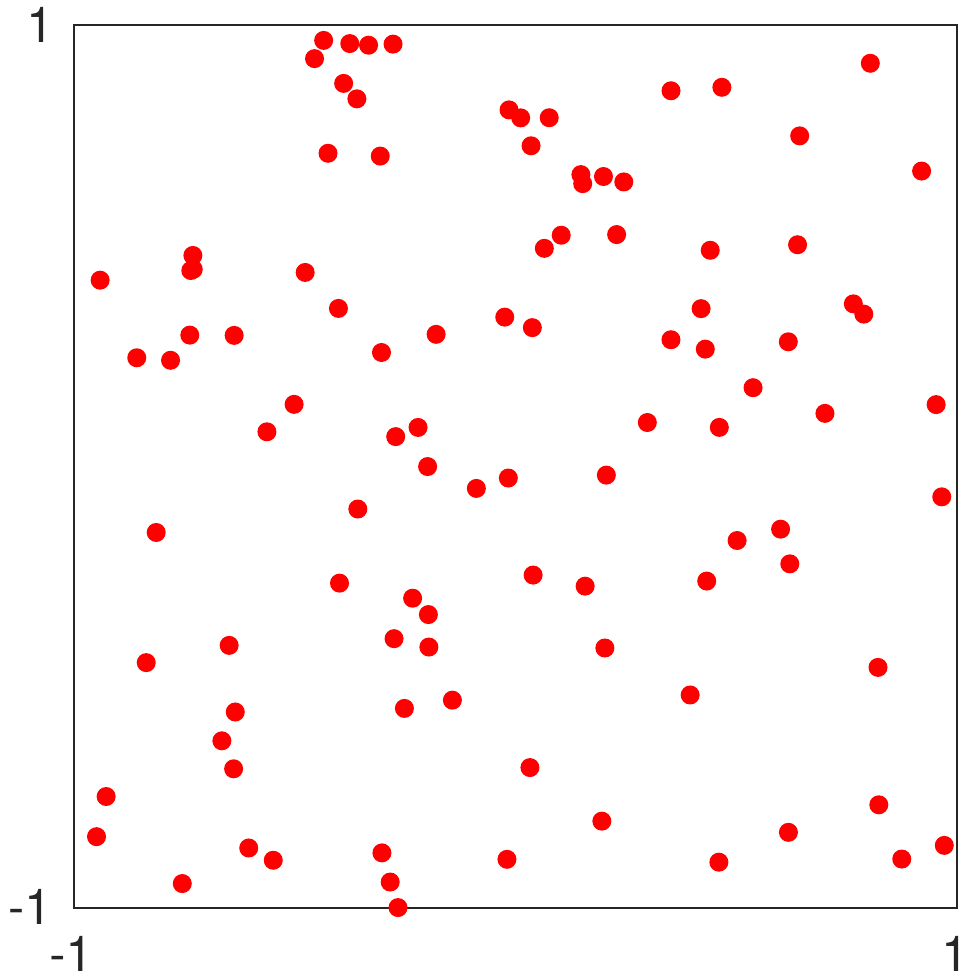}
	}
	\subfloat[][D-opt(MC)]{
		\includegraphics[width=.315\textwidth]{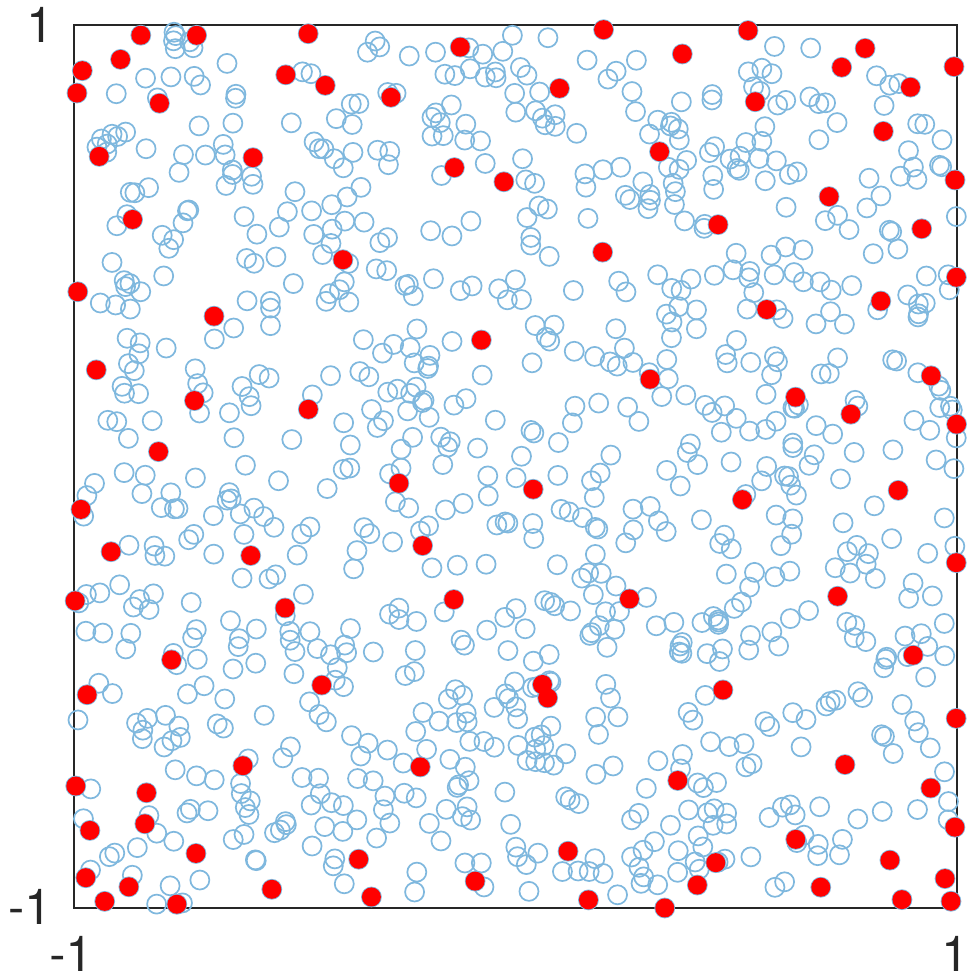}
	}
	\subfloat[][asymptotic]{
		\includegraphics[width=.315\textwidth]{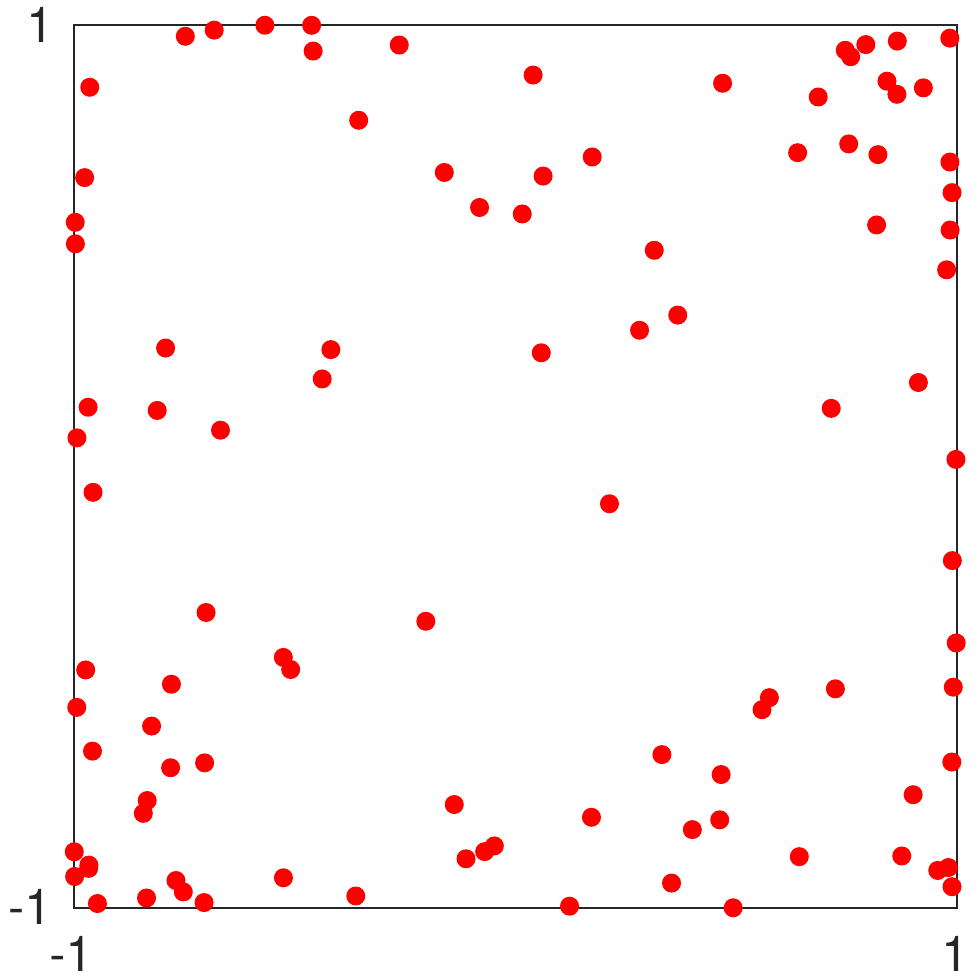}
	}
	\\
	%	\vspace{.2cm}
	\subfloat[][LHS]{
		\includegraphics[width=.315\textwidth]{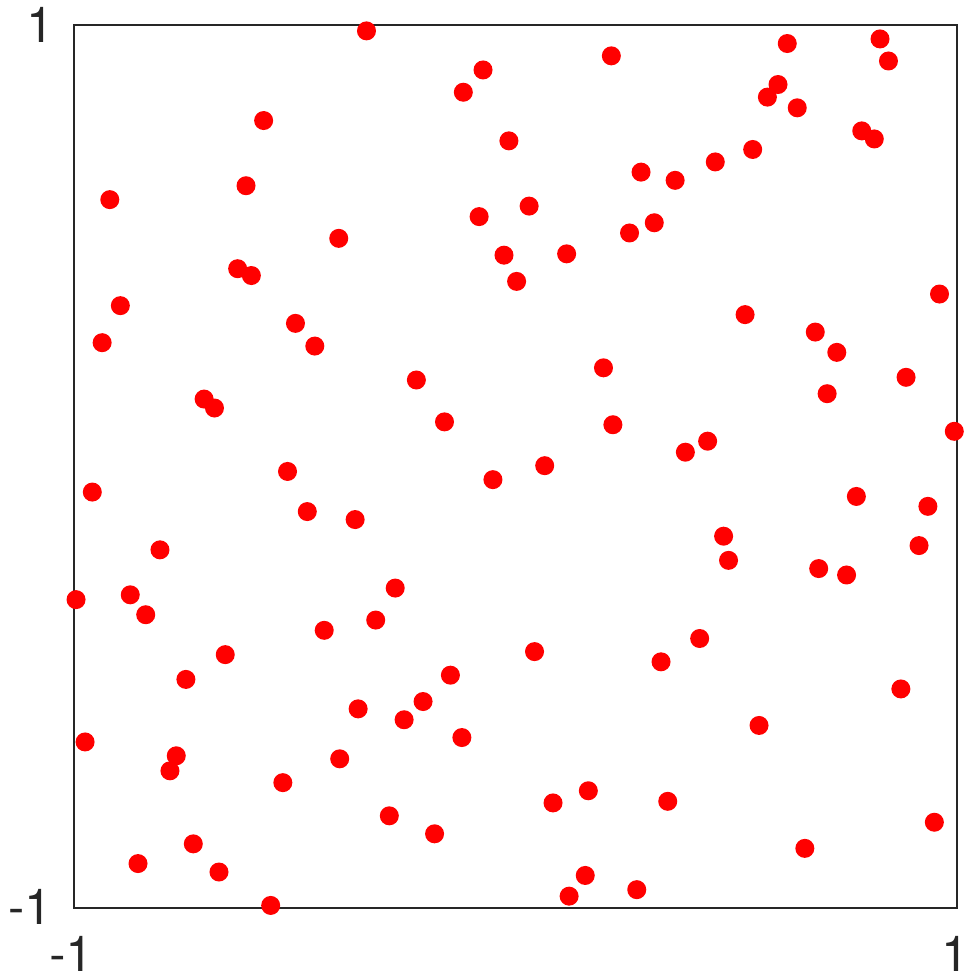}
	}
	\subfloat[][D-opt(LHS)]{
		\includegraphics[width=.315\textwidth]{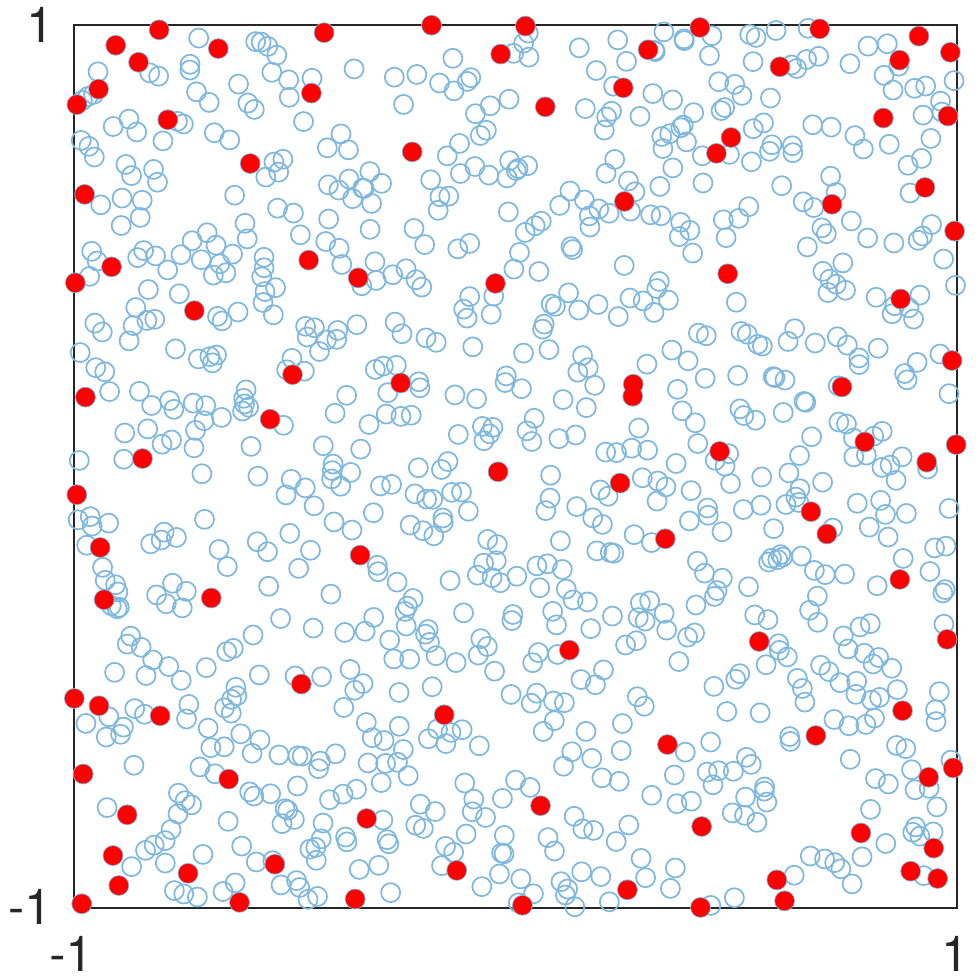}
	}
	\subfloat[][near-opt(LHS)]{
		\includegraphics[width=.315\textwidth]{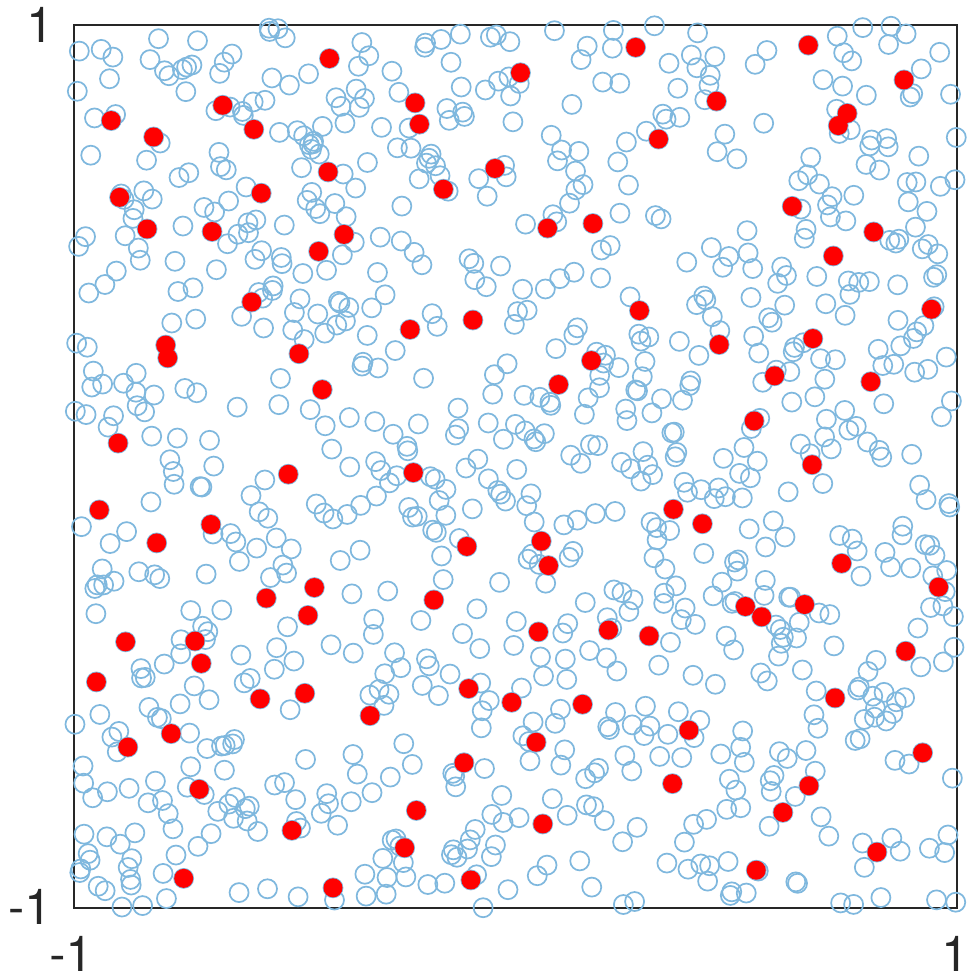}
	}
	\\
	%	\vspace{.2cm}
	\subfloat[][coh-opt]{
		\includegraphics[width=.315\textwidth]{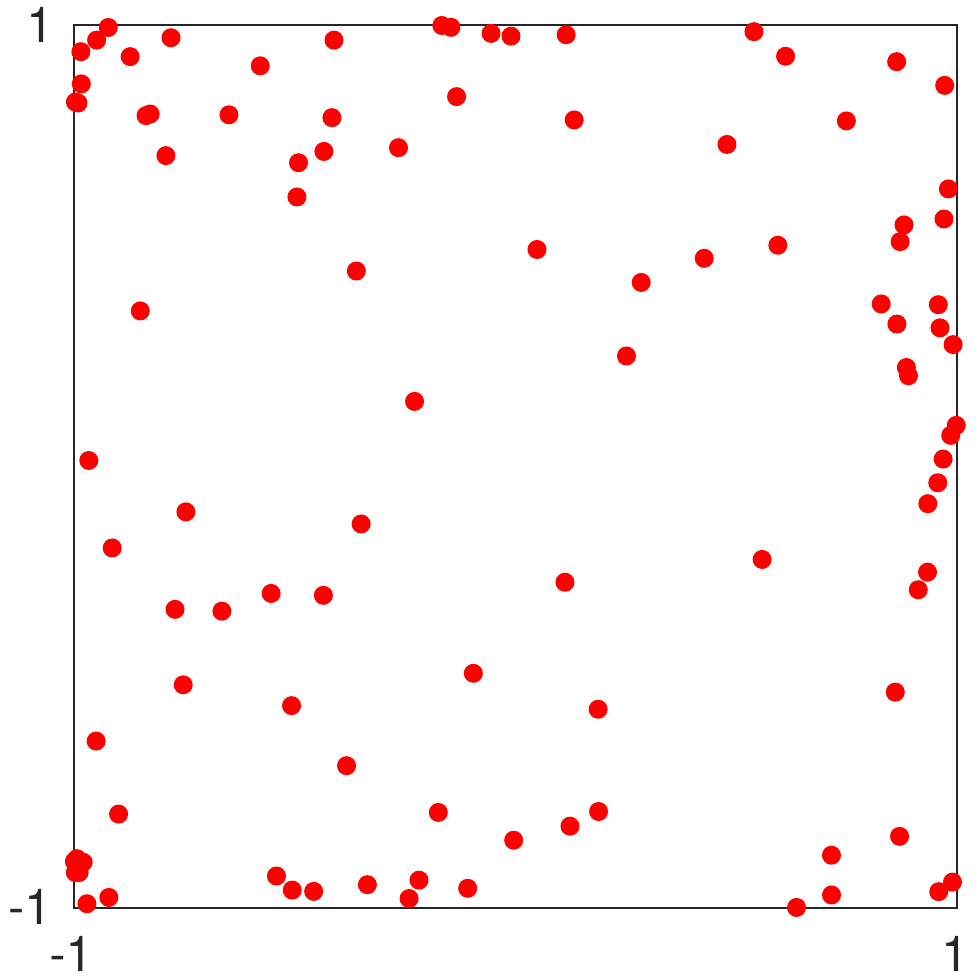}
	}
	\subfloat[][D-opt(coh-opt)]{
		\includegraphics[width=.315\textwidth]{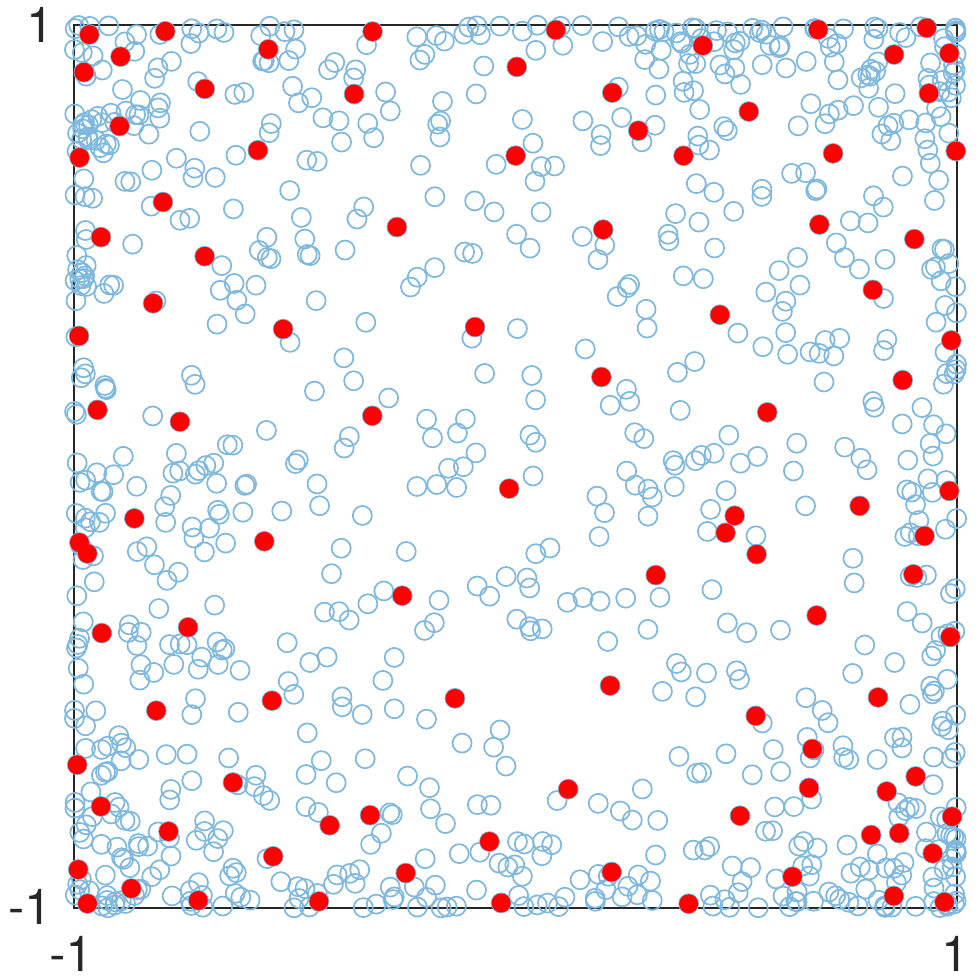}
	}
	\subfloat[][near-opt(coh-opt)]{
		\includegraphics[width=.315\textwidth]{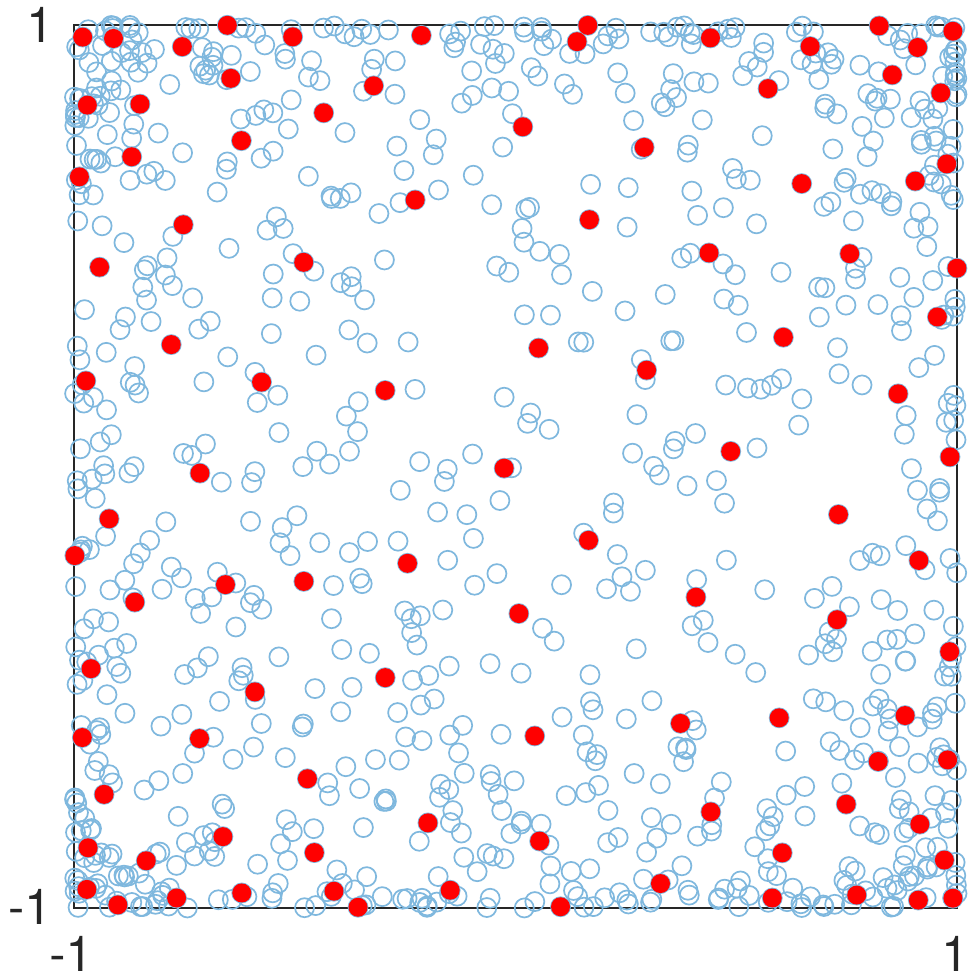}
	}
	\caption{Visualization of experimental designs constructed for uniform input in $[-1,1]^2$ for degree $p=12$. Red filled points denote the chosen experimental design, while blue circles denote the candidate set. Size of the ED: $\changedmath{N} = 100$, size of the candidate set: $M = 1000$.}
	\label{fig:candidateSamplingVisUnif}
\end{figure}

\begin{figure}[h!]
	\centering
	\subfloat[][MC]{
		\includegraphics[width=.315\textwidth]{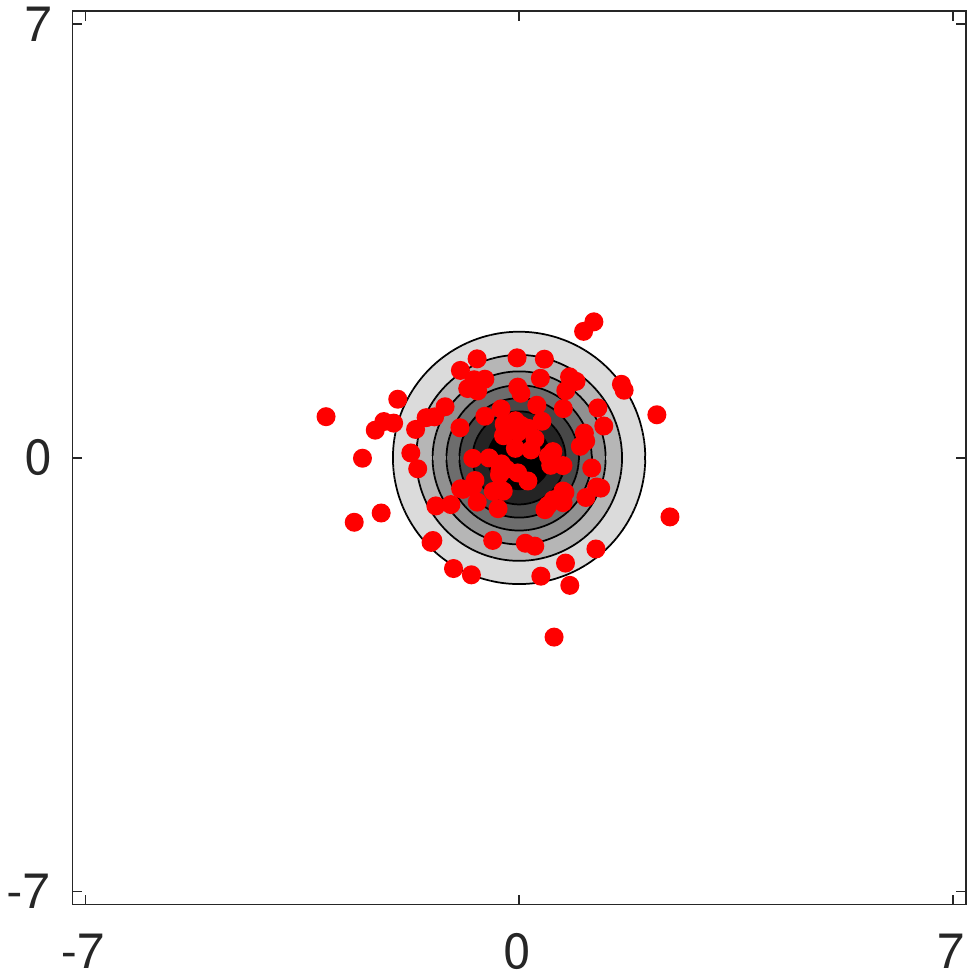}
	}
	\subfloat[][D-opt(MC)]{
		\includegraphics[width=.315\textwidth]{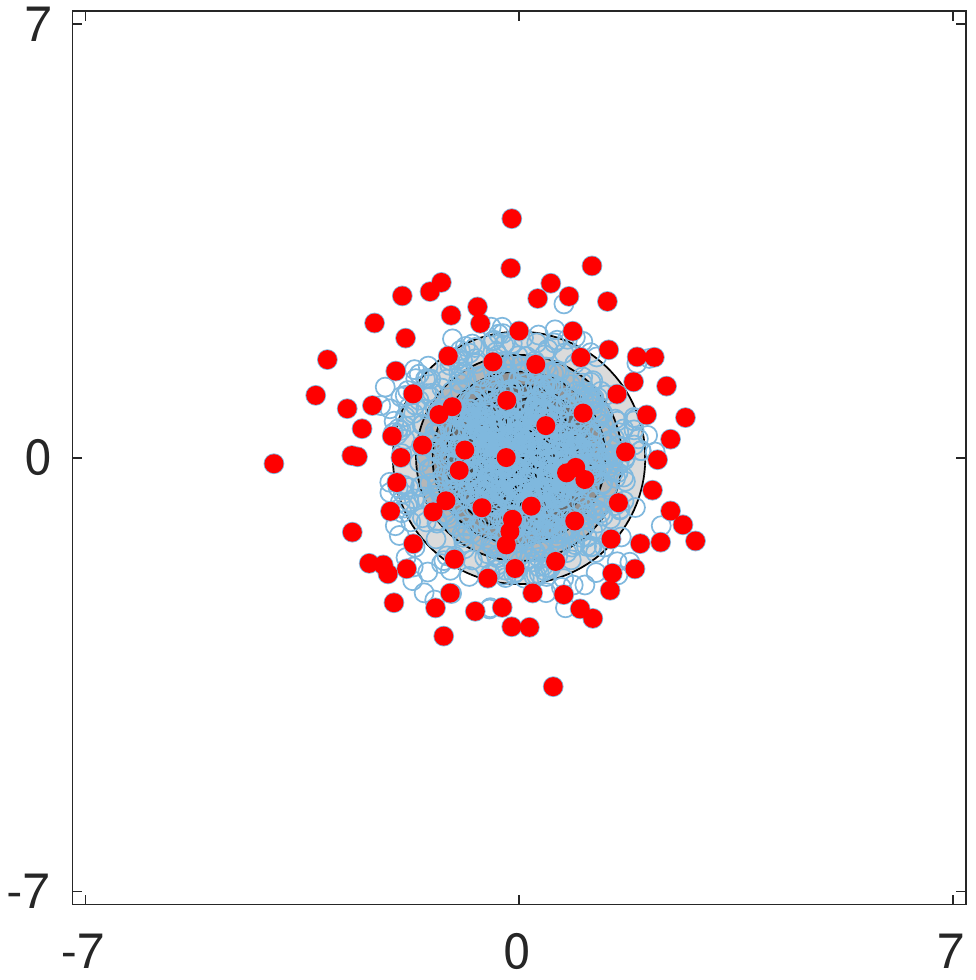}
	}
	\subfloat[][asymptotic]{
		\includegraphics[width=.315\textwidth]{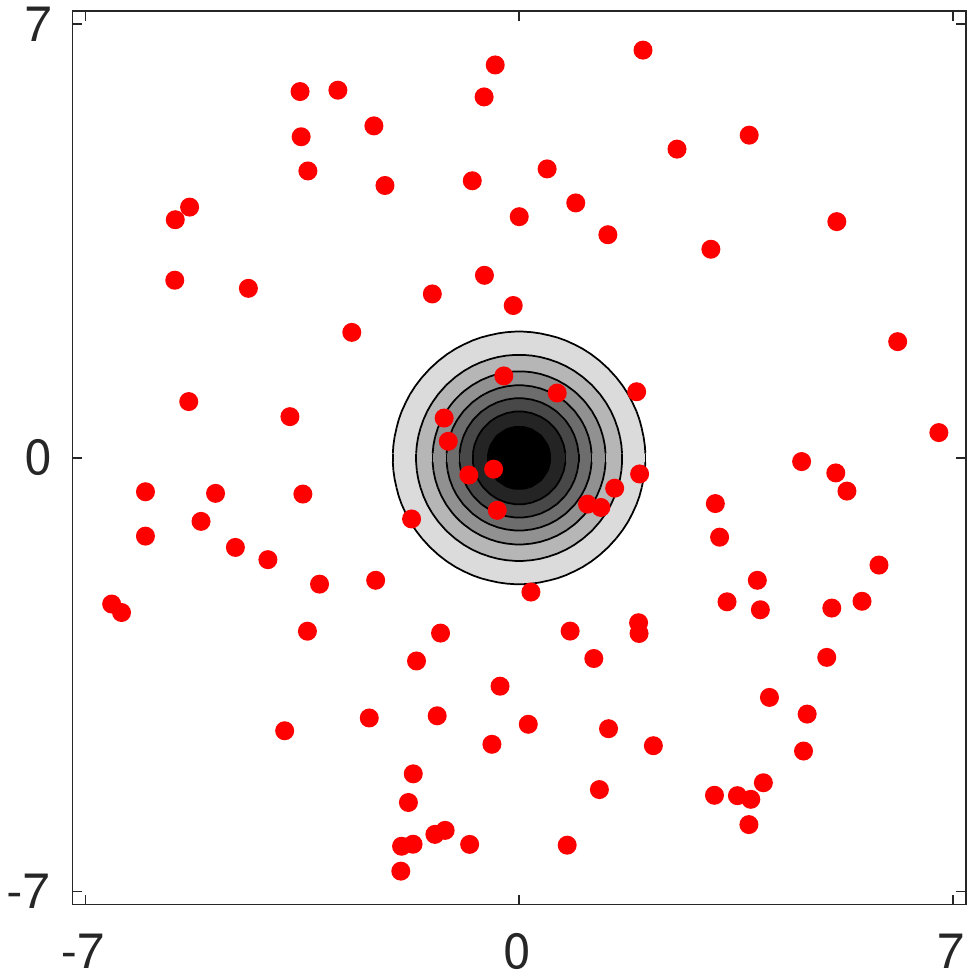}
	}
	\\
	%	\vspace{.2cm}
	\subfloat[][LHS]{
		\includegraphics[width=.315\textwidth]{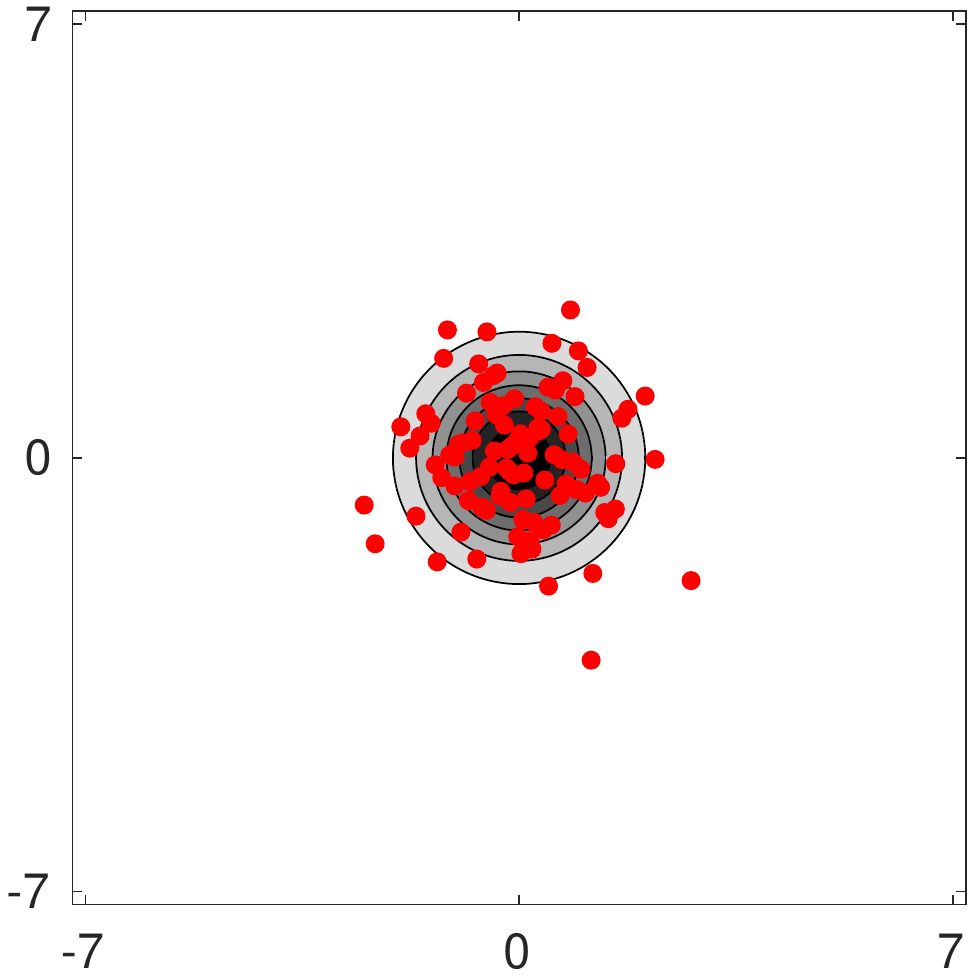}
	}
	\subfloat[][D-opt(LHS)]{
		\includegraphics[width=.315\textwidth]{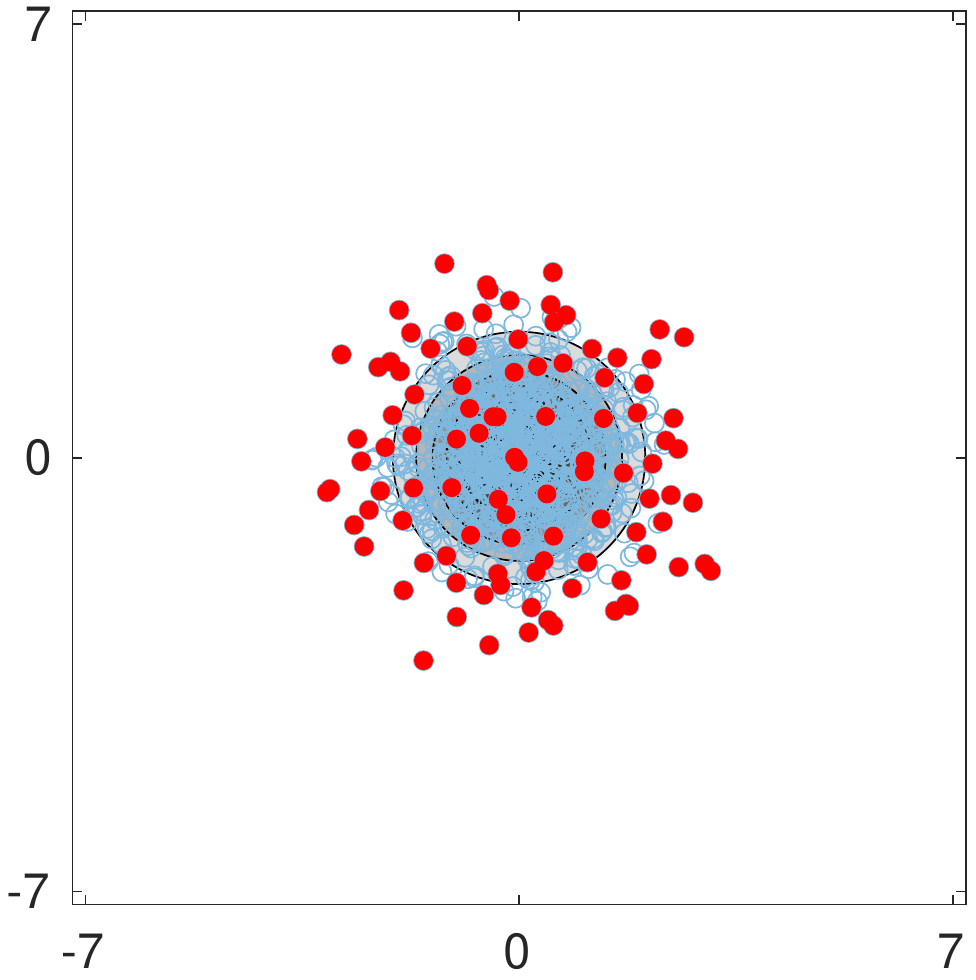}
	}
	\subfloat[][near-opt(LHS)]{
		\includegraphics[width=.315\textwidth]{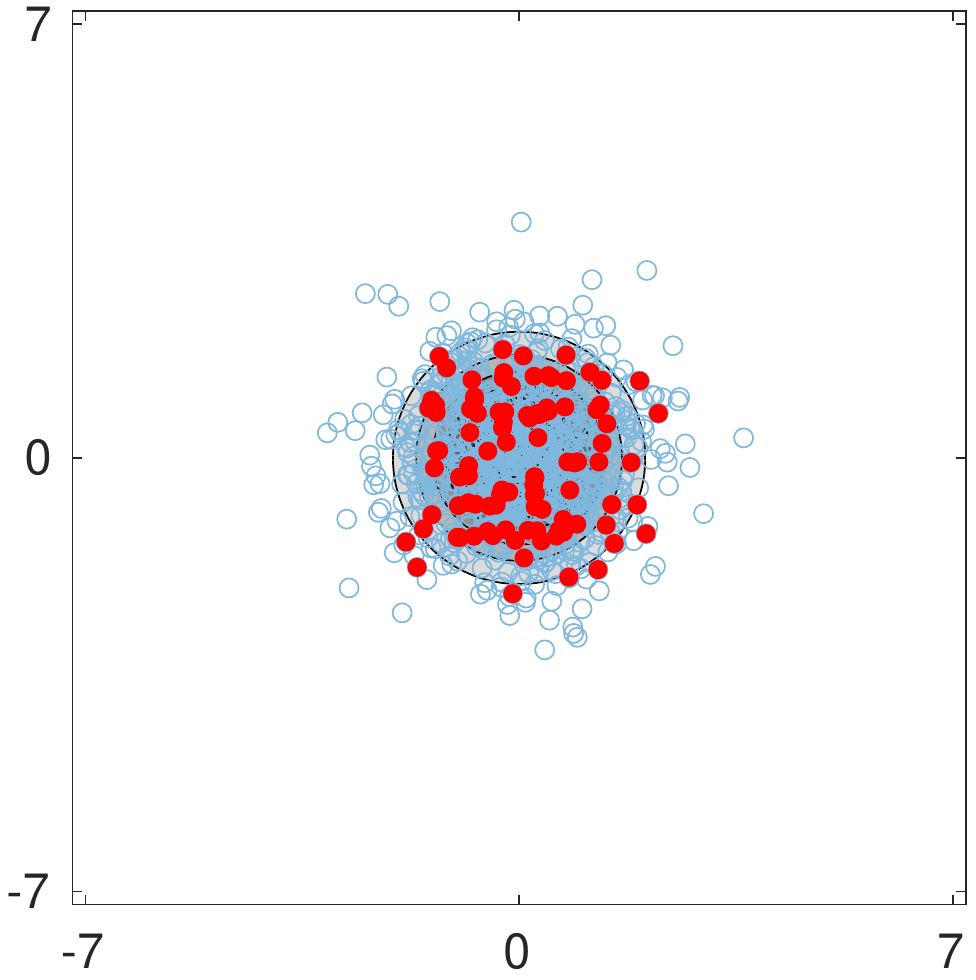}
	}
	\\
	%	\vspace{.2cm}
	\subfloat[][coh-opt]{
		\includegraphics[width=.315\textwidth]{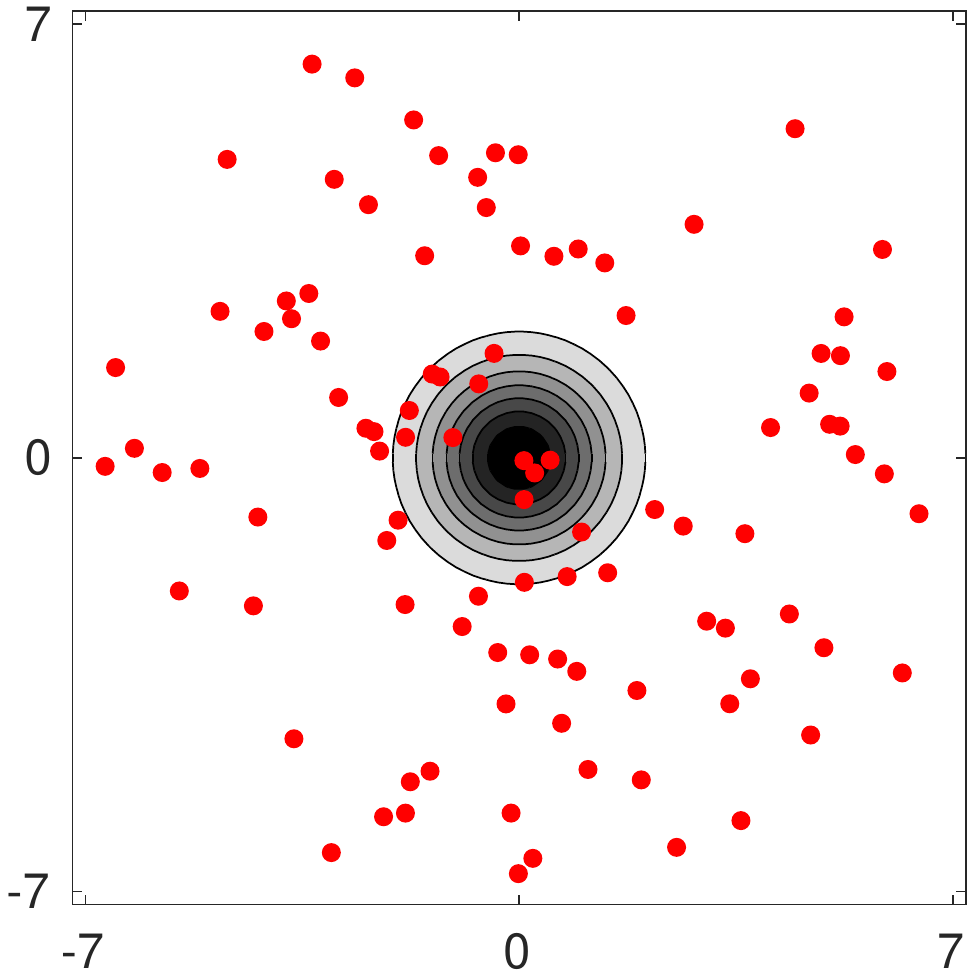}
	}
	\subfloat[][D-opt(coh-opt)]{
		\includegraphics[width=.315\textwidth]{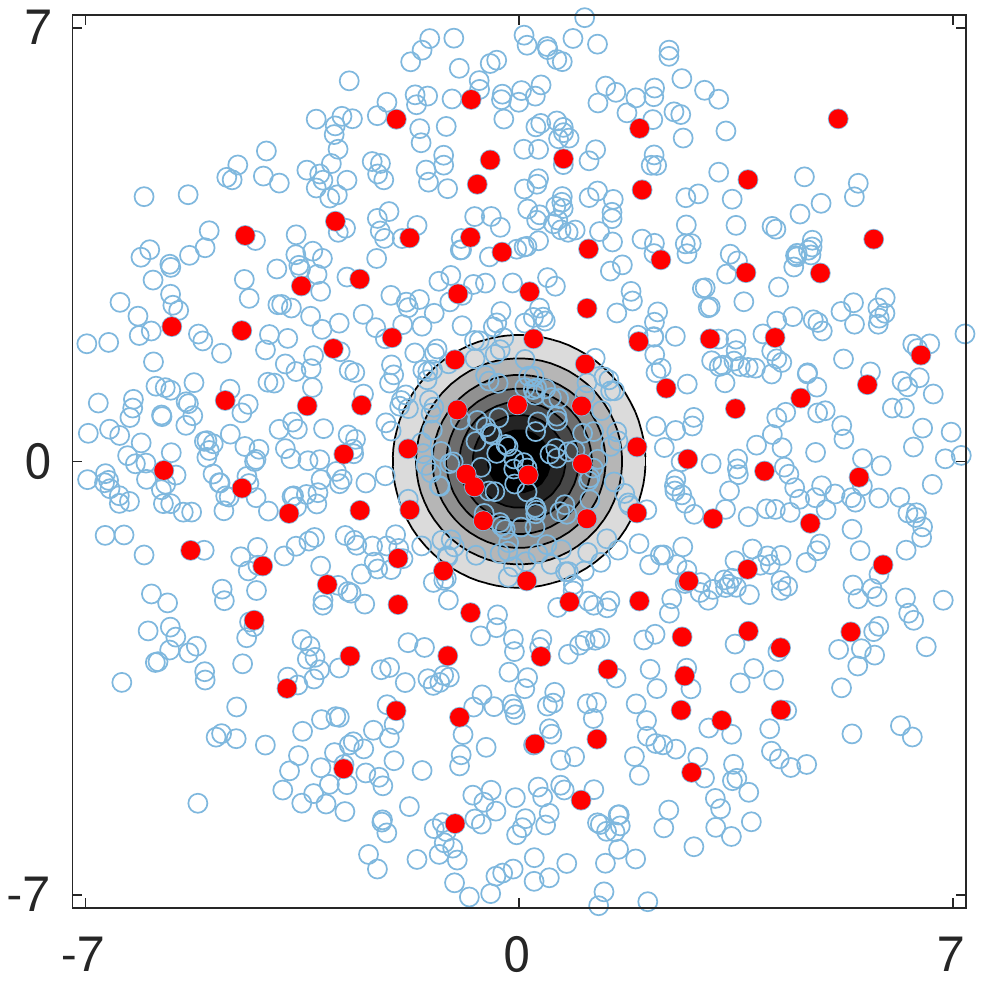}
	}
	\subfloat[][near-opt(coh-opt)]{
		\includegraphics[width=.315\textwidth]{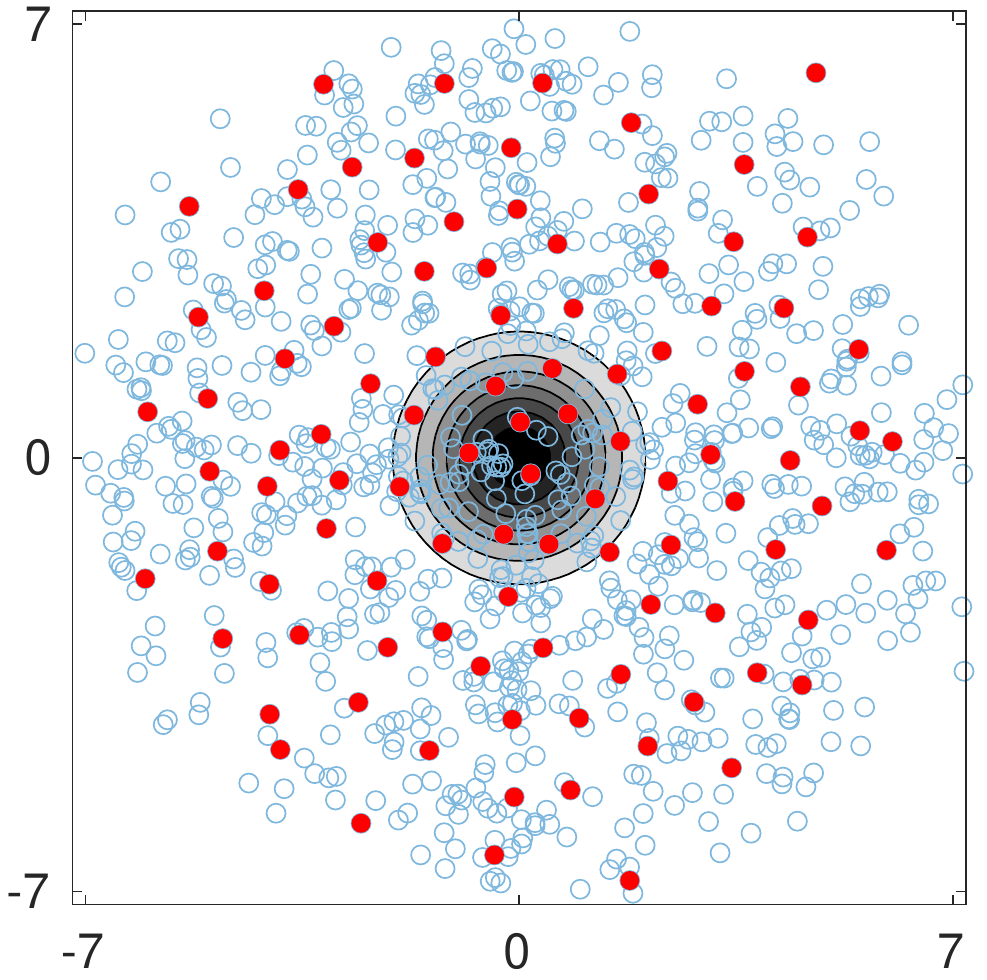}
	}
	\caption{Visualization of the experimental design constructed for standard Gaussian input in $d=2$ dimensions for degree $p=12$. The gray surface plot illustrates the Gaussian probability density function. Red filled points denote the chosen experimental design, while blue circles denote the candidate set. Size of the ED: $N = 100$; size of the candidate set: $M = 1000$. The support of the asymptotic and the coherence-optimal distribution is the ball of radius $r = \sqrt{2}\sqrt{2p+2} \approx 7.2$. Note that engineering models may be less accurate in regions where the input distribution has negligible mass.}
	\label{fig:candidateSamplingVisGaussian}
\end{figure}

In Figures~\ref{fig:candidateSamplingVisUnif} and \ref{fig:candidateSamplingVisGaussian}, we show illustrations of experimental designs in $d=2$ dimensions with $N = 100$ and $p = 12$ for selected sampling techniques. The candidate set has a size of $M = 1000$.
Figure~\ref{fig:candidateSamplingVisUnif} presents experimental designs for uniform input in the interval $[-1,1]$, while Figure~\ref{fig:candidateSamplingVisGaussian} presents experimental designs for standard Gaussian input.

Note that in the standard Gaussian case, the asymptotic distribution, the coherence-optimal distribution, and the matrix-optimal distributions based on a coherence-optimal candidate set all have a very large spread that grows with the total degree of the basis. For degree $p = 12$, some points are seven standard deviations away from the mean. Engineering models are typically calibrated only for a certain region of the input domain corresponding to nonnegligible probability, and they may be less accurate (or even fail) outside of this region.

\clearpage

\section{Details on sparse regression solvers}
\label{app:solvers}

In this appendix, we describe the sparse solvers used in our benchmark in more detail: LARS, OMP, subspace pursuit, SPGL1 and FastLaplace (BCS).
In addition, we present an overview of greedy stepwise regression solvers for sparse PCE.

There exist various formulations for the sparse regression problem. 
The typical form minimizes the $\ell^2$-norm of the empirical error under an additional constraint that is designed to enforce sparsity.

Sparsity is measured by the number of nonzero entries in a vector, formally denoted by $\norme{\ve c}{0} = \sum_{i} \mathbf{1}_{\{c_i \neq 0 \}}$ (even though this expression is not a norm).
This results in the sparse regression problem
\begin{equation}
\hat{\ve c} = \arg\min_{\ve c} \norme{\ve\Psi \ve c - \ve y}{2}^2  + \lambda \norme{\ve c}{0}
\label{eq:l0-min}
\end{equation}
called $\ell^0$-minimization. The only way to solve this problem exactly is by a combinatorial search through all possible nonzero patterns for $\ve c$, which is infeasible for large problem sizes. 

The convex relaxation of this problem is $\ell^1$-minimization, where $\norme{\ve c}{0}$ is replaced by $\norme{\ve c}{1} = \sum_{i} |c_i|$. 
There are several equivalent formulations of the relaxed problem, namely
\begin{align}
	&\hat{\ve c} = \arg\min_{\ve c} \norme{\ve\Psi \ve c - \ve y}{2}^2  + \lambda \norme{\ve c}{1}
	\label{eq:Lagrangian}
	\\
	&\hat{\ve c} = \arg\min_{\ve c} \norme{\ve c}{1} \text{ s.t. } \norme{\ve\Psi \ve c - \ve y}{2} \leq \sigma
	\label{eq:BPDN} 
	\\
	&\hat{\ve c} = \arg\min_{\ve c} \norme{\ve\Psi \ve c - \ve y}{2} \text{ s.t. } \norme{\ve c}{1} \leq \tau
	\label{eq:LASSO}
\end{align}
called Lagrangian formulation, basis pursuit denoising (BPDN), and least absolute shrinkage and selection operator (LASSO), respectively.
It has been shown that under certain conditions, the solutions to $\ell^0$-minimization and $\ell^1$-minimization coincide \citep{Bruckstein2009}.
However contrary to \eqref{eq:l0-min}, formulations \changed{\eqref{eq:Lagrangian}--\eqref{eq:LASSO}} are convex problems and allow a numerical solution with considerably smaller cost.
The three formulations \eqref{eq:Lagrangian}, \eqref{eq:BPDN}, and \eqref{eq:LASSO} are equivalent in the sense that if $\hat{\ve c}$ is solution to one of the formulations, there exists a value of constraint parameter $\sigma, \tau$, or $\lambda$ so that $\hat{\ve c}$ is also a solution to the other formulations. However, the relationship between the parameters $\sigma, \tau, \lambda$ that makes the problems equivalent depends on $\ve\Psi$ and $\ve y$ and is not known in advance \citep{Vandenberg2008}.

There exist other sparsity-enforcing formulations, such as $\ell^p$-norms \citep{Bruckstein2009}, $\ell^1 - \ell^2$-mini\-mi\-za\-tion \citep{Yin2015}, 
or elastic net \citep{Tarakanov2019}. One example that we will describe is Bayesian compressive sensing, where a sparsity-enforcing prior is used for the coefficients of the PCE, resulting in a formulation that is equivalent to a sparse regression problem with a different kind of sparsity constraint, e.g., one related to the Student-t distribution \citep{Tipping2001}.

In the following descriptions of the algorithms, $\ca \subset \N^d$ with various sub- or superscripts denotes a set of multi-indices, which by definition of PCE can be identified with a set of basis polynomials.
With the notation \changed{of} Section \ref{sec:PCE}, $\ve y \in \R^N$ denotes the vector of model responses, $\ve \Psi \in \R^{N \times P}$ denotes the regression matrix of basis polynomials evaluated at the $N$ experimental design points, and $\ve c \in \R^P$ denotes the coefficients of a PCE. The residual is defined by $\ve r = \ve y - \ve \Psi \ve c$, so that the norm of the residual is the empirical error. 

The polynomials used for building the PCE are sometimes also called basis functions, regressors, or predictors. 
A sparse PCE is a PCE for which only some of the basis functions have nonzero coefficients: these basis functions are called \textit{active}. 
We assume the regressors are normalized, so that the correlation between two regressors is equivalent to their inner product (and to the cosine of the angle between them).

% ====================================================================================================

\subsection{Orthogonal matching pursuit (OMP)}
\label{sec:OMP}
Orthogonal matching pursuit (OMP), also called forward stepwise regression, is a classical greedy technique for finding approximate solutions to the $\ell^0$-minimization problem \eqref{eq:l0-min} \citep{Pati1993,Bruckstein2009}. Despite its being a heuristic method, under certain assumptions there are theoretical guarantees for the solutions returned by OMP \citep{Tropp2007, Bruckstein2009}.
OMP is an iterative algorithm that starts out with an empty model and adds the regressors one by one to the set of active regressors. In each iteration, OMP selects the regressor that is most correlated with the current residual, adds it to the set of active regressors, and then updates the coefficients of all active regressors to make sure the new residual is orthogonal to all of them and has smallest possible norm.
The updating of the coefficients can be done through an update formula \citep{berchier2015OMP} or by computing the least-squares solution to the system of equations involving only the active regressors \citep{UQdoc_13_104}. 

The technique is presented in Algorithm \ref{alg:OMP}. The iterations are continued until $\min\{N,P\}$ basis functions are in the active set (then either all polynomials are selected, or there are not enough points in the experimental design to use least-squares anymore). 

\begin{algorithm}
	\caption{Orthogonal matching pursuit (OMP) \citep{Pati1993, Tropp2007, UQdoc_13_104}.}
	\label{alg:OMP}
	\begin{algorithmic}[1]
		\State Given a set of candidate basis functions $\ca_\text{cand}$
		\State Initialize all coefficients to zero: $\ve c_0 = 0$. Set $\ca_0 = \emptyset$
		\State Set the residual vector $\ve r := \ve y$
		\For{$i = 1 \enum m$}
		\Graycomment{$m \leq \min\{N,P\}$}
		\Statex \Graycomment{OMP can be stopped early when the error did not decrease anymore for a while}
		\State Find $\alp^* \in \ca_\text{cand} \setminus \ca_{i-1}$ with maximal correlation with the residual 
		by solving 	
		\[
		\alp^* = \arg \max_{\alp \in \ca_\text{cand} \setminus \ca_{i-1}} |\ve r^T \ve \psi_\alp|
		\]
		\Statex \Graycomment{The entries of vector $\ve \psi_\alp$ are  evaluations of the basis function $\psi_\alp$ at the ED}
		\State $\ca_{i} = \ca_{i-1} \cup \{\alp\}$
		\Graycomment{Current set of active predictors}
		\State Compute the coefficients $\ve c_i$ by least-squares using only the active indices $\ca_i$
		\Statex \Graycomment{This can be done in $\co(iN)$ when maintaining a QR factorization \citep{Tropp2007}}
		\State Update the residual $\ve r = \ve y - \ve\Psi_{\ca_i} \ve c_i$
		\EndFor
	\end{algorithmic}
\end{algorithm}
OMP does not, per se, return a sparse solution. If a desired level of sparsity $K$ is known a priori, the algorithm can be stopped after $K$ iterations. Another possibility is to stop the algorithm as soon as the residual norm is smaller than some error threshold \citep{Bruckstein2009,Doostan2011, Jakeman2015}, where the best error threshold is determined through cross-validation.
A third possibility is to determine the best number of active basis functions through a model selection criterion, e.g.\ the LOO error \citep{UQdoc_13_104}. Since the coefficients are computed by OLS on the active basis, the LOO can be computed cheaply \citep{Chapelle2002, BlatmanJCP2011}. Typically, for an increasing sequence of basis functions the LOO error first decreases (reduction of underfitting), then increases (overfitting). This can be utilized to terminate the algorithm early once the LOO error starts rising (early-stop criterion) \citep{UQdoc_13_104}. 

The computational complexity of OMP is $\co(mNP)$
\citep{Tropp2007, Dai2009},
where $m \leq \min\{N,P\}$ is the number of iterations. 
The computation of the correlations of the current residual with all regressors is $\co(NP)$ and has to be performed $m$ times. The computation of the least-squares solution in step $i$ can be done in $\co(iN)$, e.g., by maintaining a QR factorization of the information matrix \citep{Tropp2007},
or by using Schur's complement to update the information matrix inverse whenever a new regressor is added.

From the authors' experience, OMP often suffers from overfitting and can produce an unreliable LOO error estimate, which can be detrimental in basis-adaptive settings (see also \citep{LuethenIJUQ2021}).

OMP is available in many software packages, among them UQLab \citep{MarelliUQLab2014}.

% ====================================================================================================

\subsection{Least angle regression (\changed{LARS})}
\label{sec:LAR}
Least-angle regression (\changed{LARS, sometimes also abbreviated LAR}) is a greedy technique that finds an approximate solution to the $\ell^1$-minimization problem \citep{Efron2004}. 
It is similar to OMP in that the algorithm starts out with an empty model and adds regressors one by one based on their correlation with the residual. However, unlike OMP, which updates the coefficients using least-squares (making the residual orthogonal to all active regressors in each step), LARS updates the coefficients in such a way that all active regressors have equal correlation with the residual.
LARS can be interpreted as producing a path of solutions to \eqref{eq:LASSO}, corresponding to increasing $\tau$. The coefficients are increased in the equiangular direction until a nonactive regressor has as much correlation with the residual as all the active regressors. This regressor is then added to the set of active regressors and the new equiangular direction is computed. The optimal stepsize between the addition of subsequent regressors can be computed analytically \citep{Efron2004}. 
This algorithm solves \eqref{eq:LASSO} approximately. A slightly modified version of LARS, called LARS-LASSO, removes regressors whenever the sign of their coefficient changes, and it has been proven to solve \eqref{eq:LASSO} (or its noiseless counterpart) exactly under certain conditions \citep{Efron2004, Bruckstein2009}.

\begin{algorithm}
	\caption{Least angle regression (LARS) \citep{Efron2004, BlatmanJCP2011}.}
	\label{alg:originalLAR}
	\begin{algorithmic}[1]
		\State Given a set of candidate basis functions $\ca_\text{cand}$
		\State Initialize all coefficients to zero: $\ve c_0 = 0$. Set $\ca_0 = \emptyset$
		\State Set the residual vector $\ve r := \ve y$
		\For{$i = 1 \enum m$}
		\Graycomment{$m = \min\{N,P\}$}
		\State Find $\alp \in \ca_\text{cand} \setminus \ca_{i-1}$ with maximal correlation with the residual
		\Statex \Graycomment{For $i > 1$, $\alp_i$ is the element "responsible for" $\gamma_i$}
		\State $\ca_{i} = \ca_{i-1} \cup \{\alp\}$
		\Graycomment{Current set of active predictors}
		\State Compute $\ve c_{i}$ 
		\label{line:equiangular}
		\Graycomment{Equiangular direction for all $\alp \in \ca_i$}
		\Statex \Graycomment{$\ve c_1$ is equal to the first selected predictor}
		\Statex \Graycomment{\cite[Eq.\ 2.6]{Efron2004}}
		\State Compute $\gamma_i$
		\Graycomment{Optimal stepsize: using this, there is a new regressor that is as much correlated with $r$ as all regressors in $\ca_i$ are}
		\Statex \Graycomment{\cite[Eq.\ 2.13]{Efron2004}}
		\State Compute the new coefficients $\ve c_{i} = \ve c_{i-1} + \gamma_{i} \ve c_{i}$
		\Graycomment{Move the coefficients jointly into the direction of the least-squares solution until one of the other predictors in $\ca_\text{cand} \setminus \ca_i$ has as much correlation with the residual as the predictors in $\ca_i$ (ensured by choice of $\gamma_i$ and $\ve c_i$)}
		\State Update the residual $\ve r = \ve y - \ve \Psi_i \ve c_i$
		\EndFor
	\end{algorithmic}
\end{algorithm}

The LARS technique is presented in Algorithm \ref{alg:originalLAR}.
It returns a sequence $\ca_1 \subset \ca_2 \subset \ldots \subset \ca_m$ of sets containing indices of active basis functions, with $m = \min\{N,P\}$. 
Just like OMP, LARS can be stopped when a predefined sparsity $K$ is reached or when the norm of the residual $\norme{\ve r}{2}$ falls below a predefined error threshold.

A modified version of LARS, called \textit{ybrid LARS}, uses the equicorrelated approach to select the predictors, but computes the coefficients of the metamodel by least-squares \citep{Efron2004, BlatmanJCP2011}.
Once the \changed{LARS} algorithm has finished and returned the sequence of basis sets $\ca_1 \enum \ca_m$, the corresponding coefficients are recomputed by least squares, $\ve c_i = \ve c_i^\text{LSQ}$, which ensures minimal empirical error for every metamodel $(\ca_i, \ve c_i)$. 
Hybrid LARS facilitates another way to choose the best sparsity level: as for OMP, a model selection criterion (e.g.\ LOO) for each metamodel can be evaluated, and the best one is chosen. This procedure is detailed in Algorithm \ref{alg:hybrid-LAR}.
Cheap OLS-based computation of LOO \citep{Chapelle2002, BlatmanJCP2011} and the early-stop criterion \citep{UQdoc_13_104} can be applied as well.

\begin{algorithm}
	\caption{Hybrid \changed{LARS} with LOO-CV \citep{BlatmanJCP2011, UQdoc_13_104}}
	\label{alg:hybrid-LAR}
	\begin{algorithmic}[1]
		\State Initialization as in \changed{LARS} (Algorithm \ref{alg:originalLAR})
		\For{$i=1\enum m$}
		\State Run one step of \changed{LARS} and obtain $(\ca_i, \ve c_i)$
		\Graycomment{Algorithm \ref{alg:originalLAR}}
		\State Recompute the coefficient vector using least-squares on the selected basis $\ca_i$ only, obtaining $\ve c_i^\text{OLS}$ (the coefficients corresponding to $\ca \setminus \ca_i$ are set to zero)
		\Graycomment{Hybrid \changed{LARS}}
		\State Compute the LOO error $\epsilon_\text{LOO}(i)$ for $\ve c_i^\text{OLS}$
		\Graycomment{OLS-based LOO computation \citep{Chapelle2002, BlatmanJCP2011}} 
		\EndFor
		\Graycomment{early stopping possible by monitoring the LOO error \citep{UQdoc_13_104}}
		\State Return the metamodel $(\ca_{i^*}, \ve c_{i^*})$ with $i^* = \arg\min_i \epsilon_\text{LOO}(i)$
	\end{algorithmic}
\end{algorithm}

As for OMP, the computational complexity of LARS (in the case $N < P$) is $\co(mNP)$,
where $m \leq \min\{N,P\}$ is the number of iterations.  
This is due to matrix-vector multiplication and matrix inversion which have to be performed in every iteration. The latter can be computed in $\co(mN)$, when using techniques such as Schur's complement to update the information matrix inverse whenever a new regressor is added. 

LARS is available in many software packages, e.g., as MATLAB implementation in UQLab \citep{MarelliUQLab2014}.

% ====================================================================================================

\subsection{Subspace pursuit (SP)}
\label{sec:SP}
Another formulation of the $\ell^0$-minimization problem is
\begin{equation}
\min_{\ve c \in \R^P} \norme{\ve \Psi \ve c - \ve u}{2} \text{ s.t. } \norme{\ve c}{0} = K
\label{eq:l0-min-K}
\end{equation}
which is equivalent to \eqref{eq:l0-min} for a certain choice of $\lambda$. 

Subspace pursuit (SP) seeks to identify a solution to \eqref{eq:l0-min-K} by iteratively and greedily enlarging and shrinking the set of active basis functions \citep{Dai2009}. 
As with LARS and OMP, regressors are added to the set of active basis functions according to their correlation with the residual. However, the regressors are not added one by one, but batchwise. More precisely, SP maintains at all times an active basis of size $K$, where $K$ denotes the desired sparsity. In each iteration, it adds $K$ regressors at once and computes the coefficients of the active regressors by OLS.
Then, it removes the $K$ regressors with the smallest-in-magnitude coefficients.
This is continued until convergence.
Under certain assumptions, there are theoretical guarantees for the solution that SP returns \citep{Dai2009}.
To make the augmentation of the basis and the OLS regression feasible, it must hold that $2K \leq \min\{N,P\}$.

The technique is described in Algorithm \ref{alg:SP} for a fixed value of sparsity $K$. 
The residual of a vector and a regression matrix is defined as
\begin{equation}
\text{residual}(\ve y, \ve \Psi) = \ve y - \ve \Psi \ve \Psi^\dagger \ve y
\end{equation} 
where $\ve \Psi^\dagger$ denotes the pseudoinverse of $\ve \Psi$ and $\ve \Psi^\dagger \ve y = \ve c$ is the least-squares solution to $\ve \Psi \ve c \approx \ve y$ (the case of an overdetermined system).
The algorithm returns a set $\ca$ containing $K$ multi-indices.

\begin{algorithm}
	\caption{Subspace pursuit (SP) \citep{Dai2009}.}
	\label{alg:SP}
	\begin{algorithmic}[1]
		\State Given desired sparsity $K \leq \min\{\frac{N}{2}, \frac{P}{2}\}$
		\State Given the experimental design and the candidate basis $\ca_\text{cand}$, compute the associated regression matrix $\ve\Psi$ and the right-hand-side $\ve y$ 
		\State $\ca^0 = \{K \text{ indices corresponding to the largest magnitude entries in } \ve\Psi^T\ve y\}$
		\Statex\Graycomment{Scalar product of columns of $\ve\Psi$ with $\ve y$}
		\State $\ve y_\text{res}^0 = \text{residual}(\ve y, \changed{\ve \Psi_{\ca^0}})$
		\Graycomment{Residual of least-squares solution based on \changed{active basis}}
		\For{$l = 1,2, \ldots$}
		\State $\cs^l = \ca^{l-1} \cup \{K \text{ indices corresponding to the largest magnitude entries in } \ve \Psi^T \ve y_\text{res}^{l-1}\}$
		\Statex \Graycomment{Augment by indices of full basis that correlate best with the residual}
		\State $\ve c = \ve \Psi_{\cs^l}^\dagger \ve y$
		\Graycomment{Least-squares solution based on set $\cs^l$ of size $2K$}
		\State $\ca^l = \{ K \text{ indices corresponding to the largest magnitude entries in } \ve c\}$
		\label{line:newbasis}
		\State $\ve y_\text{res}^l = \text{residual}(\ve y, \ve \Psi_{\ca^l})$
		\Graycomment{Residual of least-squares solution based on $\ca^l$}
		\If{$\norme{\ve y_\text{res}^l}{2} \geq \norme{\ve y_\text{res}^{l-1}}{2}$}
		\label{line:SP-stop}
		\Graycomment{if new $K$-sparse approx.\ is worse than the previous one}
		\State STOP iteration and \Return $\ca^{l-1}$.	
		\EndIf		
		\EndFor
	\end{algorithmic}
	\textit{Remark: In line \ref{line:SP-stop}, the original publication \citep{Dai2009} uses ``$>$'' instead of ``$\geq$'', but we also want to stop when the set has converged.}
\end{algorithm}

For arbitrary sparse vectors, the computational complexity is $\co(N(P + K^2)K)$ \citep{Dai2009}. For very sparse vectors with $K^2 \in \co(P)$, the complexity thus becomes $\co(NPK)$, comparable to the runtime of OMP. The number of iterations that the SP algorithm performs can be shown to be $\co(K)$ in general and even $\co(\log K)$ in certain cases \citep{Dai2009}.

When the optimal sparsity level $K$ is unknown, it can be determined e.g.\ by cross-validation: \citet{Diaz2018} suggest running Algorithm \ref{alg:SP} for a range of $N_K = 10$ different values for $K$ and choosing the one with the smallest 4-fold cross-validation error.
\changed{In this paper, we propose to use leave-one-out cross-validation instead of 4-fold cross-validation, resulting in the SP variant \SPloo{}.}

A related algorithm is CoSAMP \citep{Needell2009}, which differs from SP mainly in the number of regressors added in each iteration.

Subspace pursuit is available as MATLAB implementation in the software package \texttt{DOPT\_PCE} \citep{Diaz2018,SPcode}.

% ====================================================================================================

\subsection{SPGL1}
\label{sec:convex}
$\ell^1$-minimization is a convex problem, since both the objective function and the constraint are convex functions. Therefore, convex optimization methods can be used to find a solution. In this section, we describe the algorithm SPGL1 \citep{Vandenberg2008}.

For a given value of $\tau$, formulation \eqref{eq:LASSO} (LASSO) can be solved by spectral projected gradient (SPG) descent \citep{Birgin2000, Vandenberg2008}.%
\footnote{SPG is a gradient-based optimization algorithm with several enhancements (Barzilai--Borwein spectral step length and the Grippo--Lampariello--Lucidi scheme of nonmonotone line search) and projection onto the feasible set $\Omega_\tau = \{\ve c \in \R^P: \norme{\ve c}{1} \leq \tau\}$}
However, for real-world problems, we often do not know a priori an appropriate value for $\tau$. 
On the other hand, a sensible range of values for $\sigma$ in formulation \eqref{eq:BPDN} (BPDN) can typically be estimated based on the noise level in the data and the expected model fit. 
In the case of PCE metamodelling, $\sigma$ can be related to an estimate of the relative MSE through RelMSE $= \frac{\sigma^2}{N \Varhat{\ve y}}$, whose values are for engineering models typically between $10^{-10}$ and $10^0=1$.

The main idea of the solver SPGL1 is to solve BPDN through a detour over LASSO. 
Let $\ve c_\tau$ be the solution to LASSO for a given $\tau$. Define a function $\phi:\R_+ \to \R_+$ by
\begin{equation}
\phi(\tau) := \norme{\ve\Psi \ve c_\tau - \ve y}{2}.
\end{equation}
Then the solution to BPDN with $\sigma := \phi(\tau)$ is $\ve c_\tau$. 
In other words, $\phi$ is the functional relationship between $\sigma$ and $\tau$ that makes the two formulations BPDN and LASSO equivalent for given $\ve \Psi$ and $\ve y$.
$\phi$ is the Pareto front of LASSO and BPDN and shows the trade-off between the minimal achievable $\ell^1$-norm of the coefficients and the minimal $\ell^2$-norm of the corresponding residual. 
The Pareto front is convex, nonincreasing and differentiable with an analytically computable derivative \citep{Vandenberg2008}.

To find a solution to BPDN with a given $\sigma$, LASSO is solved with SPG several times for a sequence of $\tau$ until one is found with $\phi(\tau) = \sigma$. The sequence of $\tau$ is created by performing Newton's root finding algorithm on the function $f(\tau) = \sigma - \phi(\tau)$.

Each SPG iteration has a computational complexity of $\co(NP + P \log P)$ (from matrix-vector multiplication and $\ell^1$-projection). Multiplying this with the number of SPG steps and the number of Newton steps yields the computational complexity of SPGL1.

This algorithm is available as MATLAB package \texttt{SPGL1} \citep{Vandenberg2008, SPGL1}. 

In our numerical benchmarks computing sparse PCE for compressible models, SPGL1 was among the slowest solvers and often returned rather dense solutions.

% ====================================================================================================

\subsection{Sparse Bayesian learning}
\label{sec:SBL}
Methods from the class of Bayesian compressive sensing (BCS), also known as sparse Bayesian learning (SBL), 
embed the regression problem in a probabilistic framework \changed{\citep{Tipping2001, Ji2008, Babacan2010, Sargsyan2014, Tsilifis2020}}. 
The goal is to compute, for a given model response vector $\ve y$ and a regression matrix $\ve \Psi$, the coefficient vector $\ve c^\text{MAP}$ which maximizes the posterior distribution $p(\ve c | \ve y)$. Another quantity of interest could be the most probable value $y^*$ at a new point $\ve x^*$ maximizing $p(y^* | \ve y)$.

In BCS, it is assumed that the ``measurements'' $\ve y$ are generated by adding zero-mean, finite-variance noise to the evaluations of the true model. This noise is often assumed to be Gaussian white noise \changed{with standard deviation $\sigma$}, which, for a given input $\ve x$, results in a Gaussian distribution for its output $\ve y$ with mean $\ve \Psi \ve c$ and covariance matrix $\sigma^2 \mathbb{1}$
\changed{, i.e., $\ve y | \ve c, \ve x, \sigma \sim \cn(\ve\Psi\ve c, \sigma^2\mathbb{1})$}.
Note that in the case of PCE, this is generally not a valid assumption: when an important term is missing from the PCE model, the discrepancy between measurements and PCE model evaluations can be highly correlated, heteroscedastic, and non-Gaussian, and have nonzero mean. 
However, even though the assumptions might not be fulfilled, this framework can still be useful for finding sparse solutions.

The class of BCS algorithms comprises several methods that differ in the assumptions on the distributions of the various hyperparameters and in the (usually iterative, approximate) techniques for computing the posterior quantities.

\begin{figure}[htbp]
	\subfloat[][General setup of BCS]{\includegraphics[width=.4\textwidth]{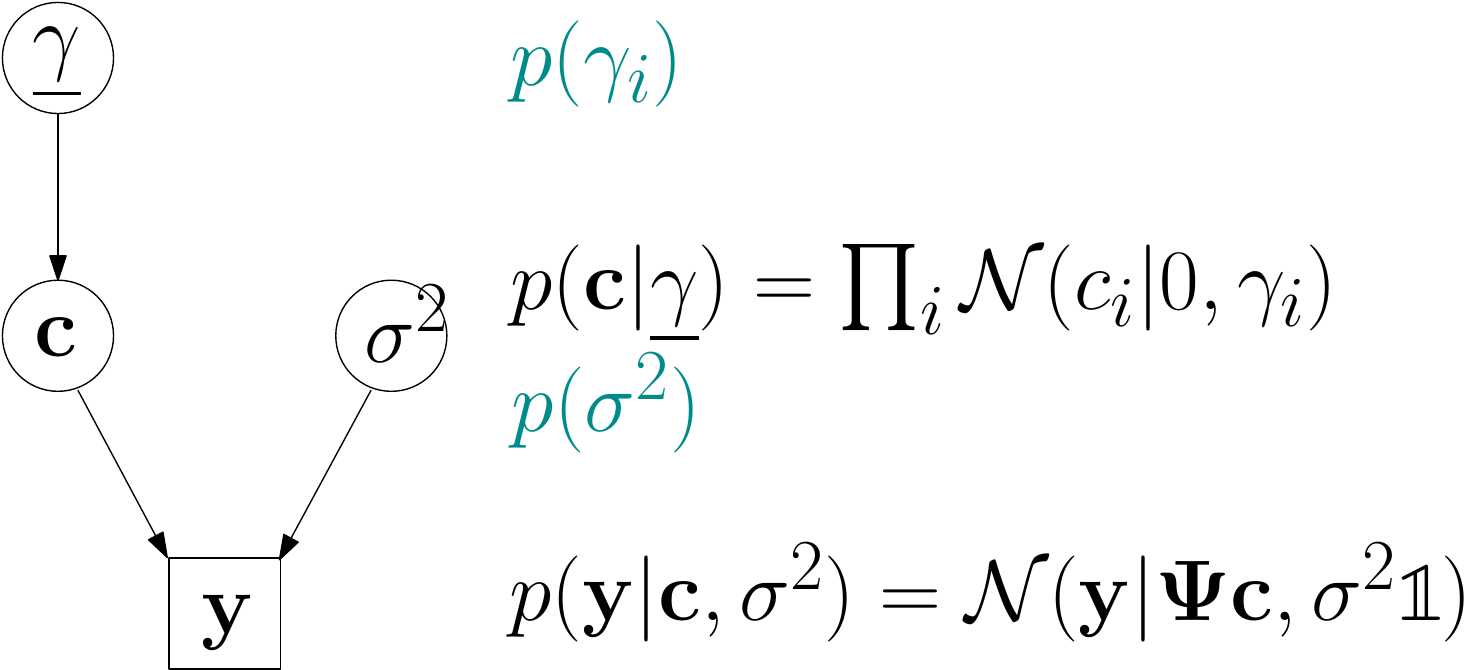}\label{fig:generalsetup}}
	\hfill
	\subfloat[][Setup by \citep{Babacan2010}]{\includegraphics[width=.4\textwidth]{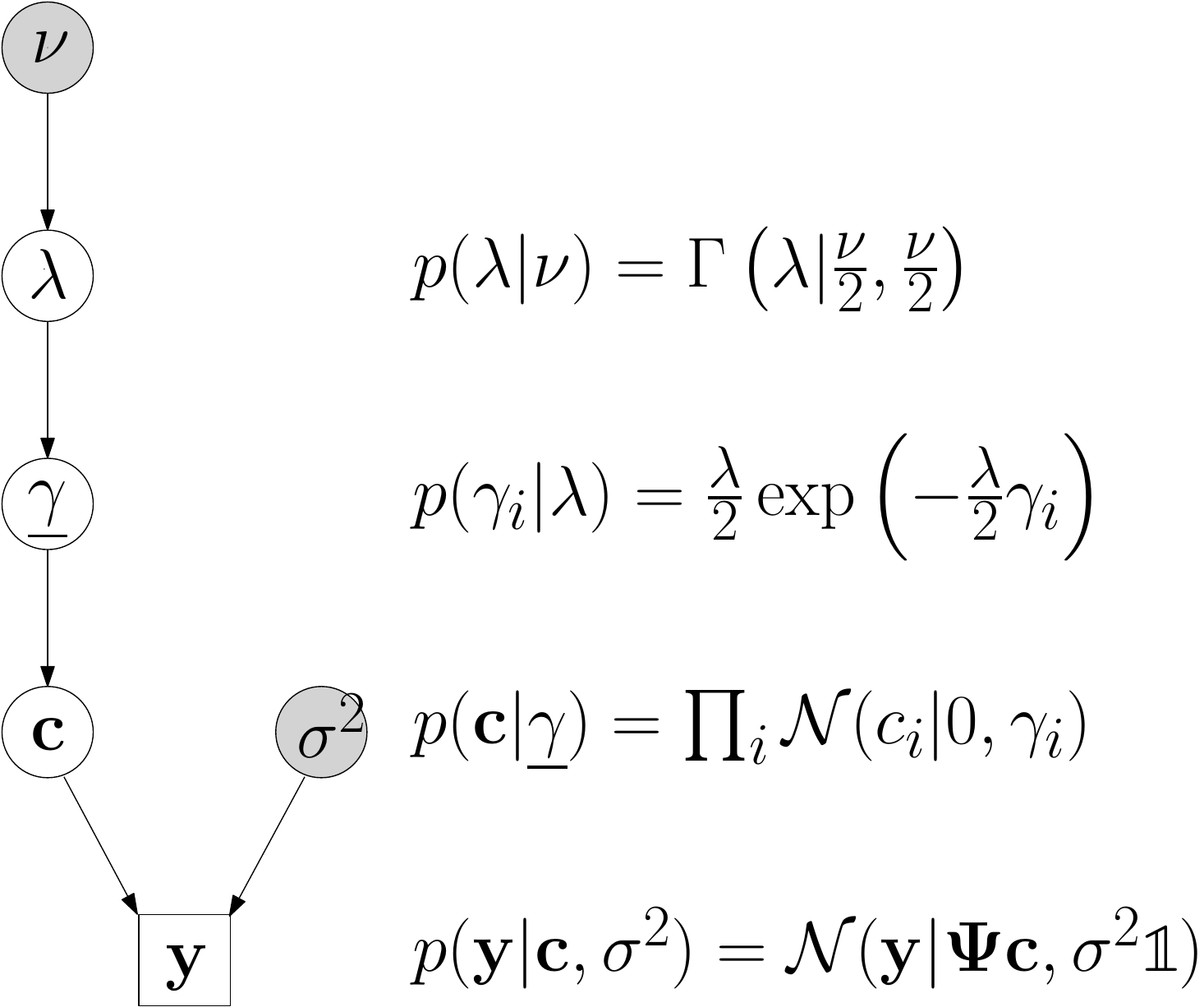}\label{fig:babacan}}
	\caption{Illustration of the general setup of BCS (a) and the hierarchical generalization of \citep{Babacan2010} (b). (a): The likelihood and the prior on the coefficients are usually Gaussian, but the choice of $p(\sigma^2)$ and $p(\ve\gamma)$ differs between publications, as well as the resulting solution algorithm. (b): \citet{Babacan2010} makes a specific choice for $p(\ve\gamma)$ and includes an additional layer of hyperparameters. Shaded variables are held fixed.}
\end{figure}

In Figure \ref{fig:generalsetup} we present the general setup of sparse Bayesian learning. 
The measurements $\ve y$ are assumed to follow a Gaussian distribution as described above. 
The noise variance $\sigma^2$ is assumed to be a random variable whose distribution has to be specified (e.g.\ fixed, uniform, or inverse-Gamma).
The coefficients $c_i$ are assumed to be random variables as well, drawn from a normal distribution with mean zero and variance $\gamma_i$, i.e., each weight has its own variance. $\ve\gamma$ is a so-called \textit{hyperparameter}, parametrizing the distribution of a parameter.

So far, the described setup with fixed $\sigma^2$ and $\gamma_i$ would yield (weighted) ridge regression.
The sparsity comes into play through an assumption on the distribution of the hyperparameter $\ve\gamma$. For specific choices of $p(\gamma_i)$, it can be shown that the resulting \textit{effective prior} on the coefficients $p(\ve c) = \int p(\ve c | \ve \gamma)p(\ve\gamma) \di{\ve\gamma}$ is a sparsity-encouraging distribution, i.e., one that has a sharp peak at zero, encouraging zero values, while at the same time a heavy tail, allowing for large coefficient values as well.
Examples are the Laplace distribution and the generalized Student-t distribution \citep{Wipf2004a, Babacan2010, Figueiredo2001}.

Such sparsity-encouraging distributions are often intractable to use, because they do not allow for analytical computation of the desired values (such as the most likely coefficients given the data, or the prediction of the measurement value at a new point). 
However, feasible algorithms can be developed based on a suitable approximation step. 
Various frameworks for sparse Bayesian learning have been proposed whose setup follows the general structure of Figure \ref{fig:generalsetup}, but which differ in the choice of priors for the hyperparameters 
and employ different solution algorithms for the MAP estimate of these hyperparameters
\changed{\citep{Tipping2001,Faul2002,Figueiredo2003,Tipping2003,Wipf2004b, Ji2008, Seeger2008, Babacan2010,Sargsyan2014,Tsilifis2020, Bhattacharyya2020}}.
The MAP estimate of the hyperparameters is inserted into the distribution for $\ve c$. Because $p(\ve c | \ve \gamma)$ and $p(\ve y | \ve c, \sigma^2 )$ are normal distributions, all subsequent computations can be carried out analytically \citep{Tipping2001,Wipf2004a, Wipf2004b}.
The sparsity of $\ve c$ is enforced because by the choice of $p(\ve\gamma)$ and the other distributions, many of the components $\gamma_i^\text{MAP}$ of $\ve\gamma^\text{MAP}$ will actually be zero, forcing the corresponding $c_i$ to be zero as well.

BCS in the implementation of \citep{Babacan2010} was suggested for sparse PCE by \citep{Sargsyan2014}.
This approach employs an additional layer of hyperparameters as displayed in Figure~\ref{fig:babacan}.
The prior on the coefficient variances is an exponential distribution $p(\gamma_i | \lambda) = \text{Exp}\left( \gamma_i \big| \frac{\lambda}{2} \right)$ with shared hyperparameter $\lambda$. 
The hyperparameter $\lambda$ follows a Gamma distribution $p(\lambda | \nu) = \Gamma(\lambda | \frac{\nu}{2}, \frac{\nu}{2})$ with hyperparameter $\nu$. $\nu \to 0$ implies $p(\lambda) \propto \frac{1}{|\lambda|}$ (improper prior) and $\nu \to \infty$ implies the certain value $\lambda = 1$. In practice, \citep{Babacan2010} find that $\nu = 0$ gives the best results.
The prior on $\beta = \sigma^{-2}$ is a Gamma distribution $p(\beta) = \Gamma(\beta|a,b)$ with hyperparameters $a,b$. In practice, the algorithm does not estimate $\beta$ well, which is, however, crucial; therefore, it is set to a fixed value (e.g.\ $\beta^{-1} = 0.01\norme{\ve y}{2}^2$ in \citep{Babacan2010}; in our benchmark, we use cross-validation to determine the best value for this parameter, similarly to the strategy for SPGL1).
The objective function is the logarithm of the joint distribution $\cl(\ve\gamma, \lambda, \beta) = \log p(\ve y, \ve\gamma, \lambda,\beta)$, which is an analytical expression. To maximize it, \citet{Babacan2010} adapt the fast approximate algorithm of \citet{Tipping2003, Faul2002} to their generalized hierarchical setting. Here, the derivatives of the objective function with respect to the hyperparameters $\lambda, \beta$, and $\gamma_i, i = 1\enum P$ are computed. This results in an iterative scheme where these parameters are optimized one at a time while the other ones are held fixed.  
The algorithm is explained in detail in \citet[Algorithm 1]{Babacan2010} and has been implemented in MATLAB under the name FastLaplace \citep{FastLaplace}.

% ====================================================================================================

\subsection{Greedy stepwise regression solvers}
\label{app:greedystepwise}
Many of the sparse regression solvers that have been proposed for computing sparse PCE belong to the class of greedy stepwise regression. 
Here, starting from an empty model, the regressors are added one by one according to a selection criterion (forward selection). Some methods also include a backward elimination step. Then the  coefficients of the selected regressors are computed. The procedure is iterated until a stopping criterion is reached. Alternatively, several models are built and one is selected in the end using a model selection criterion. We summarize some greedy stepwise regression techniques proposed for sparse PCE, together with their choices for selection criterion, coefficient computation, and stopping criterion, in Table~\ref{tab:stepwiseregression}, including the well-known methods OMP and LARS.

New greedy methods in the fashion of Table~\ref{tab:stepwiseregression} can easily be derived by pairing other methods for the regressor selection, the coefficient computation method, and the model selection criterion.
Note that except for LARS and OMP, these greedy methods are heuristic (no theoretical guarantee of convergence) and often depend on a number of tuning parameters.

\begin{sidewaystable}[tbhp]
	\tiny
	\caption{A selection of greedy stepwise regression algorithms proposed for sparse PCE. F = forward selection, FB = forward selection and backward elimination. }
	\label{tab:stepwiseregression}
	\begingroup
	\renewcommand{\arraystretch}{1.5}
	\begin{tabular}{p{.1\textwidth}  p{.04\textwidth} p{.15\textwidth} p{.14\textwidth} p{.23\textwidth}  p{.25\textwidth} }
		\hline
		\textbf{Ref.} & \textbf{F/B} & \textbf{Regressor selection} & \textbf{Computation of}\newline \textbf{the coefficients} & \textbf{Model selection /} \newline \textbf{stopping criterion} & \textbf{Comment}
		\\ \hline %\hline
		OMP \newline \citep{Tropp2007} & F & correlation with \newline residual & OLS & \textit{several choices:} \newline given number of regressors \citep{Tropp2007}; (modified) LOO \citep{UQdoc_13_104}; norm of residual \citep{Baptista2019}; threshold on moving average of LOO \citep{Baptista2019}
		& theoretical guarantees exist \citep{Tropp2007}
		\\ %\hline
		LARS \newline \citep{Efron2004,BlatmanJCP2011} & F & correlation with \newline residual & least angle strategy & \textit{several choices:} \newline coefficient of determination \citep{Efron2004}; (modified) LOO \citep{BlatmanJCP2011}
		& theoretical guarantees exist \citep{Efron2004, Bruckstein2009}
		\\ %\hline \hline
		\citep{BlatmanCras2008} & FB & coefficient of \newline determination  & OLS & coefficient of determination
		& degree-adaptive
		\\ %\hline
		\citep{BlatmanPEM2010} & FB & coefficient of \newline determination & OLS & LOO
		& degree-adaptive; with ED enrichment
		\\ %\hline
		\citep{HuYoun2011} & F & testing bivariate \newline interaction & stepwise moving \newline least-squares & coefficient of determination 
		& Only for interaction order $\leq 2$. ED-adaptive. Degree of interaction terms determined in inner loop.
		\\ %\hline
		\citep{AbrahamJCP2017} & FB & F: one-predictor \newline regression criterion \newline B: confidence intervals of coefficients & OLS & given number of regressors or \newline iterations 
		& -
		\\ %\hline
		\citep{Shao2017} & F & partial correlation \newline  coefficient & variant of BCS \newline (fixed priors) & KIC 
		& prior depends on interaction order and degree of the regressor; \newline degree-adaptive 
		\\ %\hline
		\citep{Cheng2018a} & FB & variance contribution of coefficient & support vector \newline regression (SVR) & coefficient of determination 
		& -
		\\ %\hline
		AFBS \newline \citep{Zhao2019a} & FB & correlation with \newline residual & OLS & threshold for residual norm 
		& Main difference with OMP: backward selection. Employs several elimination checks.
		\\ %\hline
		{\tiny HSPLSR-PCE }\newline \citep{Zhao2019b} & FB & F: partial least-squares (PLS) (regressors partitioned into \newline blocks) \newline B: soft thresholding & linear PLS on the regression matrix & modified (pseudo) LOO 
		& Note that the input is assumed to be Gaussian, but the (nonlinear) \newline regressors do not follow a Gaussian distribution.
		\\ %\hline
		\citep{Zhou2019} & F & (partial) distance \newline correlation & PLS & coefficient of determination 
		& adaptive in degree and interaction \newline order 
		\\ %\hline
		\citep{Zhou2019c} & FB & correlation with \newline residual & BCS & KIC 
		& nesting iterative BCS procedure with forward-backward-selection scheme
		\\ \hline
	\end{tabular}	
	\endgroup	
\end{sidewaystable}

\clearpage

\section{Additional results}
\label{app:additional_results}%
\changed{
	In this appendix, we display additional results that complement the results shown in Sections~\ref{sec:results_only_solvers}--\ref{sec:results_nearopt}.
	For a detailed description of the setup, we refer the reader to Section~\ref{sec:benchmark}.
}

\subsection{\changed{Comparison of sparse solvers}}
\label{app:more_models}
\begin{figure}[htbp]
	\centering
	\subfloat[][Undamped oscillator ($d = 6$)]{\includegraphics[width=.49\textwidth, height=0.25\textheight, keepaspectratio]{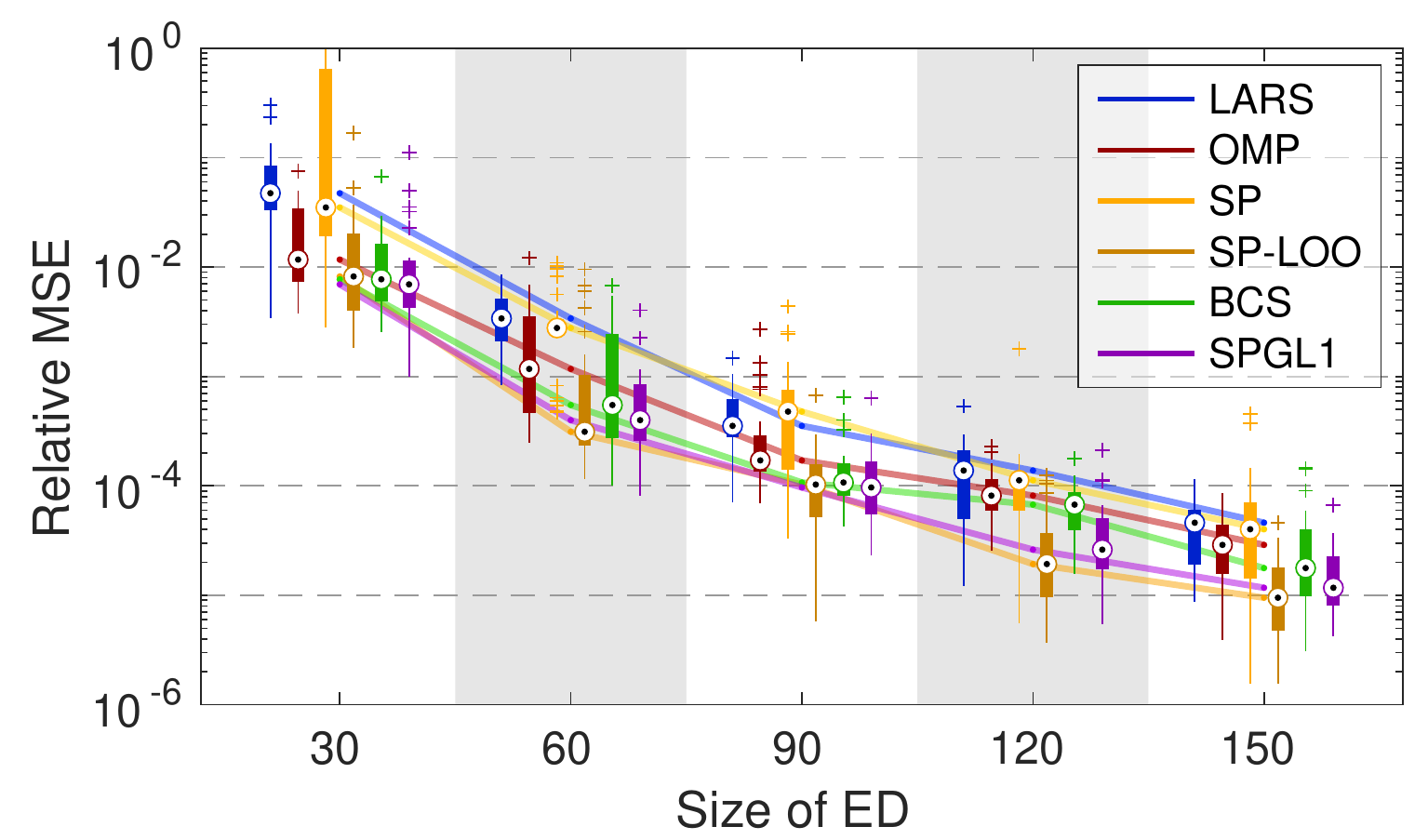}}
	\hfill
	\subfloat[][Damped oscillator ($d = 8$)]{\includegraphics[width=.49\textwidth, height=0.25\textheight, keepaspectratio]{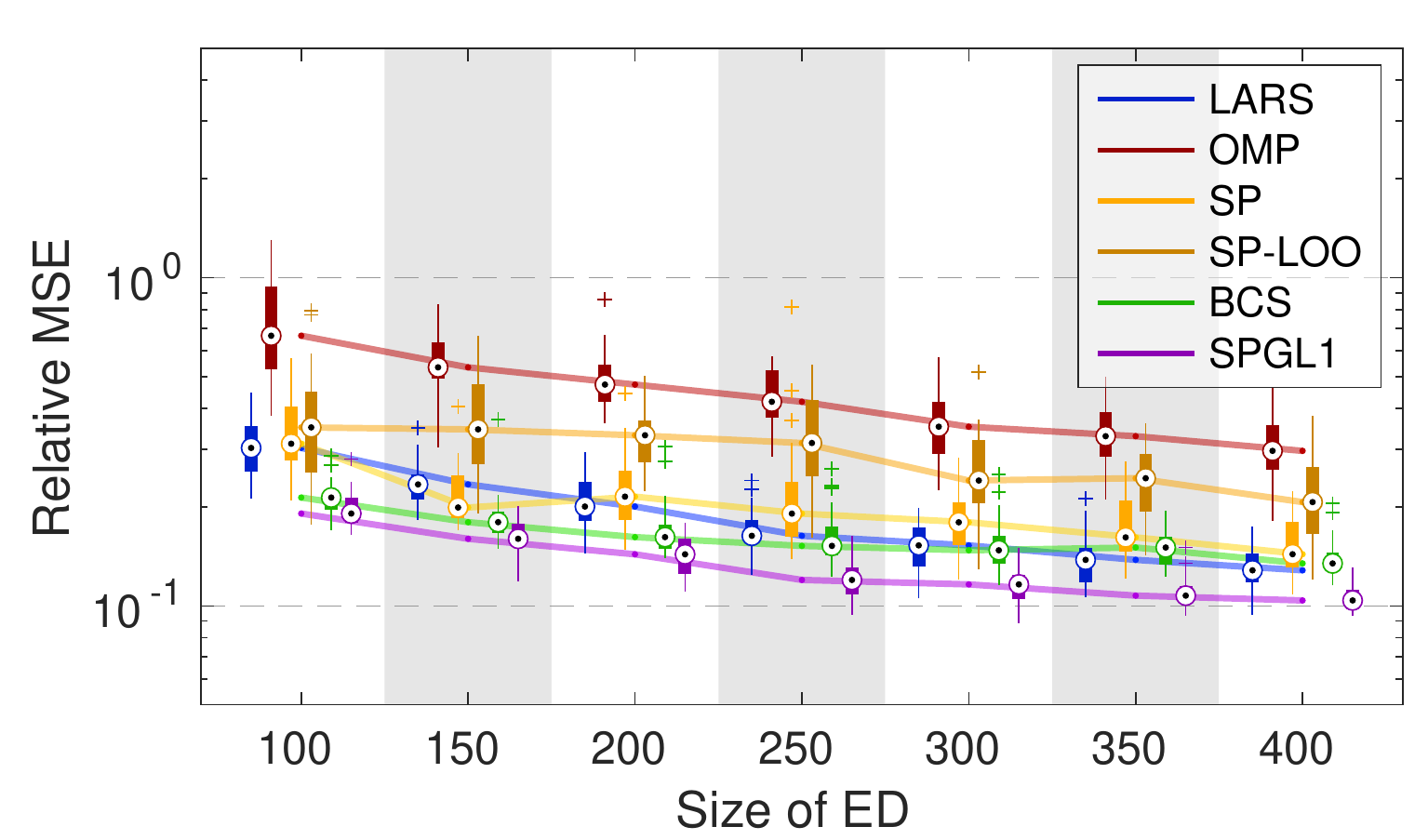}}
	\\
	\subfloat[][Wingweight function ($d = 10$)]{\includegraphics[width=.49\textwidth, height=0.25\textheight, keepaspectratio]{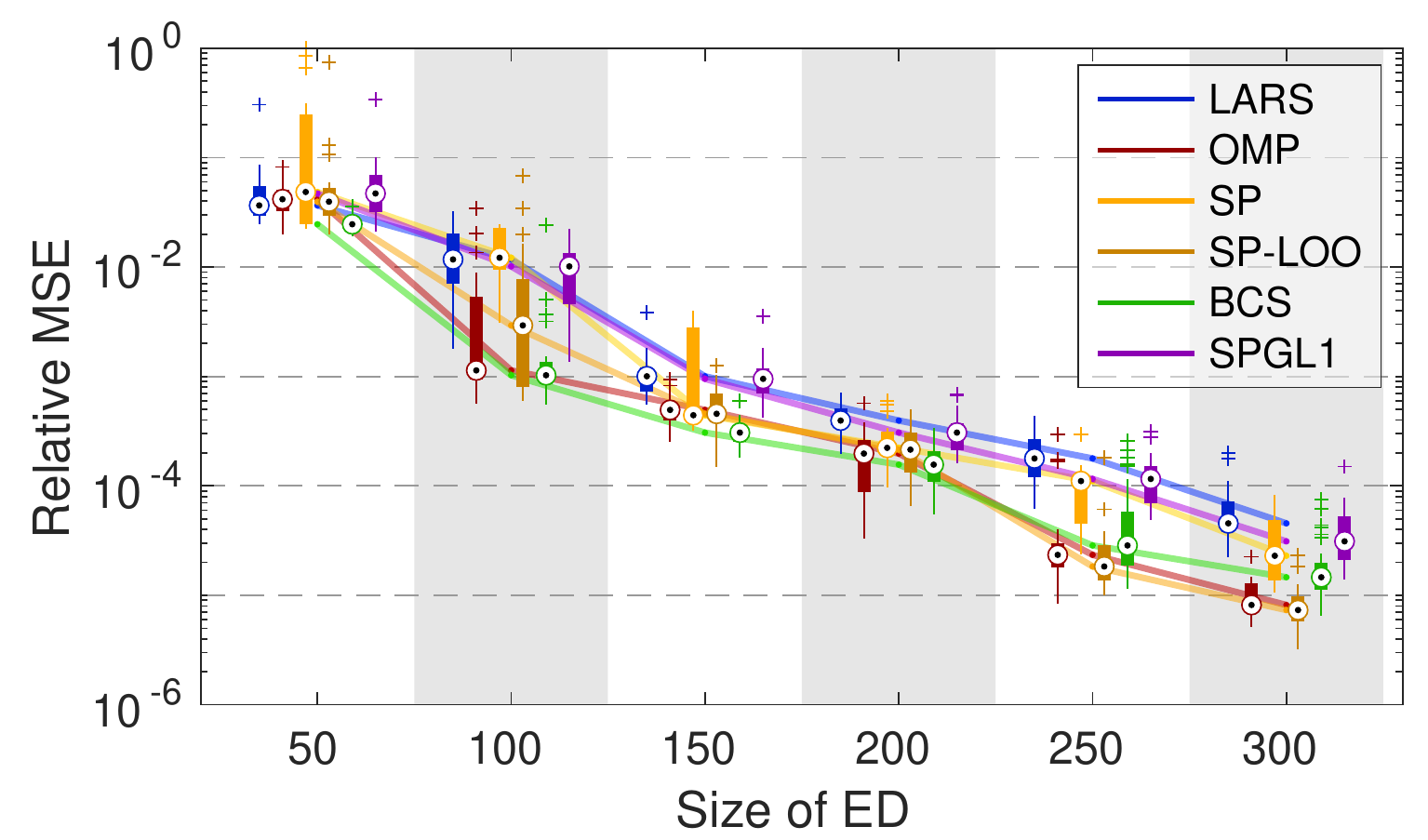}}
	\hfill
	\subfloat[][Truss model ($d = 10$)]{\includegraphics[width=.49\textwidth, height=0.25\textheight, keepaspectratio]{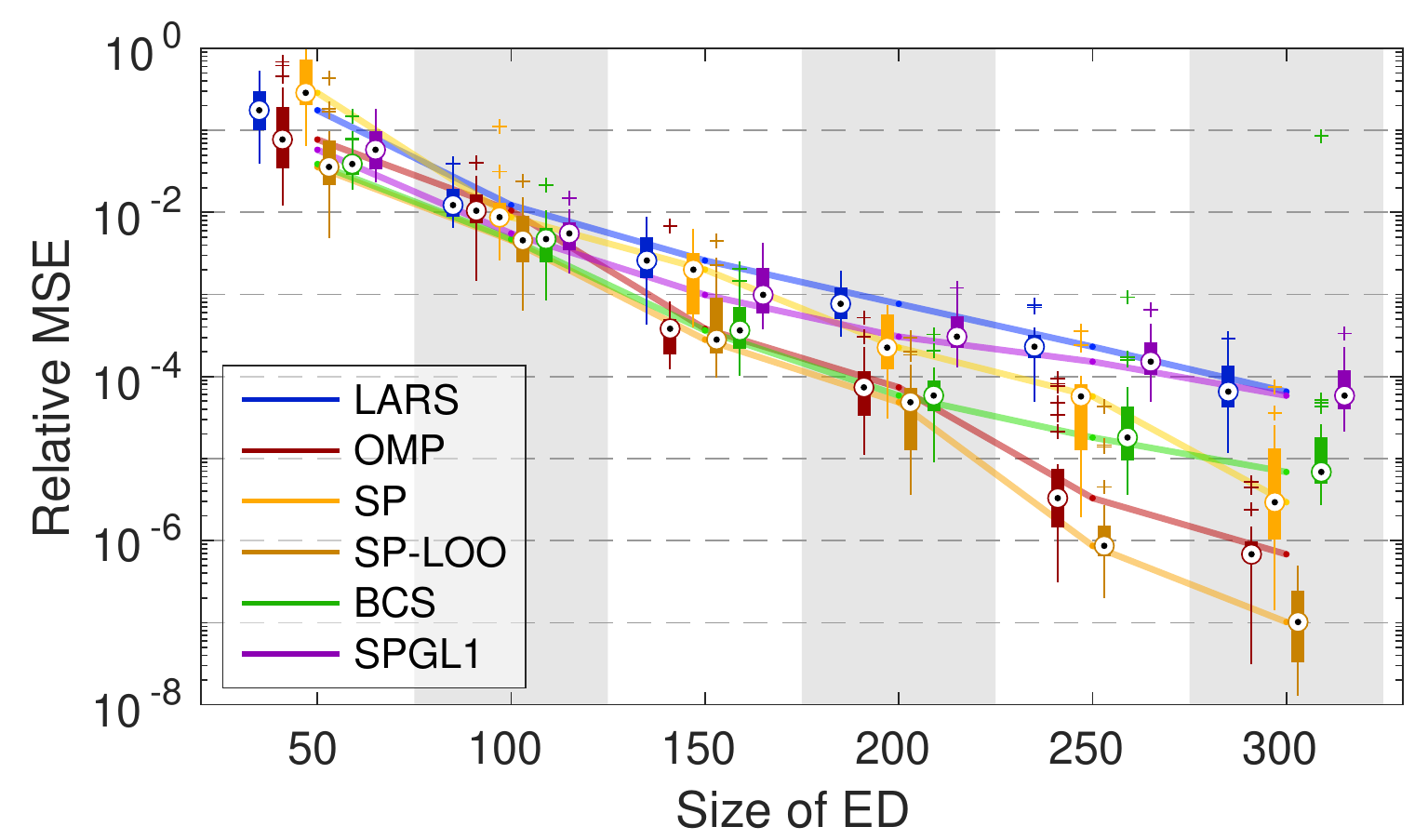}}
	\\
	\subfloat[][Morris function ($d = 20$)]{\includegraphics[width=.49\textwidth, height=0.25\textheight, keepaspectratio]{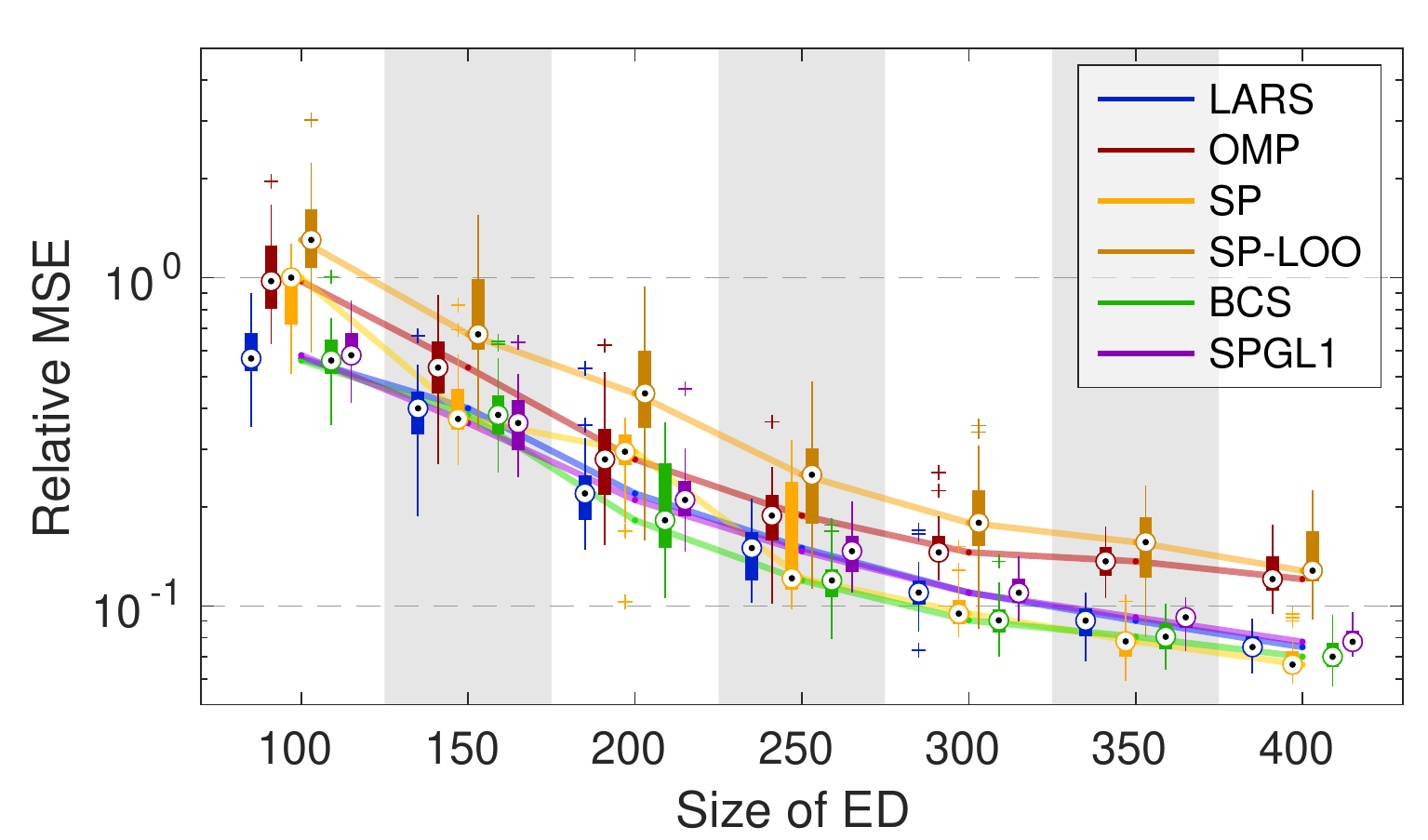}}
	\hfill
	\subfloat[][Structural frame model ($d = 21$)]{\includegraphics[width=.49\textwidth, height=0.25\textheight, keepaspectratio]{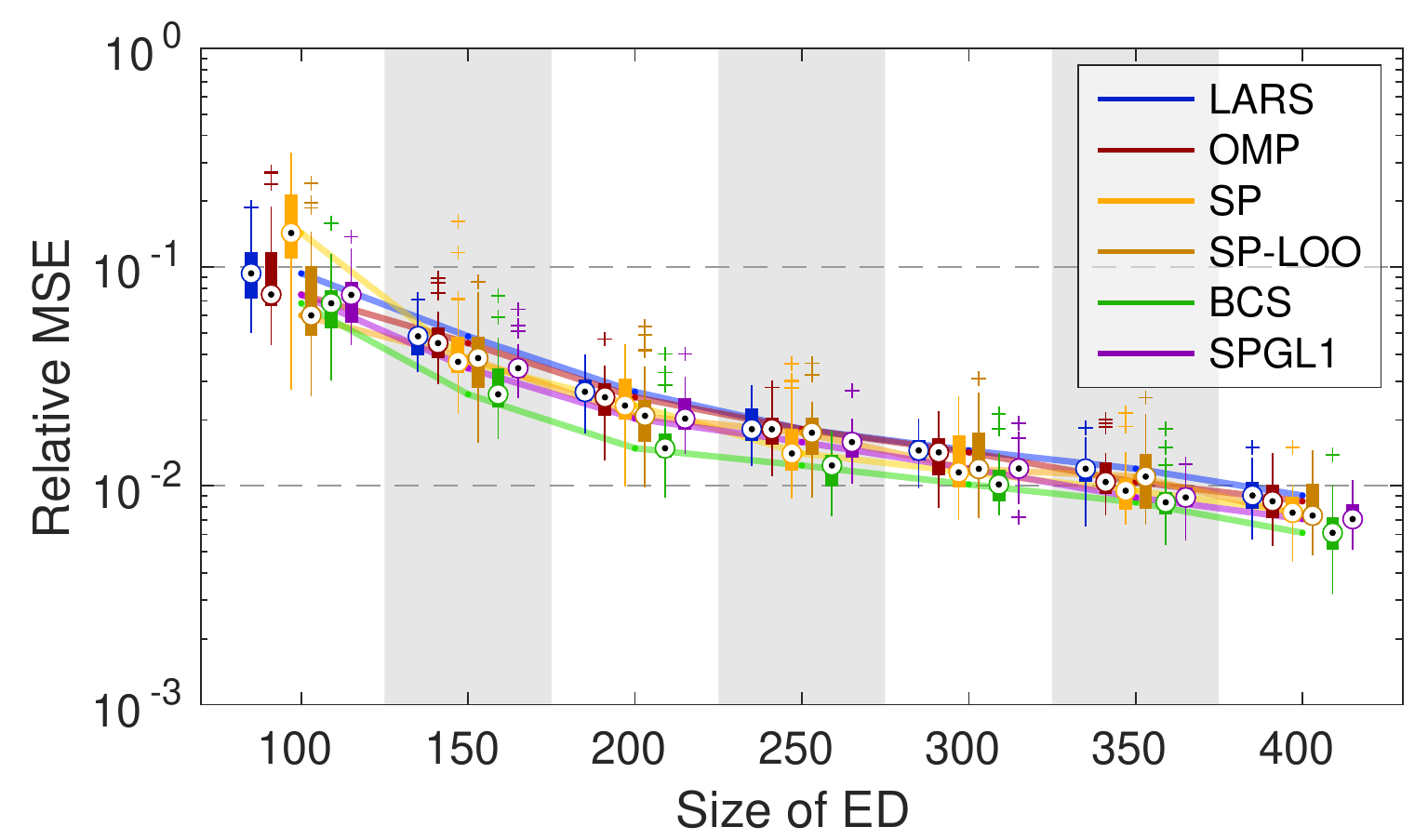}}
	\\
	\subfloat[][1-dim diffusion model ($d = 62$)]{\includegraphics[width=.49\textwidth, height=0.25\textheight, keepaspectratio]{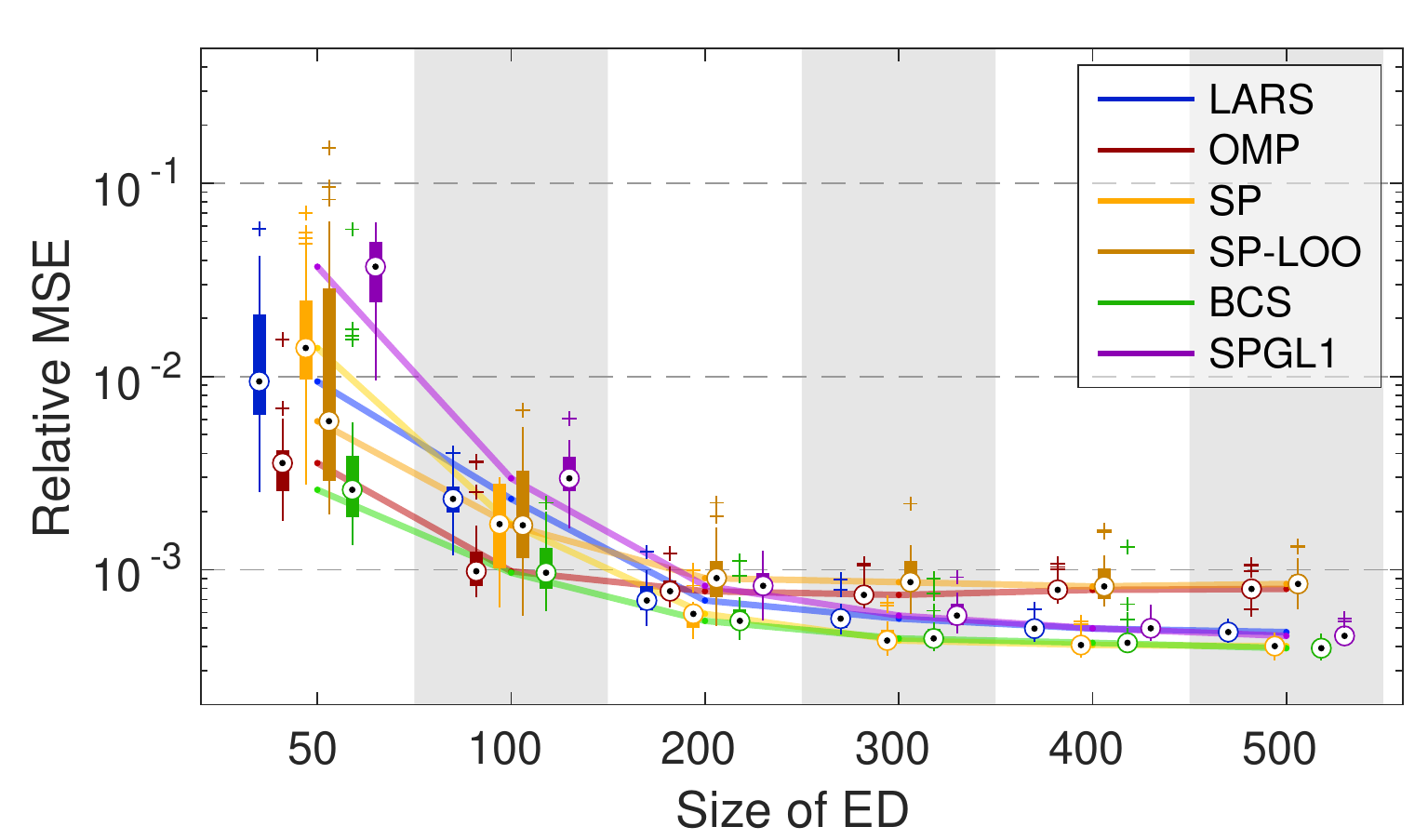}}
	\caption{\changed{Results for seven additional models (see Table~\ref{table:models} in Section~\ref{sec:results_only_solvers} for more details), complementing the results in Figure~\ref{fig:results_onlysolvers}. Boxplots of relative MSE against experimental design for six sparse solvers and LHS design. Thirty replications. Note that the damped oscillator and the Morris function are very challenging for PCE: no solver achieves a relative MSE significantly smaller than 0.1, even when large EDs are used.}}
	\label{fig:results_more_models}	
\end{figure}

\changed{In Figure~\ref{fig:results_more_models} we display the boxplots of relative MSE for the seven additional models presented in Table~\ref{table:models}. }

\subsection{\changed{Comparison of sampling schemes together with solvers}}
\begin{figure}[htbp]
	\centering
	\subfloat[][LARS]{\includegraphics[width=.49\textwidth, height=0.22\textheight, keepaspectratio]{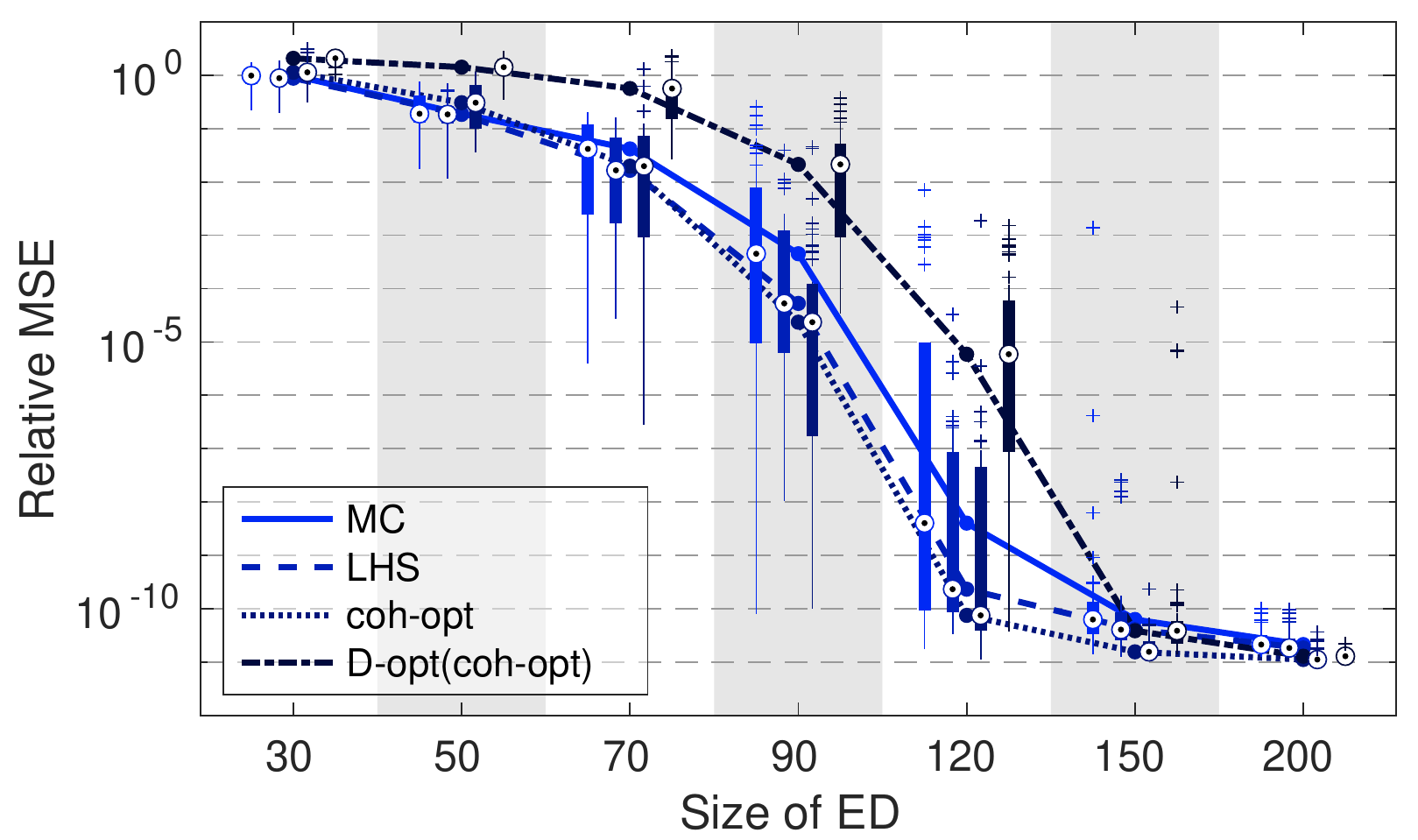}
		\label{fig:additional_results_ishigami_LARS}	}
	\hfill
	\subfloat[][OMP]{\includegraphics[width=.49\textwidth, height=0.22\textheight, keepaspectratio]{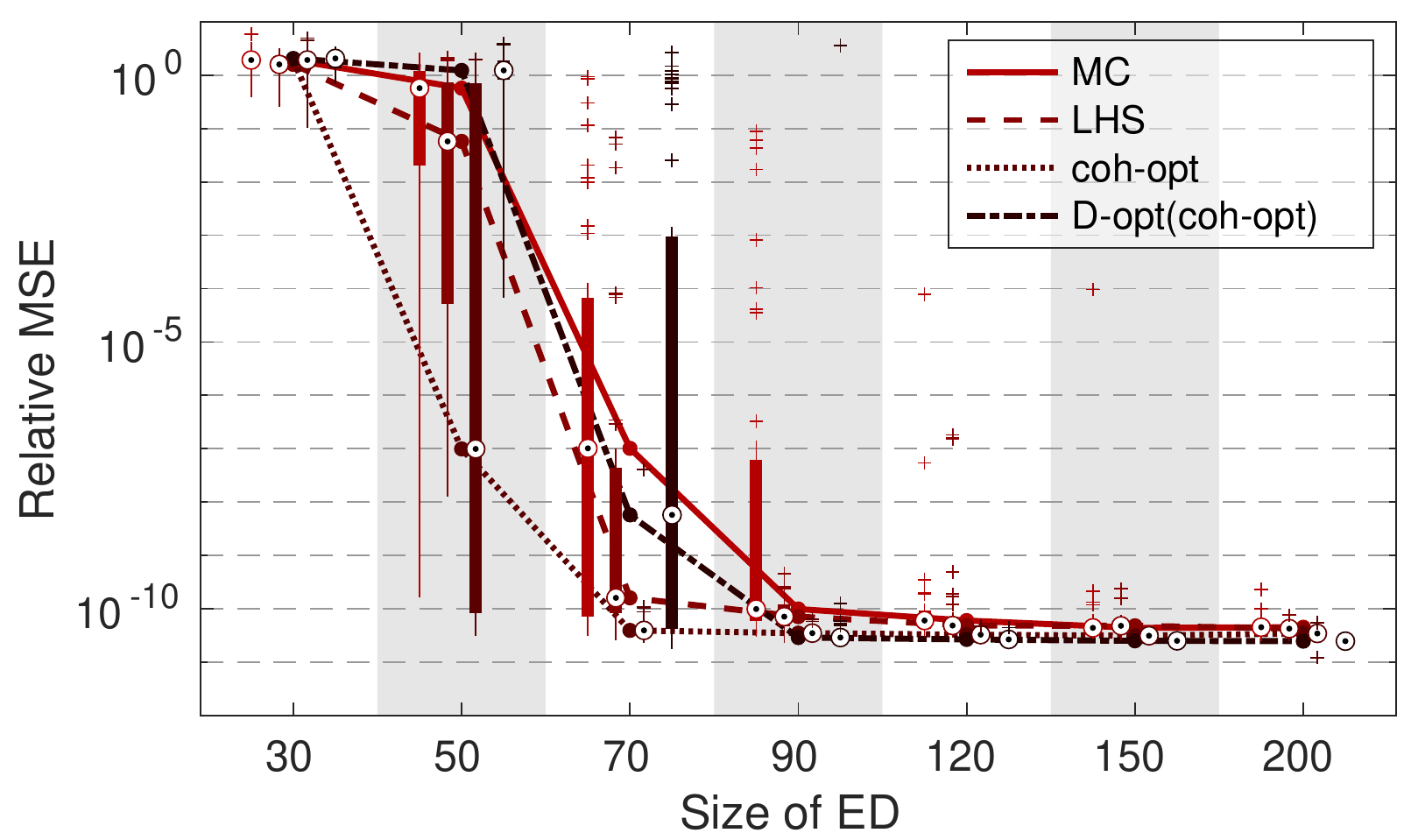}}
	\\
	\subfloat[][SP]{\includegraphics[width=.49\textwidth, height=0.22\textheight, keepaspectratio]{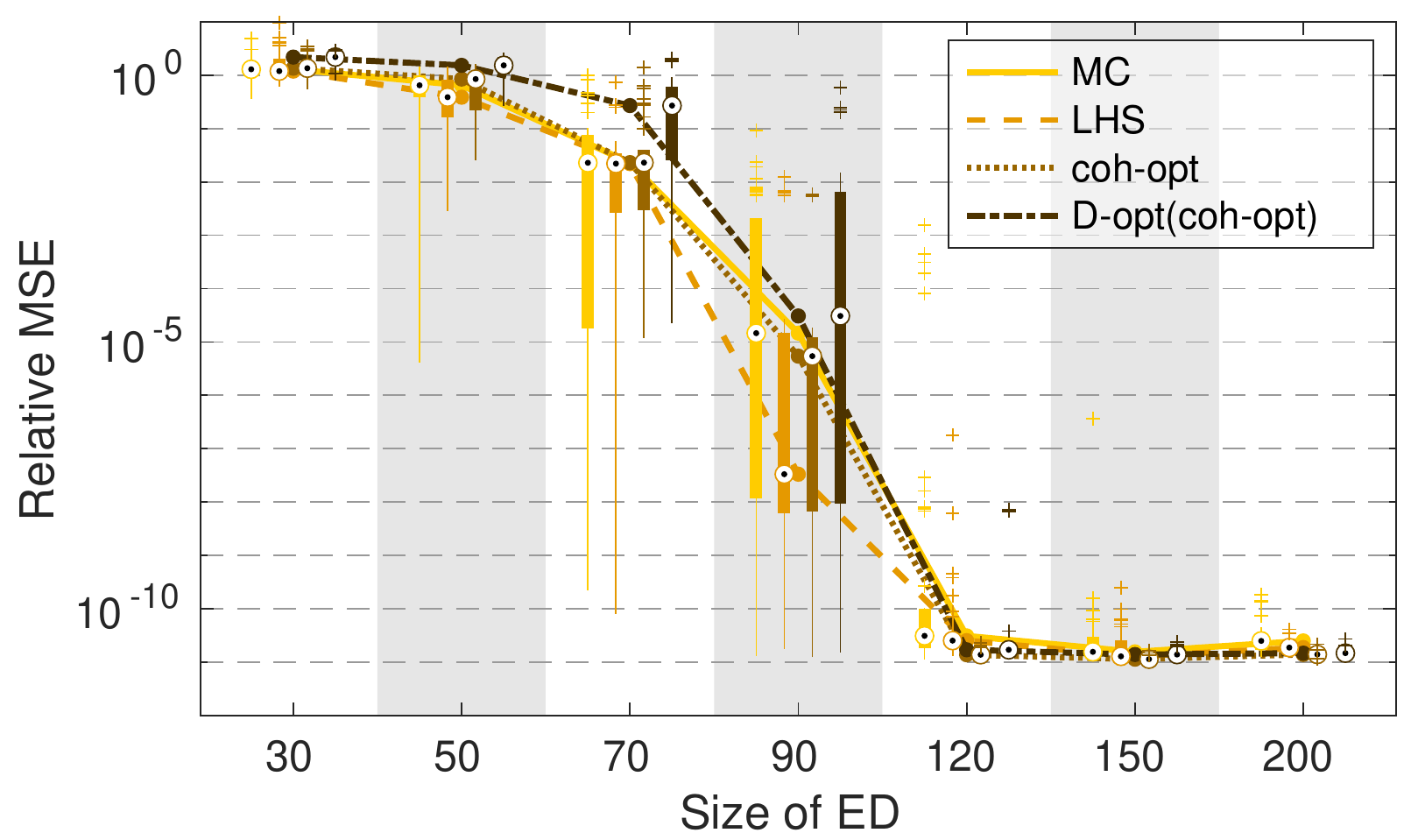}}
	\hfill
	\subfloat[][\SPloo{}]{\includegraphics[width=.49\textwidth, height=0.22\textheight, keepaspectratio]{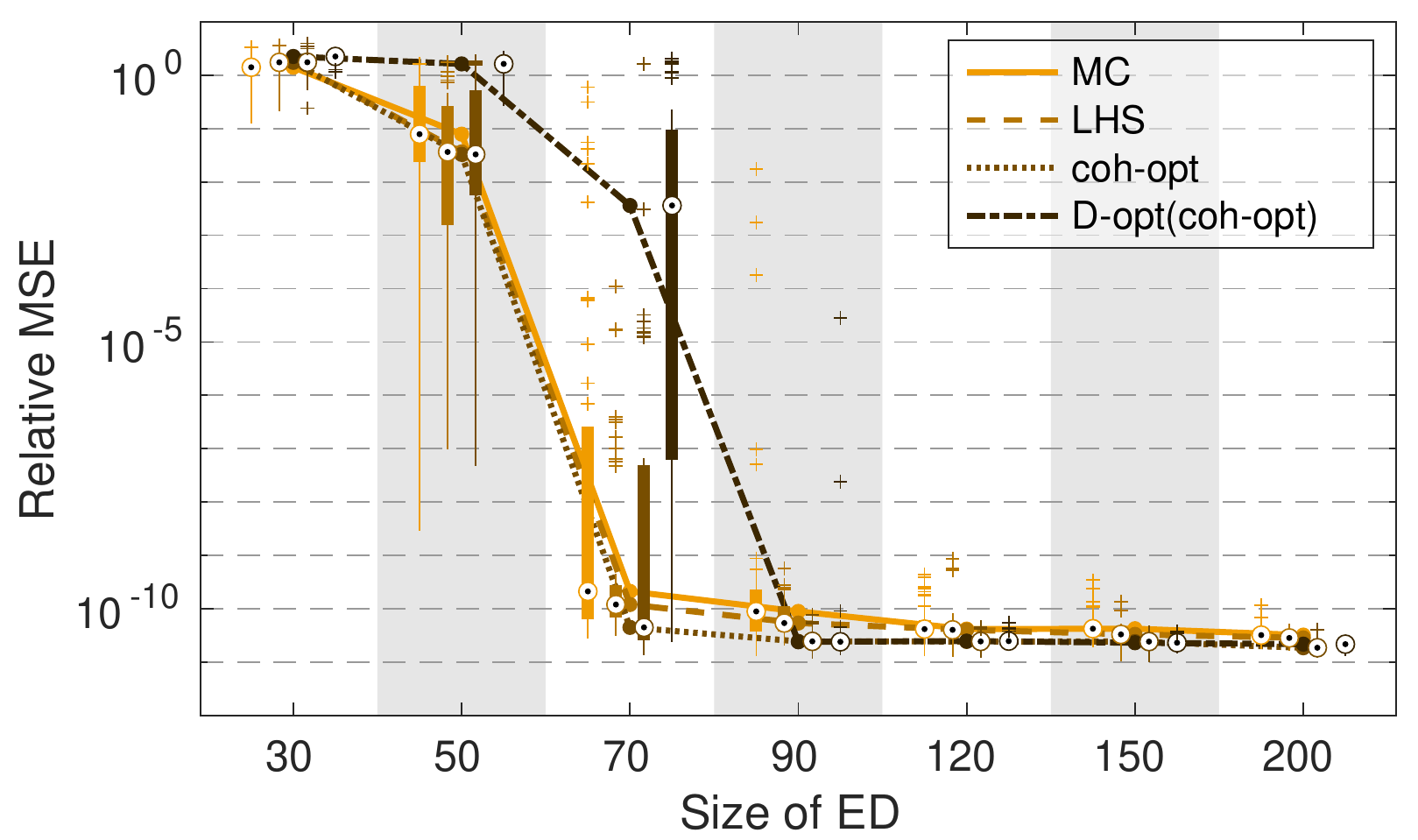}}
	\\
	\subfloat[][BCS]{\includegraphics[width=.49\textwidth, height=0.22\textheight, keepaspectratio]{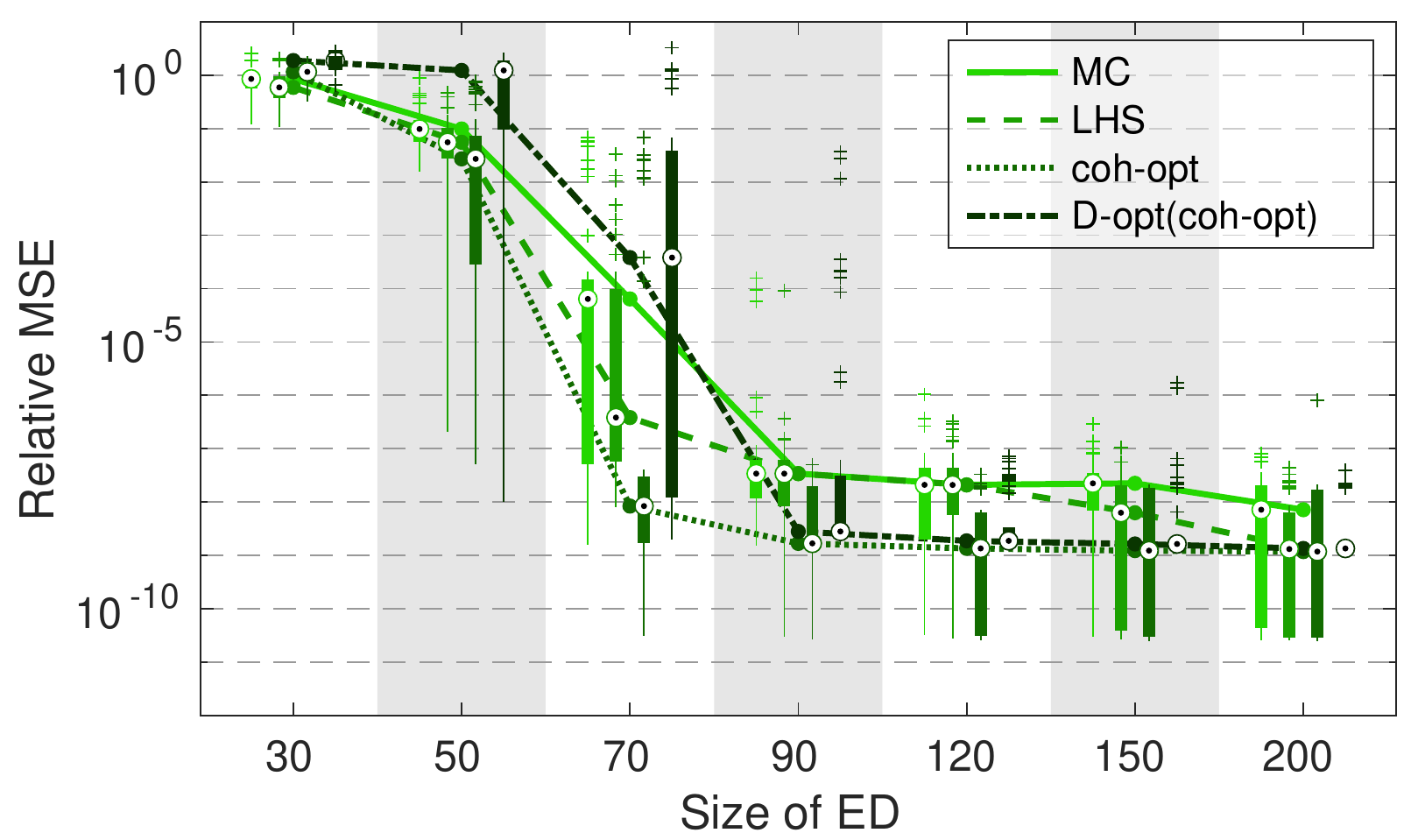}}
	\\
	\hfill
	\subfloat[][Small ED (70 points)]{\includegraphics[width=.49\textwidth, height=0.22\textheight, keepaspectratio]{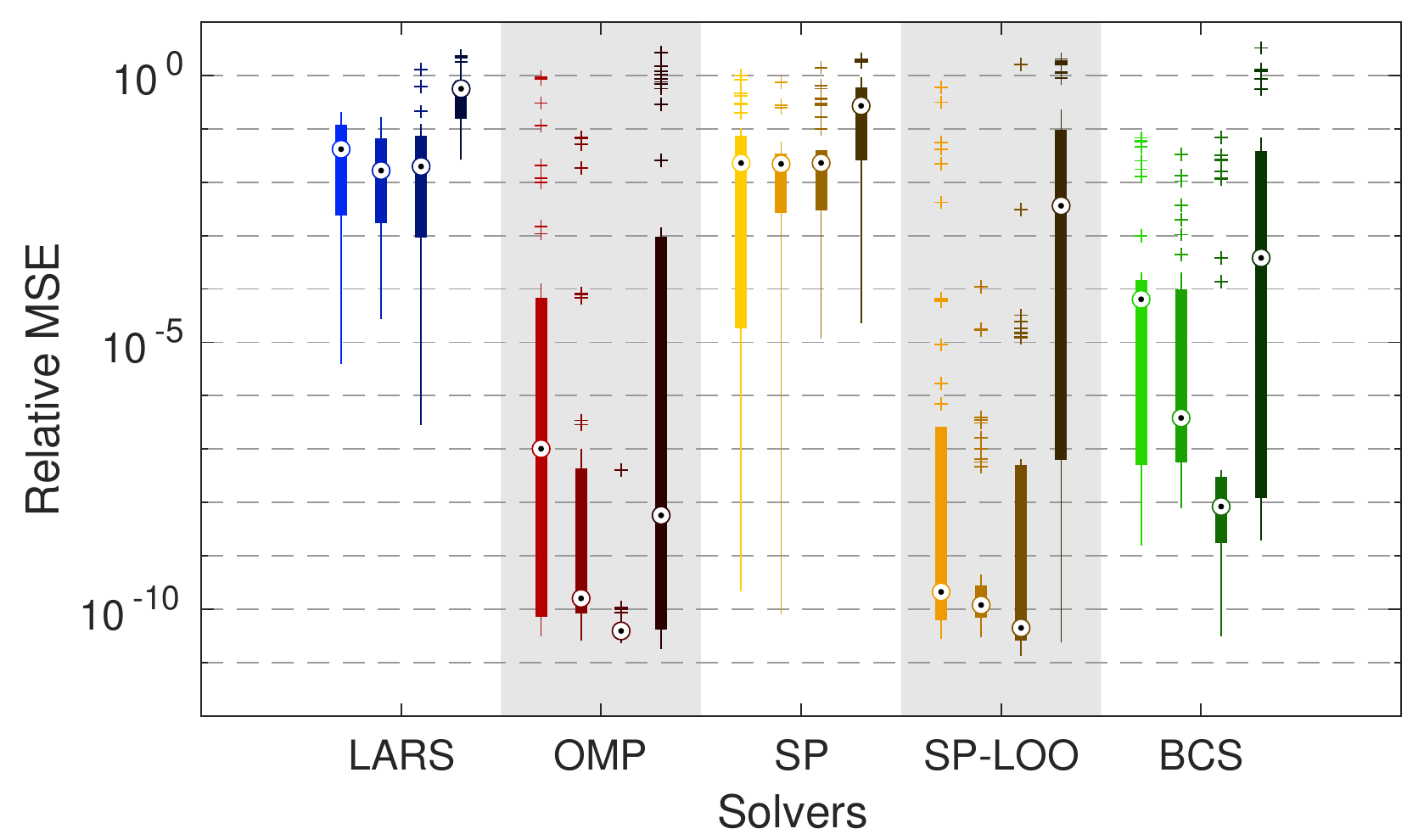}}
	\hfill
	\subfloat[][Large ED (150 points)]{\includegraphics[width=.49\textwidth, height=0.22\textheight, keepaspectratio]{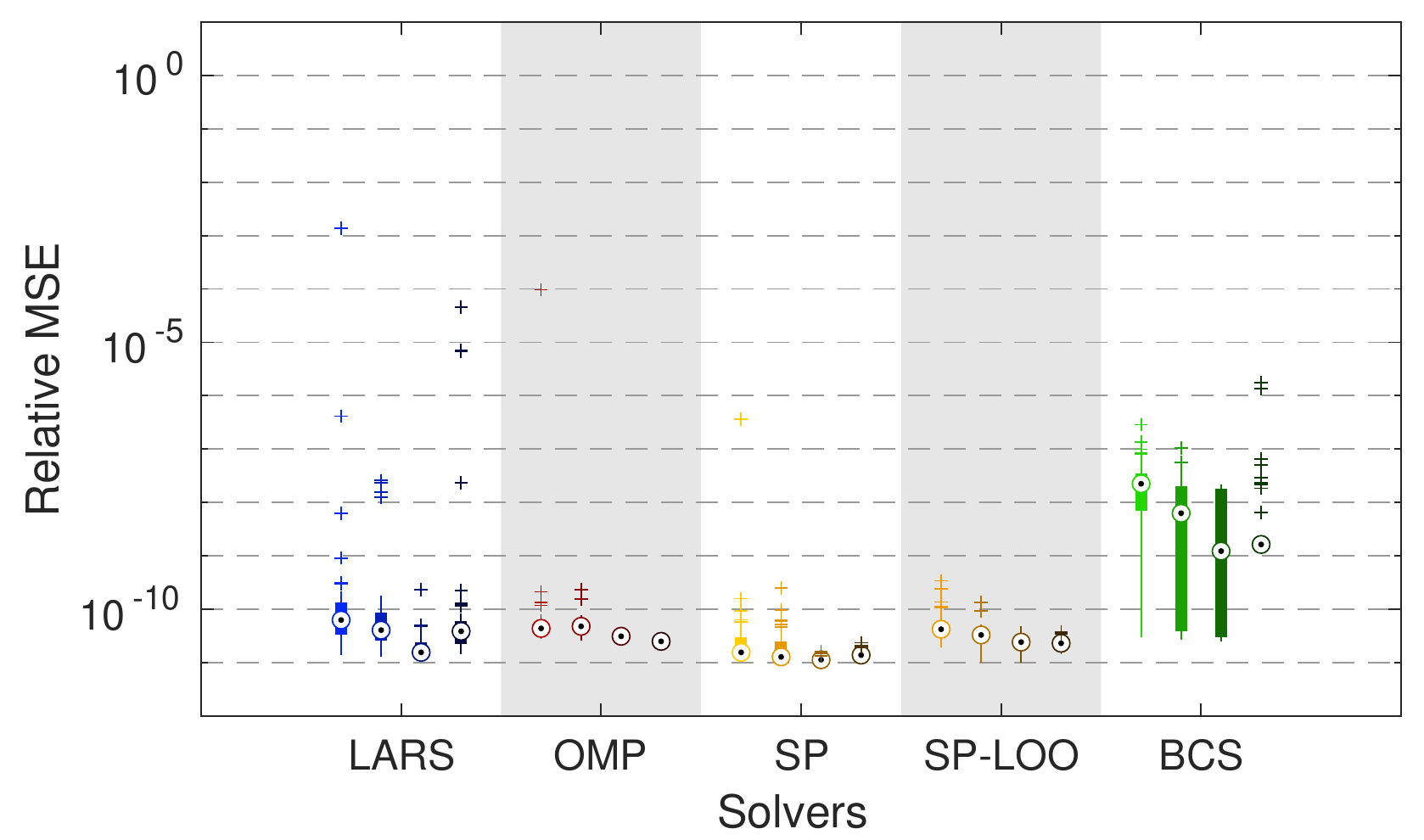}}
	\hfill
	\caption{\changed{Boxplots of relative MSE from the benchmark of five solvers and four sampling schemes for the Ishigami model ($d = 3, p = 14, q=1$). Solvers are coded by colors. Sampling schemes are shown in varying shades and line styles. In (f) and (g), we show the relative MSE of each of the solvers combined with each sampling scheme in the order MC--LHS--coh-opt--D-opt(coh-opt).}}
	\label{fig:results_ishigami_additional}	
\end{figure}

\begin{figure}[htbp]
	\centering
	\subfloat[][LARS]{\includegraphics[width=.49\textwidth, height=0.22\textheight, keepaspectratio]{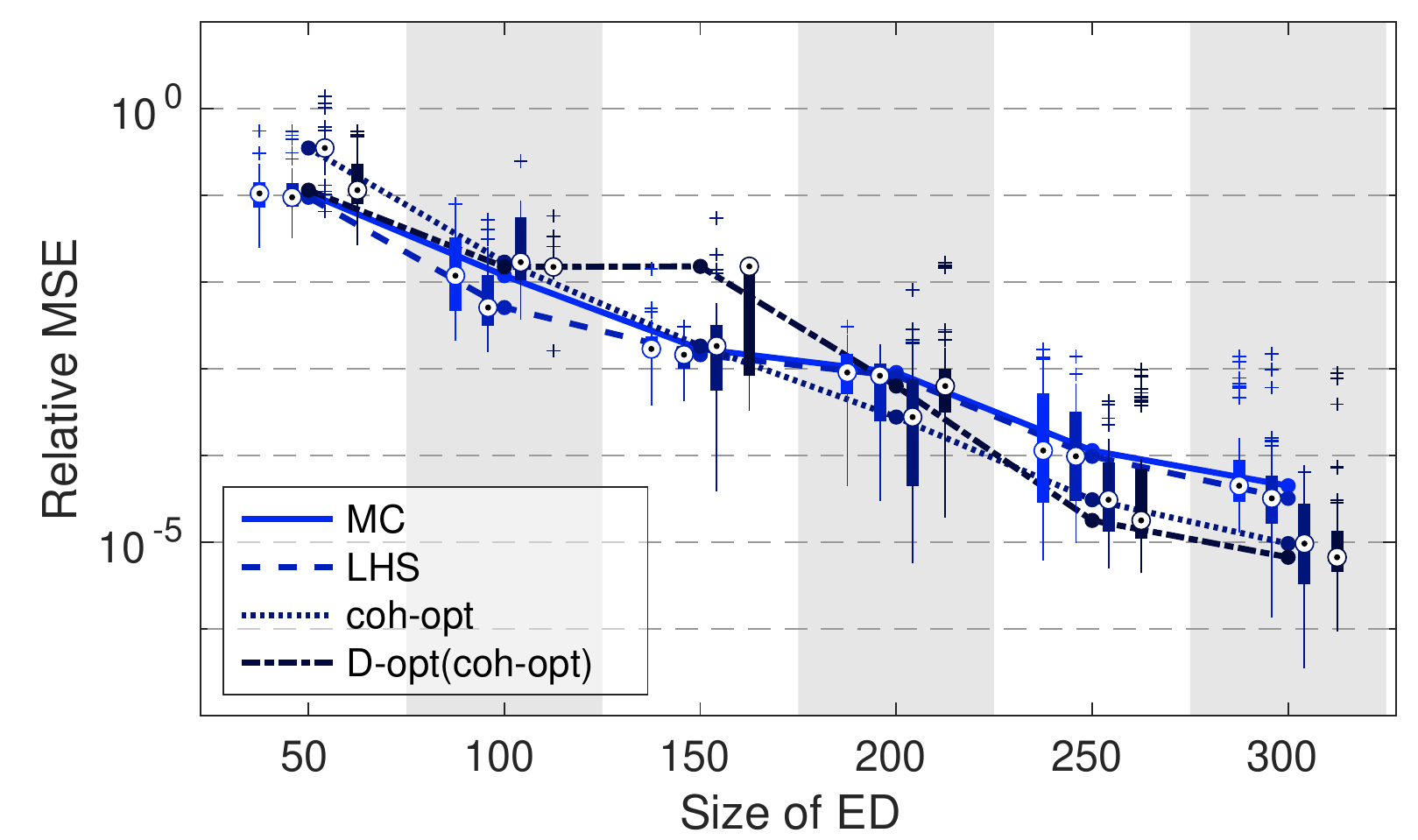}
		\label{fig:additional_results_borehole_LARS}	}
	\hfill
	\subfloat[][OMP]{\includegraphics[width=.49\textwidth, height=0.22\textheight, keepaspectratio]{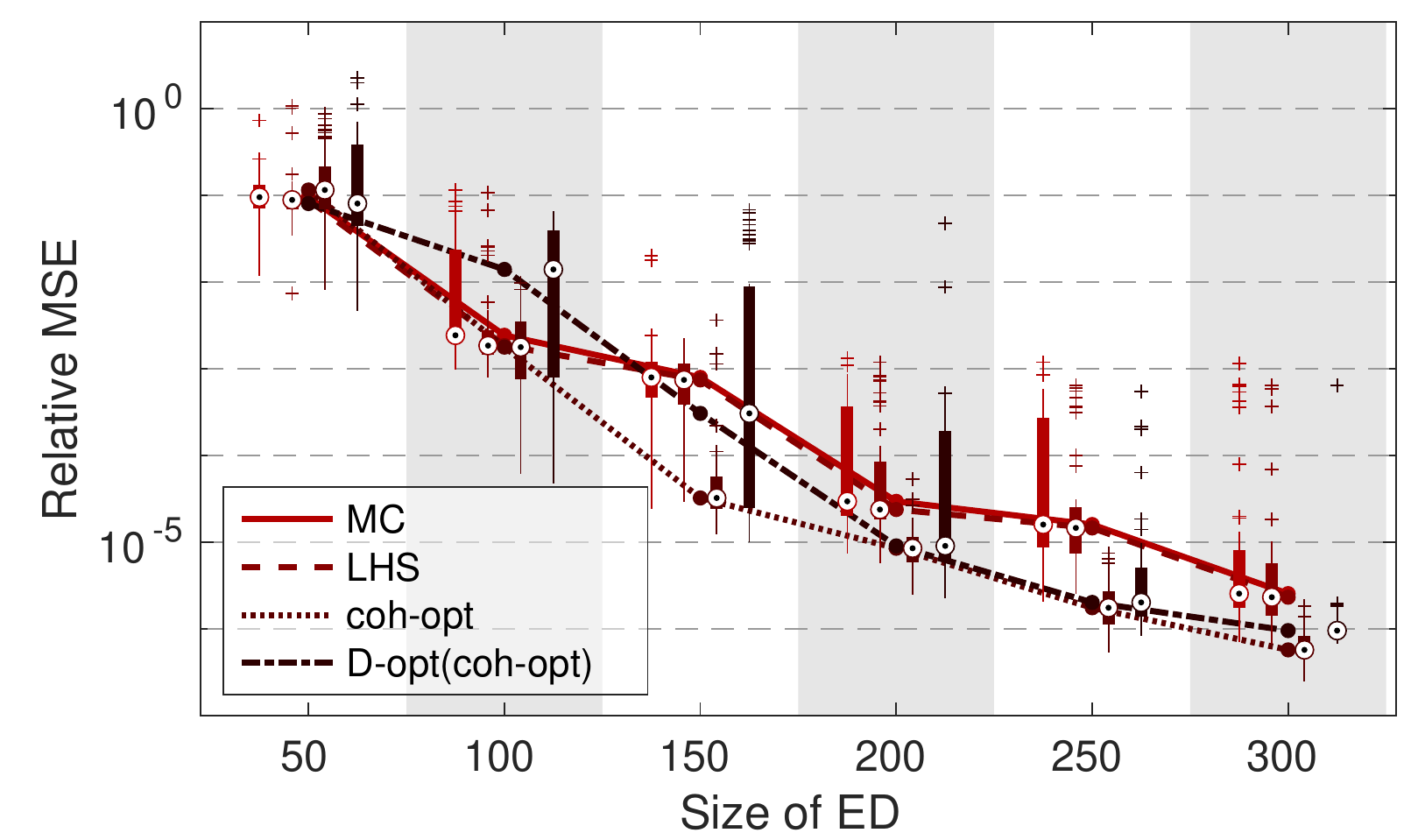}}
	\\
	\subfloat[][SP]{\includegraphics[width=.49\textwidth, height=0.22\textheight, keepaspectratio]{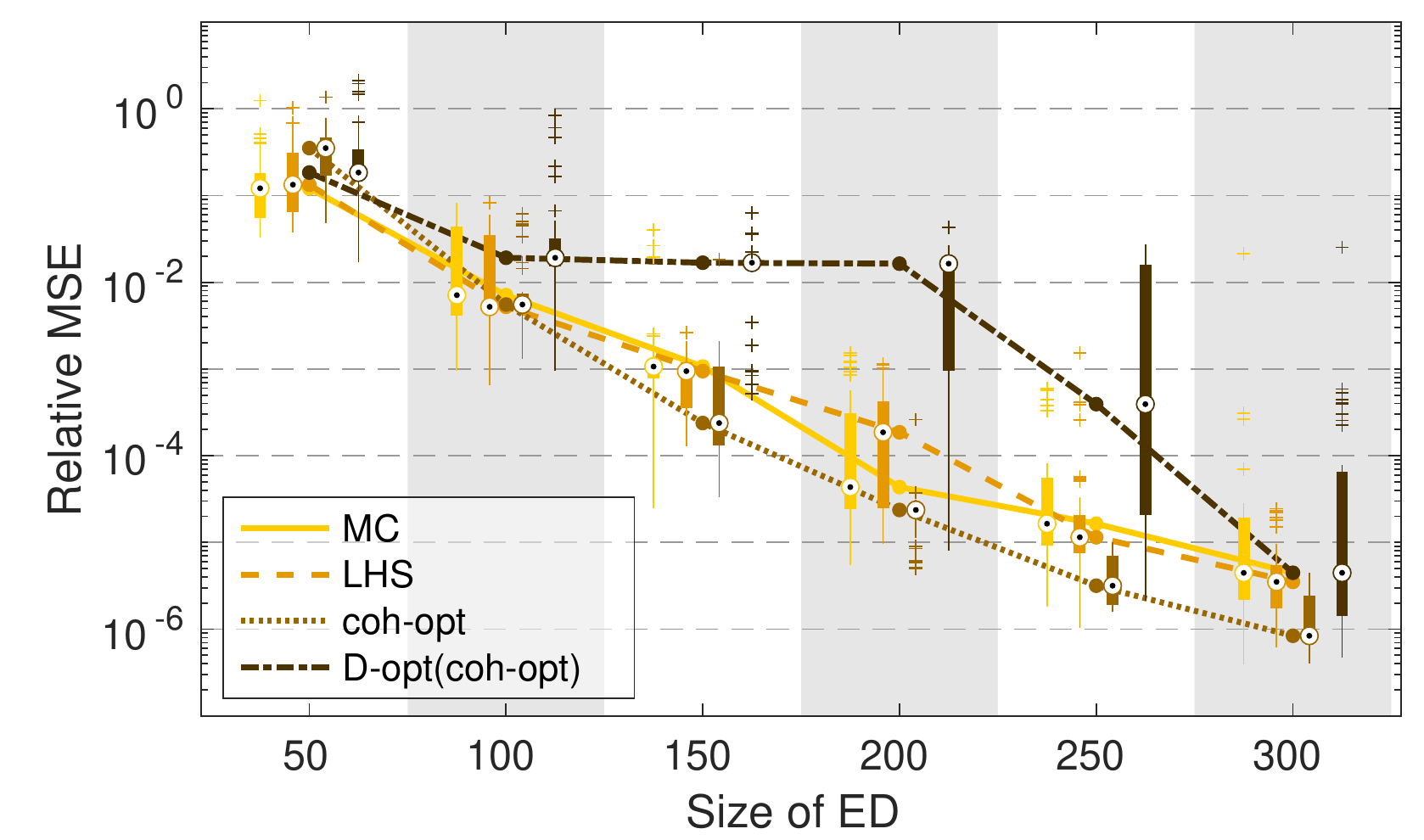}}
	\hfill
	\subfloat[][\SPloo{}]{\includegraphics[width=.49\textwidth, height=0.22\textheight, keepaspectratio]{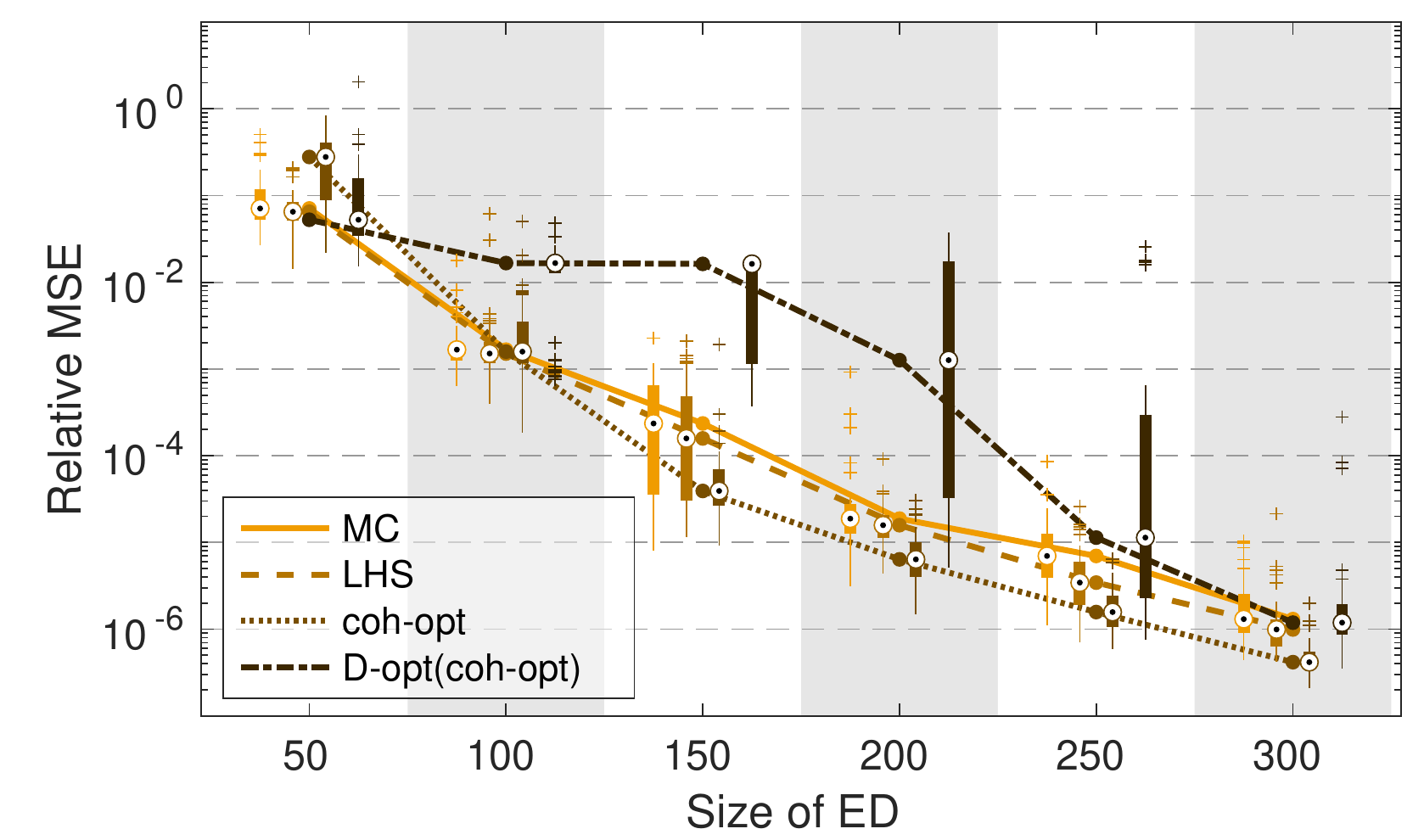}}
	\\
	\subfloat[][BCS]{\includegraphics[width=.49\textwidth, height=0.22\textheight, keepaspectratio]{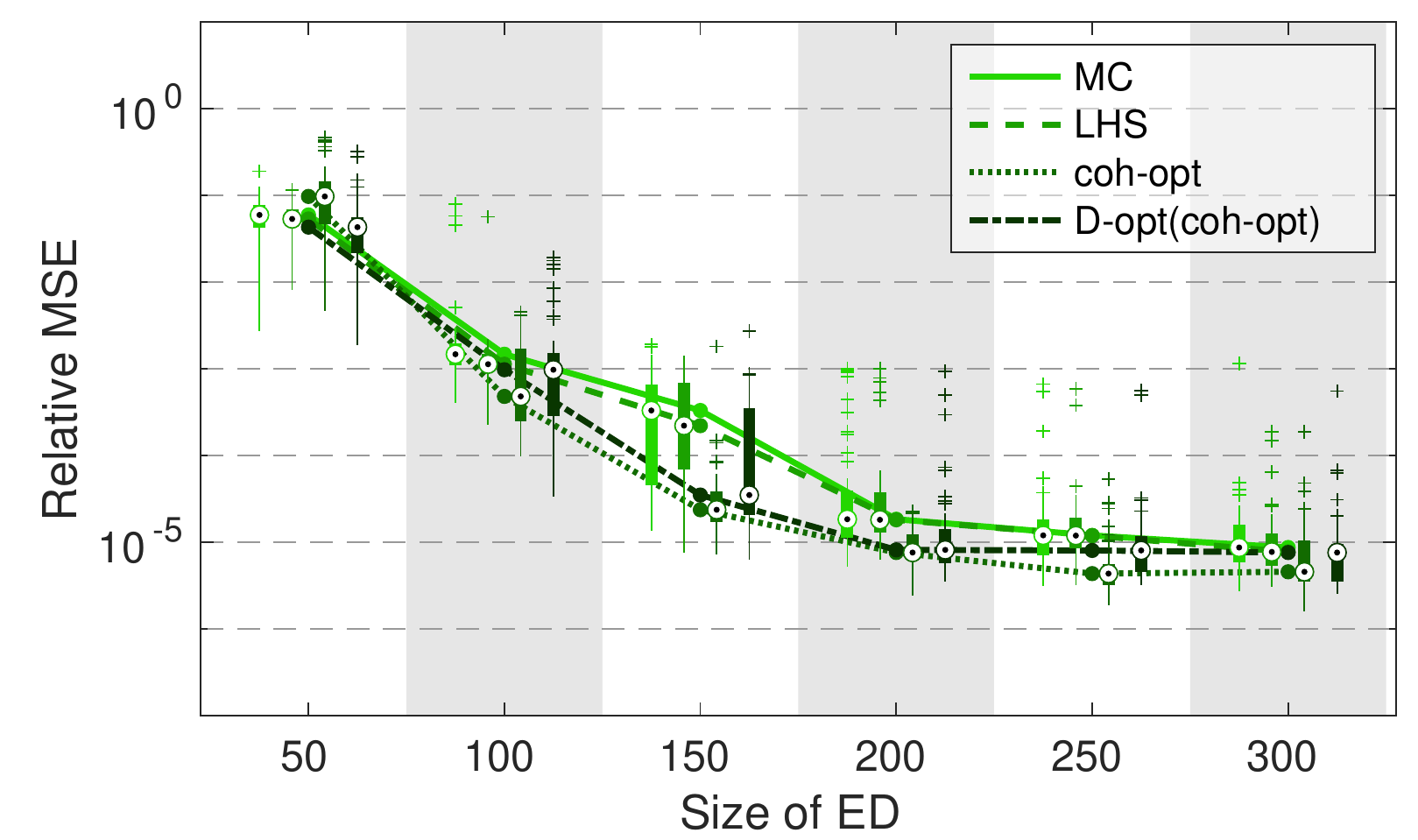}}
	\\
	\subfloat[][Small ED (100 points)]{\includegraphics[width=.49\textwidth, height=0.22\textheight, keepaspectratio]{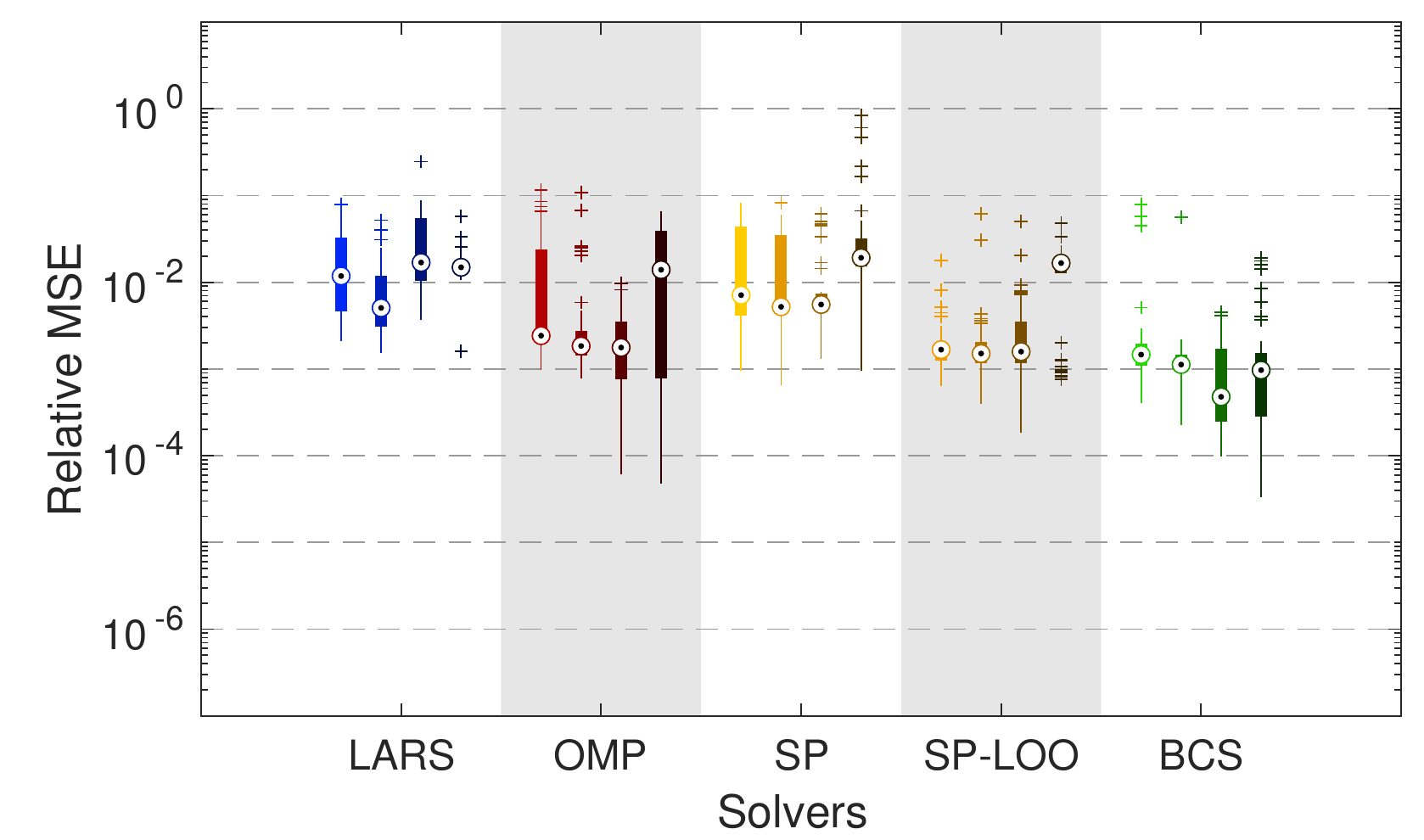}}
	\hfill
	\subfloat[][Large ED (250 points)]{\includegraphics[width=.49\textwidth, height=0.22\textheight, keepaspectratio]{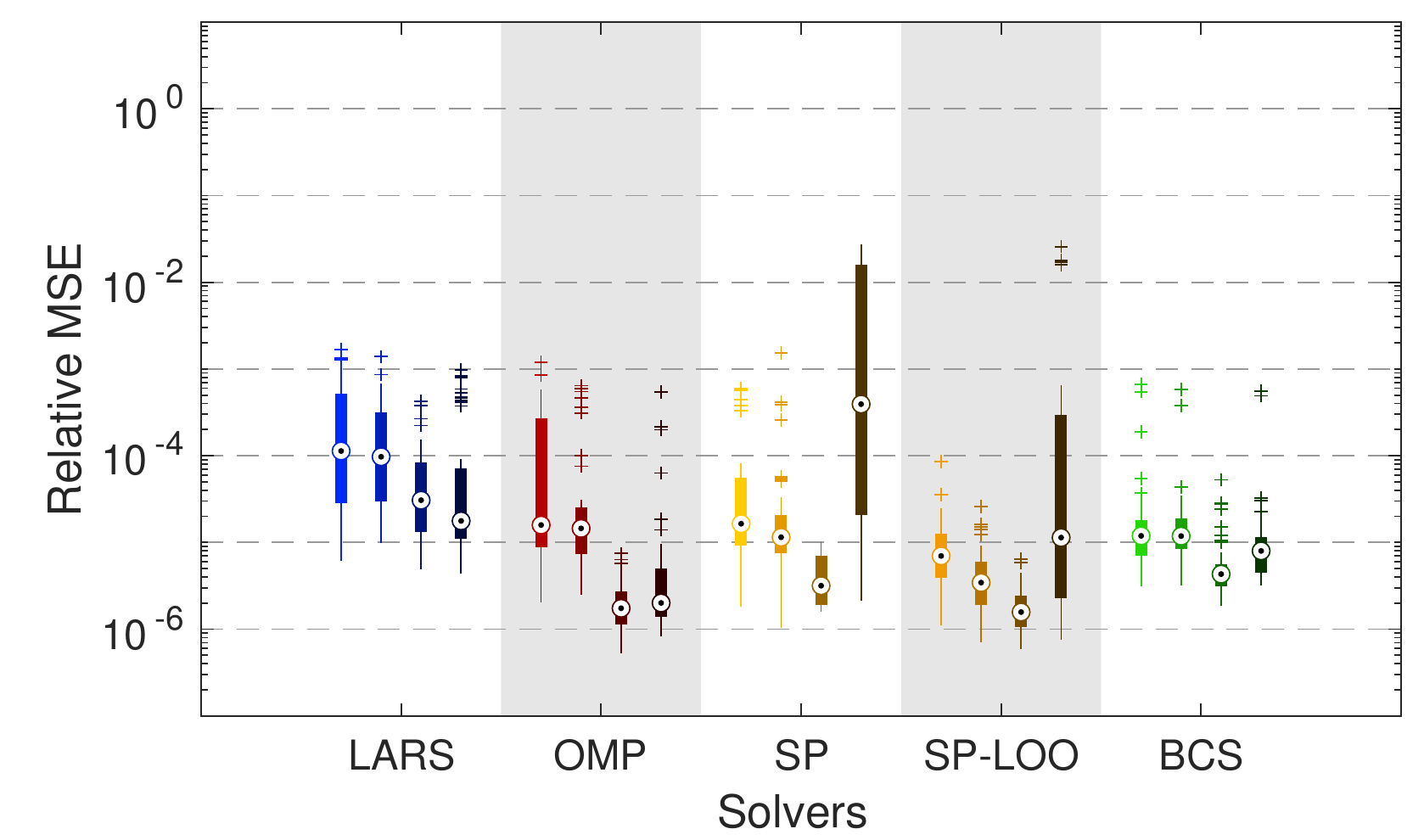}}
	\caption{\changed{Boxplots of relative MSE from the benchmark of five solvers and four sampling schemes for the borehole model ($d = 8, p = 4, q = 1$). Solvers are coded by colors. Sampling schemes are shown in varying shades and line styles. In (f) and (g), we show the relative MSE of each of the solvers combined with each sampling scheme in the order MC--LHS--coh-opt--D-opt(coh-opt). }}
	\label{fig:results_borehole_additional}	
\end{figure}

\begin{figure}[htbp]
	\centering
	\subfloat[][LARS]{\includegraphics[width=.49\textwidth, height=0.25\textheight, keepaspectratio]{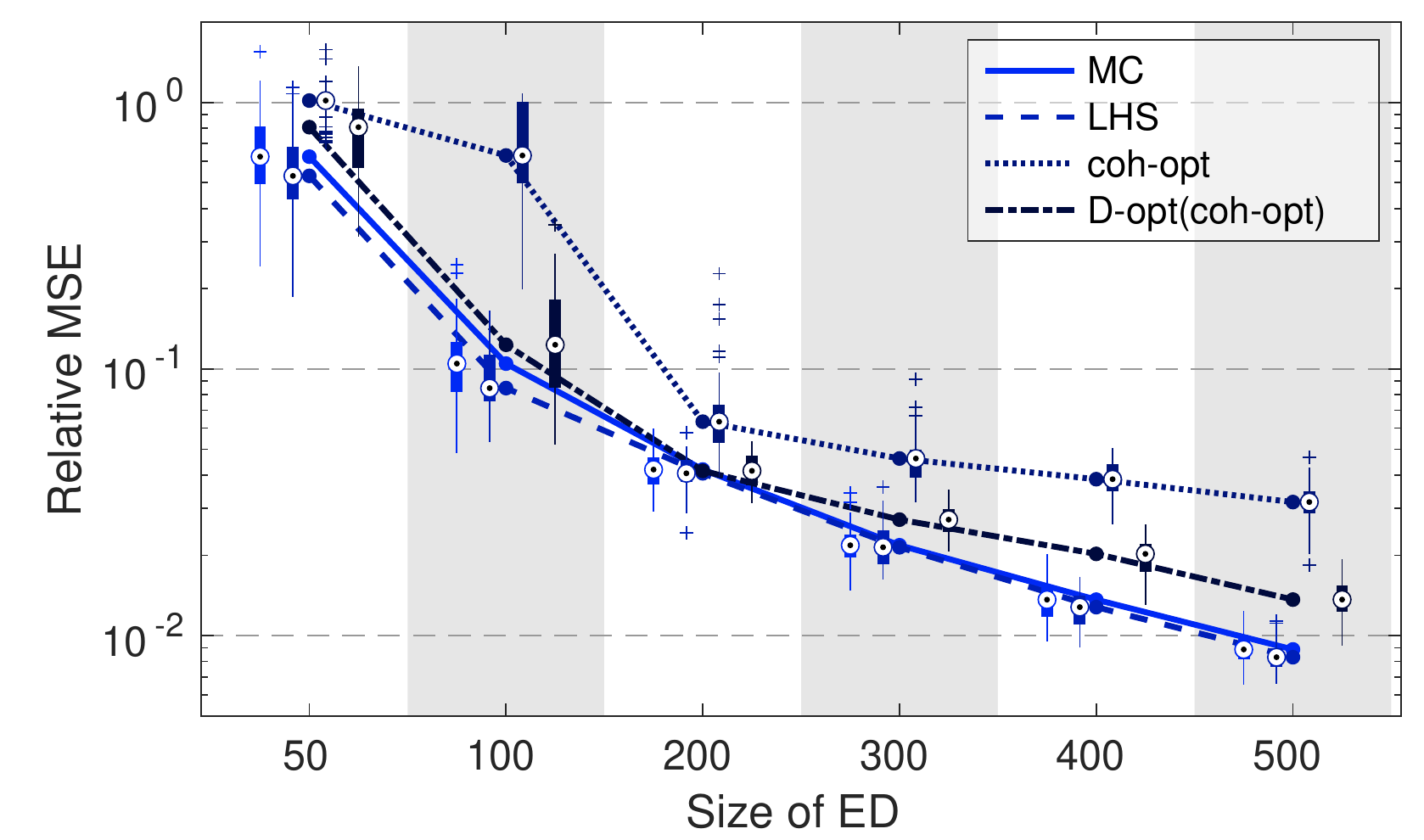}
		\label{fig:results_diffusion2D_LARS}}
	\hfill
	\subfloat[][OMP]{\includegraphics[width=.49\textwidth, height=0.25\textheight, keepaspectratio]{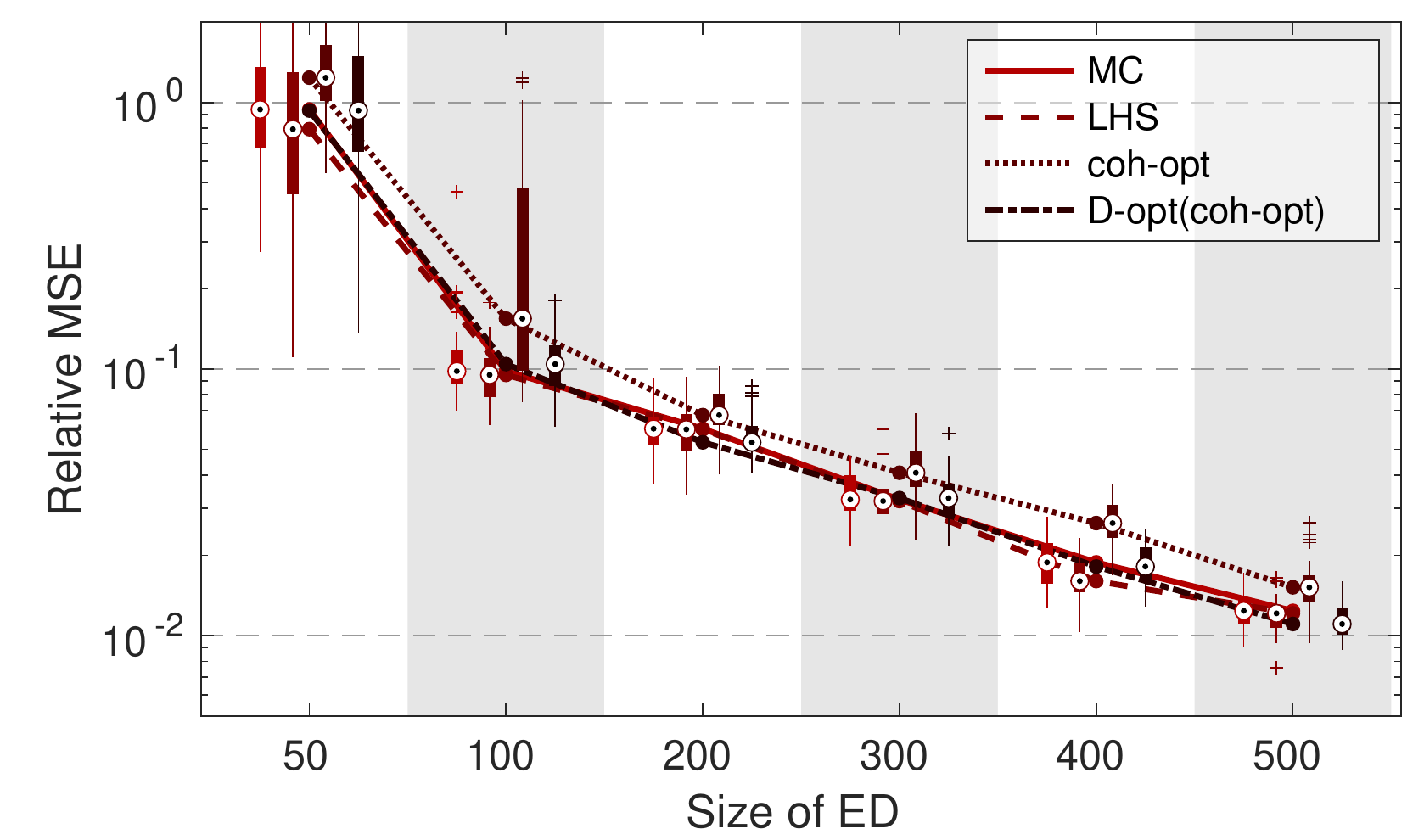}}
	\\
	\subfloat[][SP]{\includegraphics[width=.49\textwidth, height=0.25\textheight, keepaspectratio]{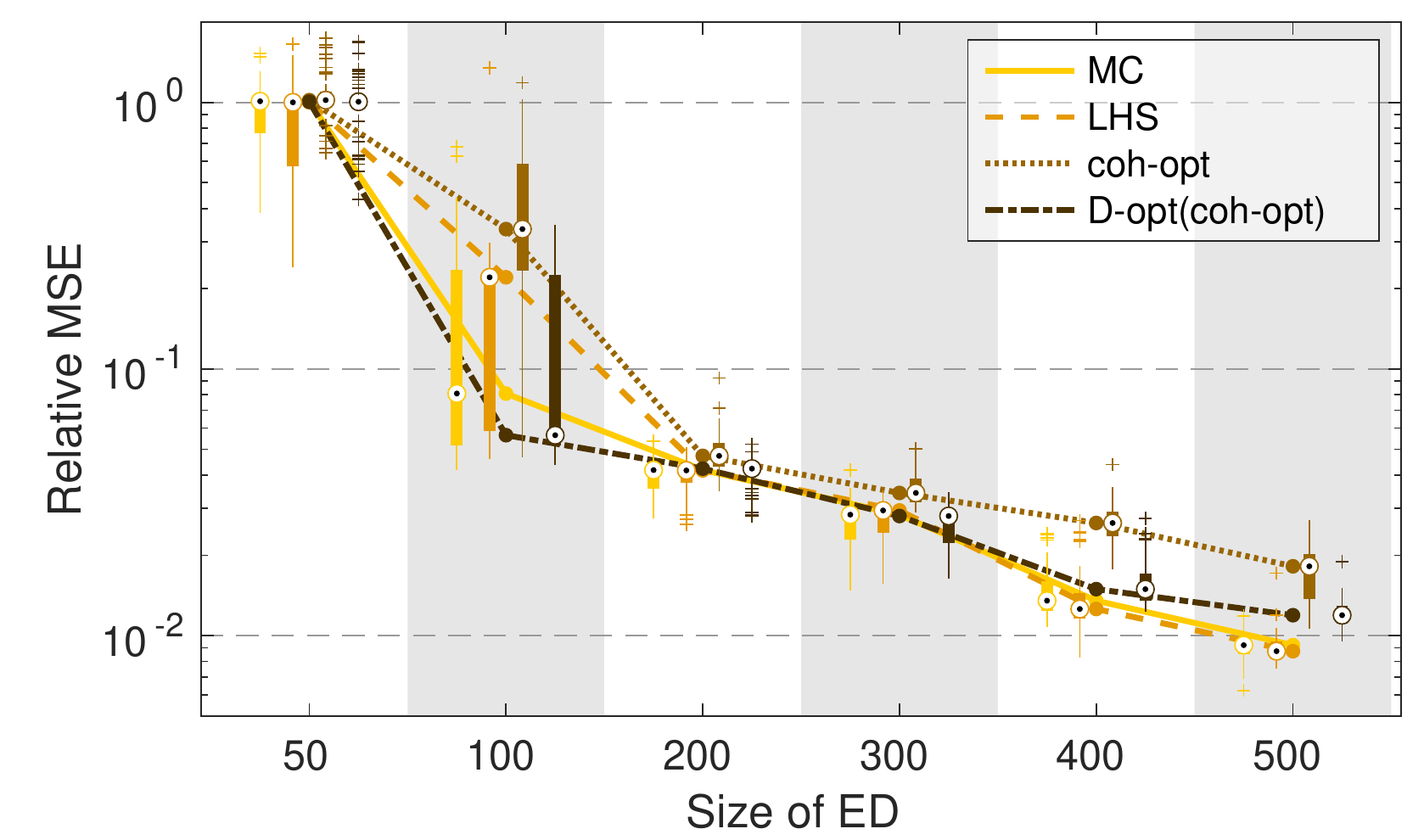}}
	\hfill
	\subfloat[][BCS]{\includegraphics[width=.49\textwidth, height=0.25\textheight, keepaspectratio]{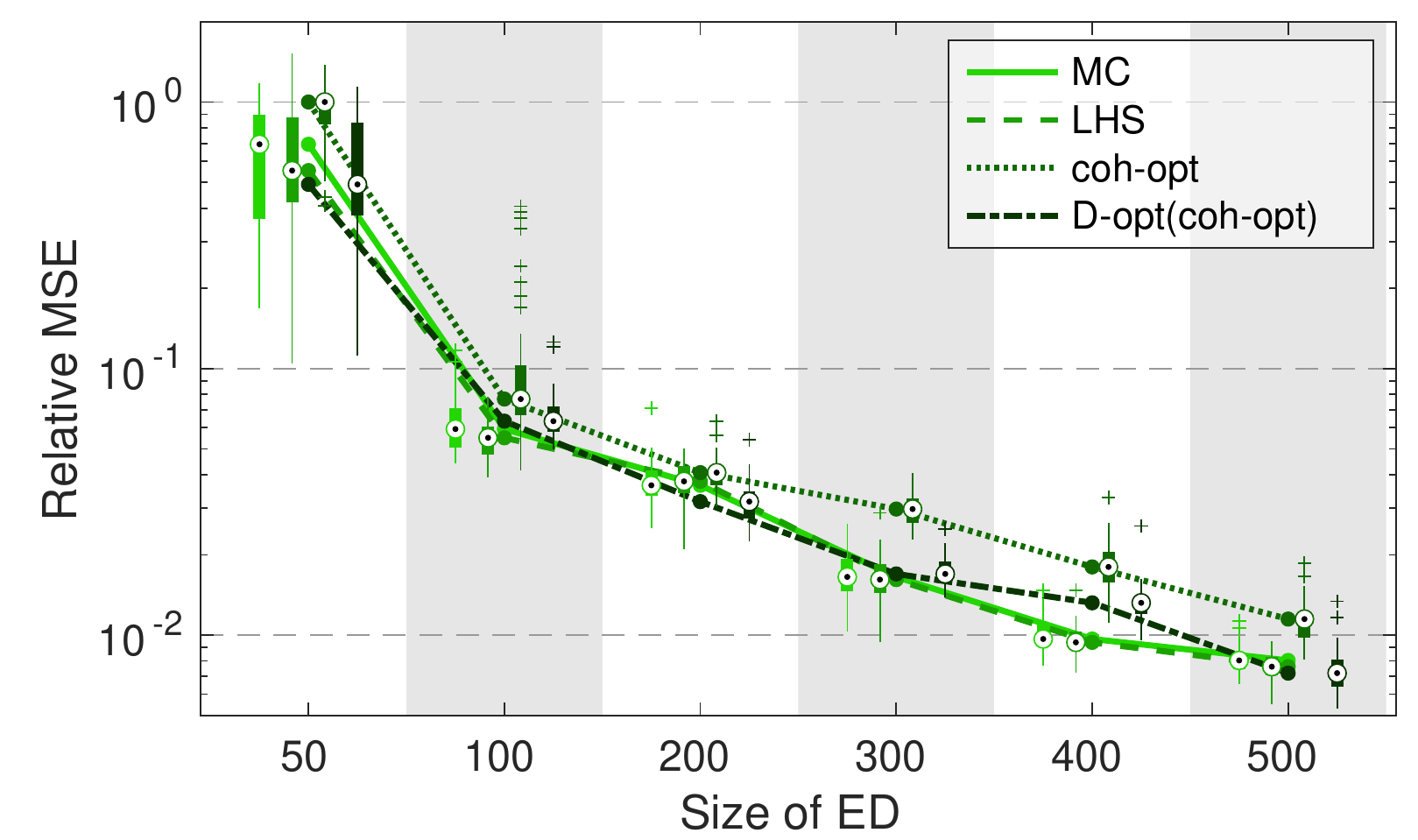}}
	\\
	\hfill
	\subfloat[][Small ED (100 points)]{\includegraphics[width=.49\textwidth, height=0.25\textheight, keepaspectratio]{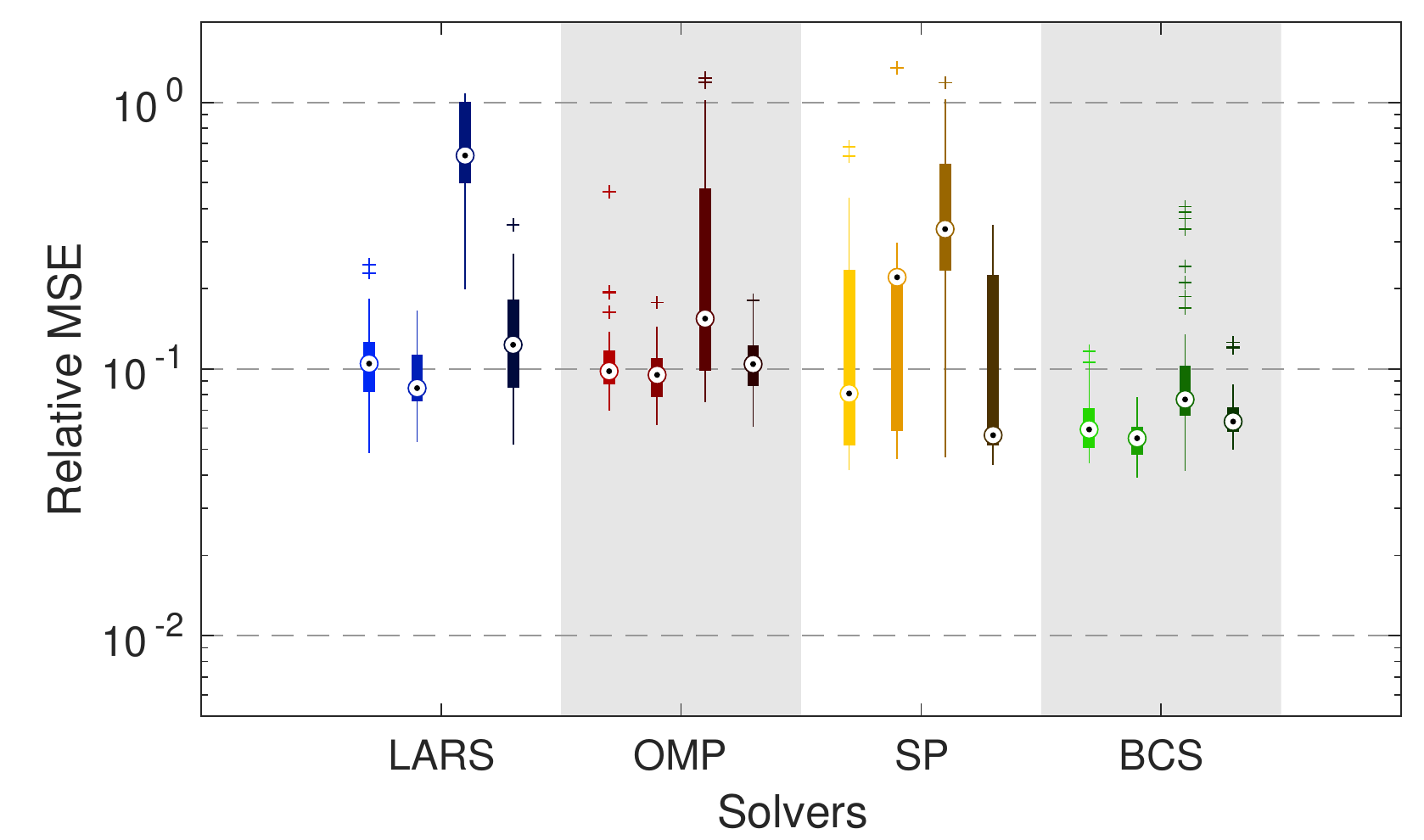}}
	\hfill
	\subfloat[][Large ED (400 points)]{\includegraphics[width=.49\textwidth, height=0.25\textheight, keepaspectratio]{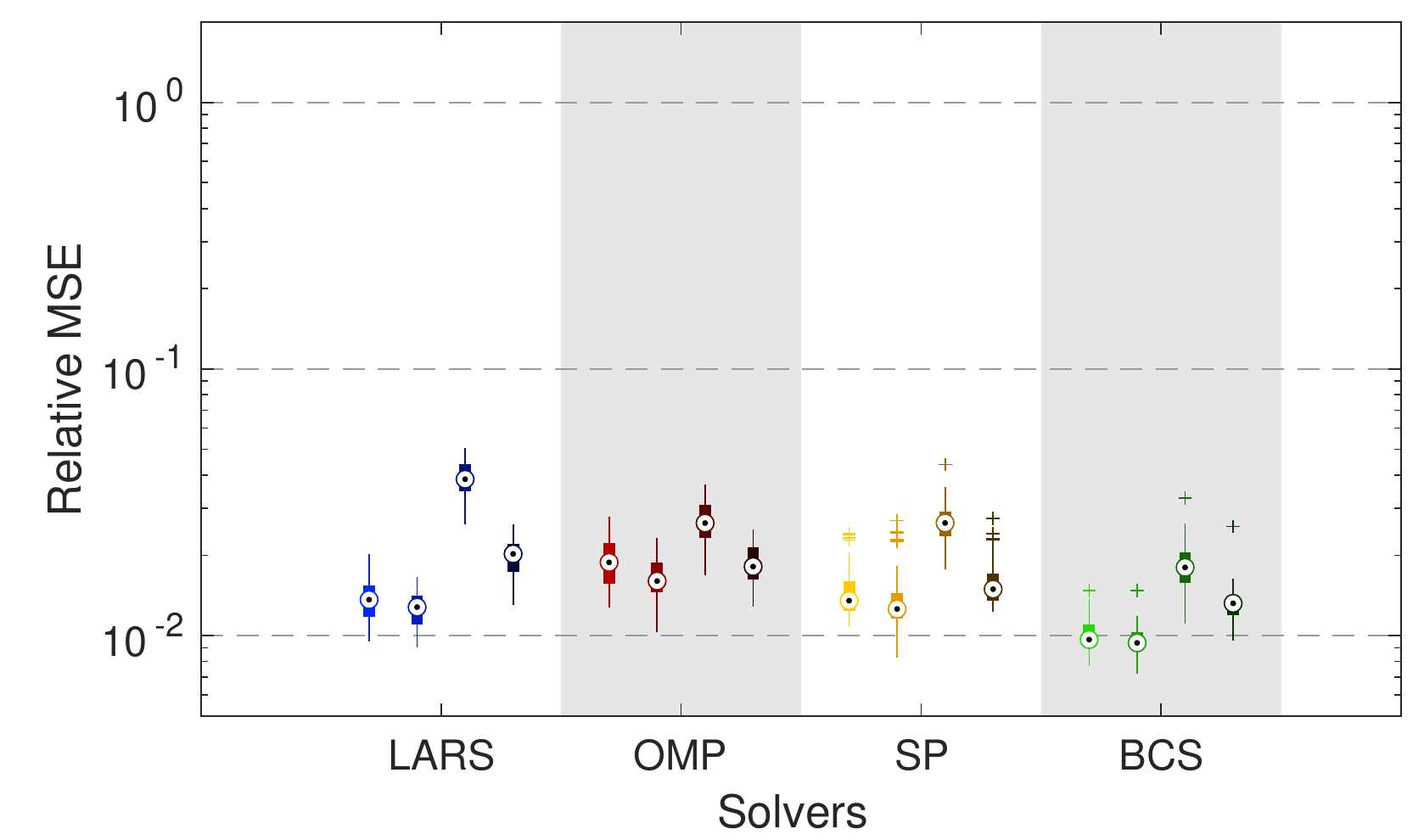}}
	\hfill
	\caption{\changed{Boxplots of relative MSE from the benchmark of four solvers and four sampling schemes for the two-dimensional diffusion model ($d = 53, p = 4, q = 0.5$). Solvers are coded by colors. Sampling schemes are shown in varying shades and line styles. In (e) and (f), we show the relative MSE of each of the solvers combined with each sampling scheme in the order MC--LHS--coh-opt--D-opt(coh-opt).}}
	\label{fig:results_diffusion2D_additional}	
\end{figure}

\begin{figure}[htbp]
	\centering
	\subfloat[][LARS]{\includegraphics[width=.49\textwidth, height=0.25\textheight, keepaspectratio]{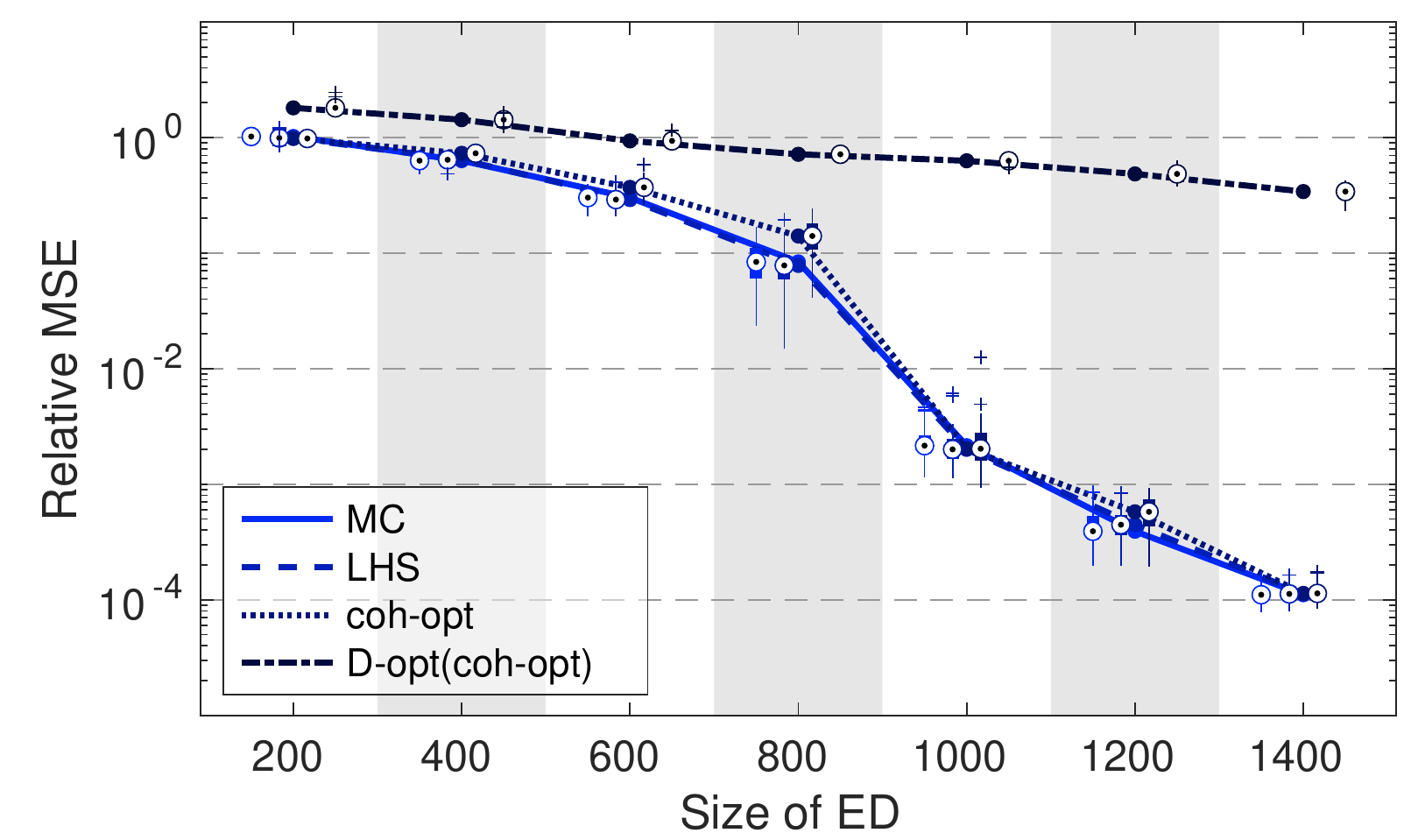}
		\label{fig:additional_results_highdimfct_LARS}	}
	\hfill
	\subfloat[][OMP]{\includegraphics[width=.49\textwidth, height=0.25\textheight, keepaspectratio]{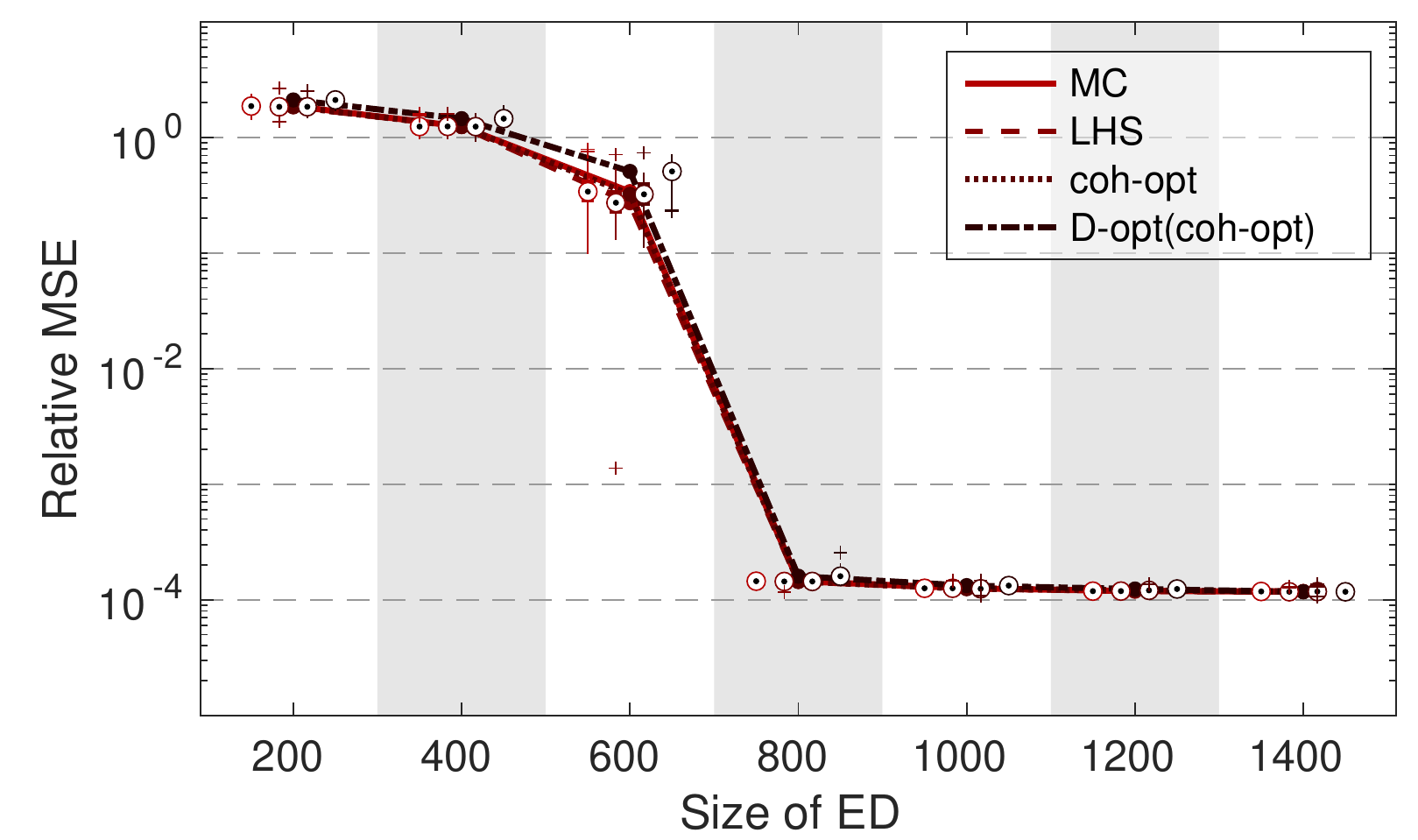}}
	\\
	\subfloat[][SP]{\includegraphics[width=.49\textwidth, height=0.25\textheight, keepaspectratio]{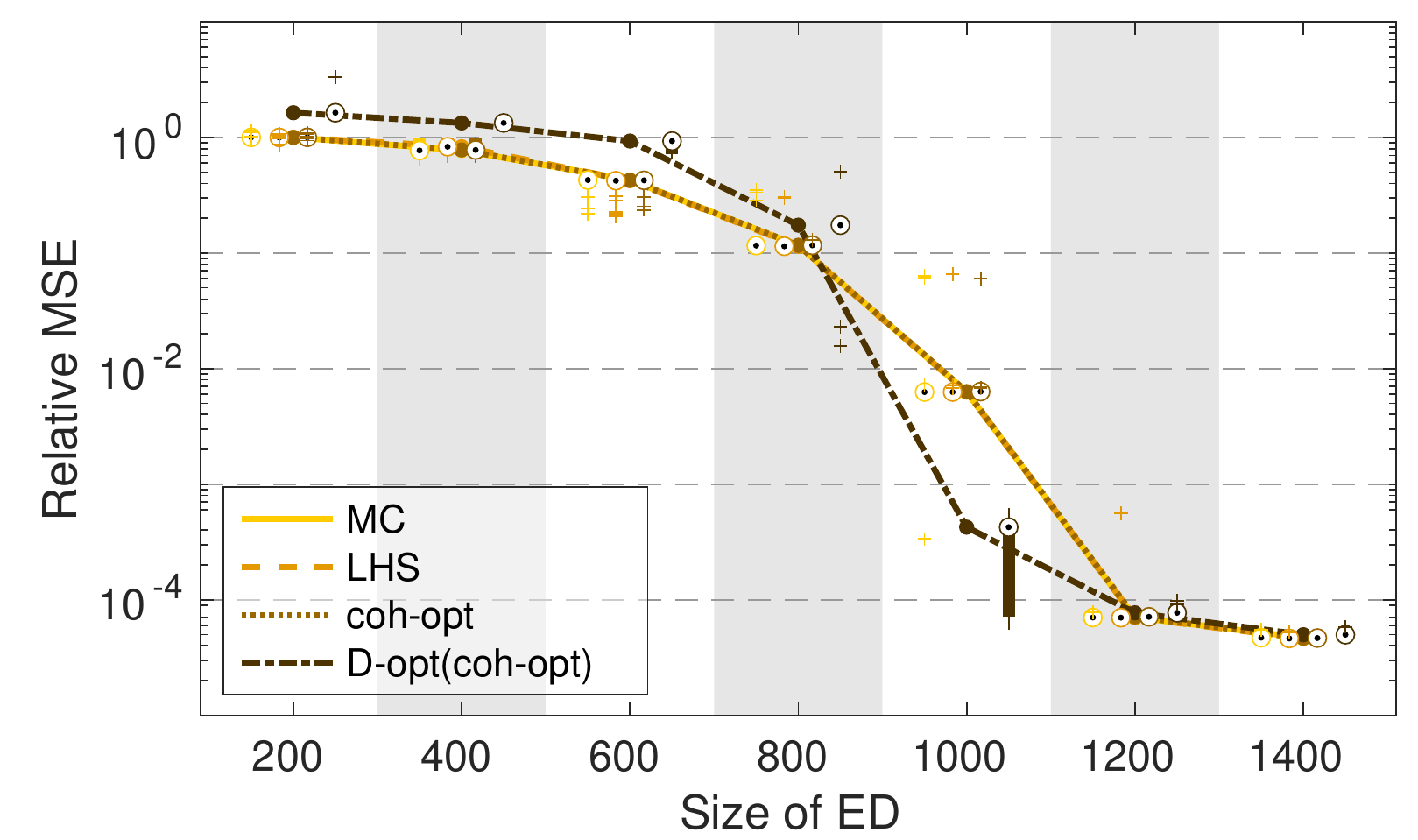}}
	\hfill
	\subfloat[][BCS]{\includegraphics[width=.49\textwidth, height=0.25\textheight, keepaspectratio]{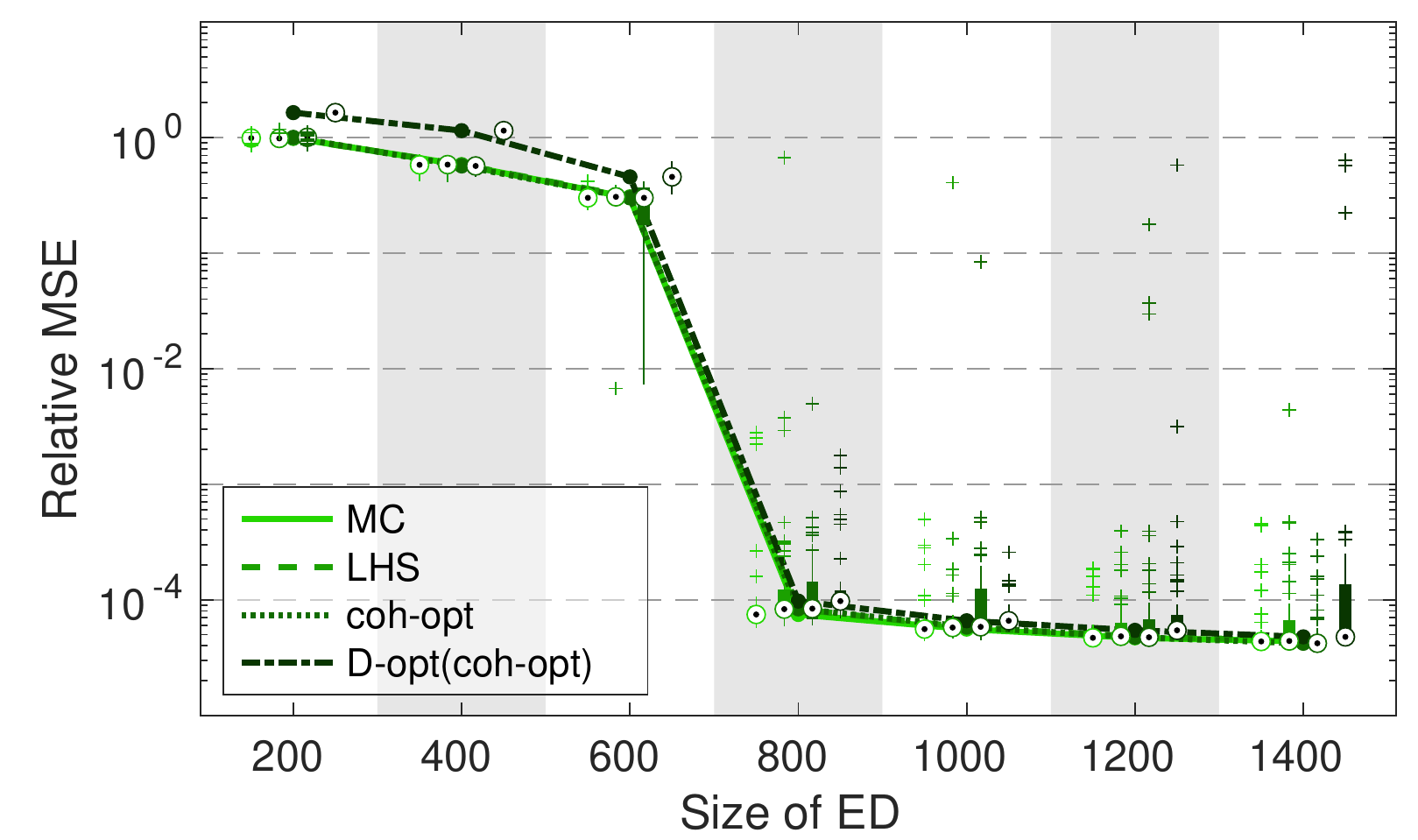}}
	\\
	\hfill
	\subfloat[][Small ED (400 points)]{\includegraphics[width=.49\textwidth, height=0.25\textheight, keepaspectratio]{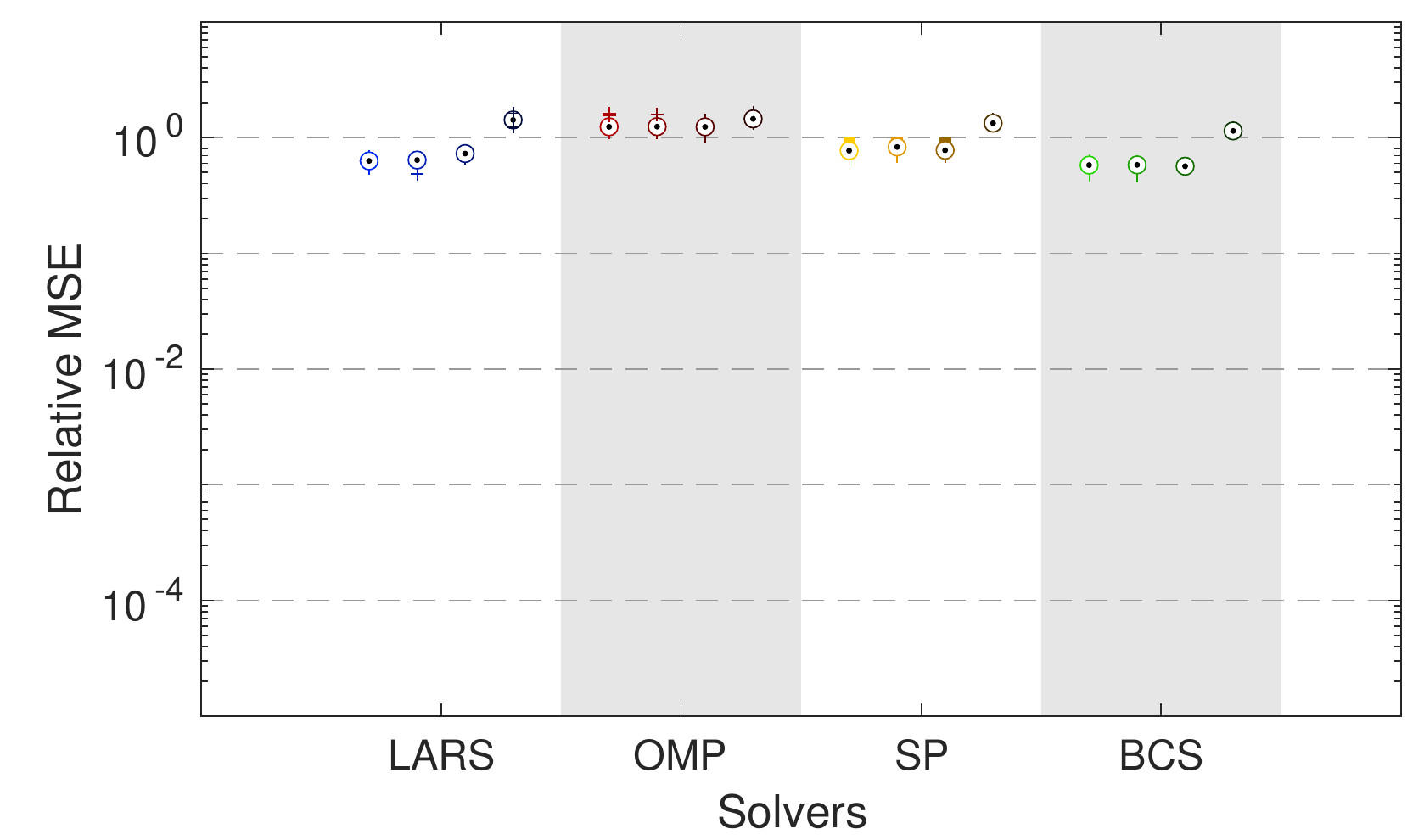}}
	\hfill
	\subfloat[][Large ED (1200 points)]{\includegraphics[width=.49\textwidth, height=0.25\textheight, keepaspectratio]{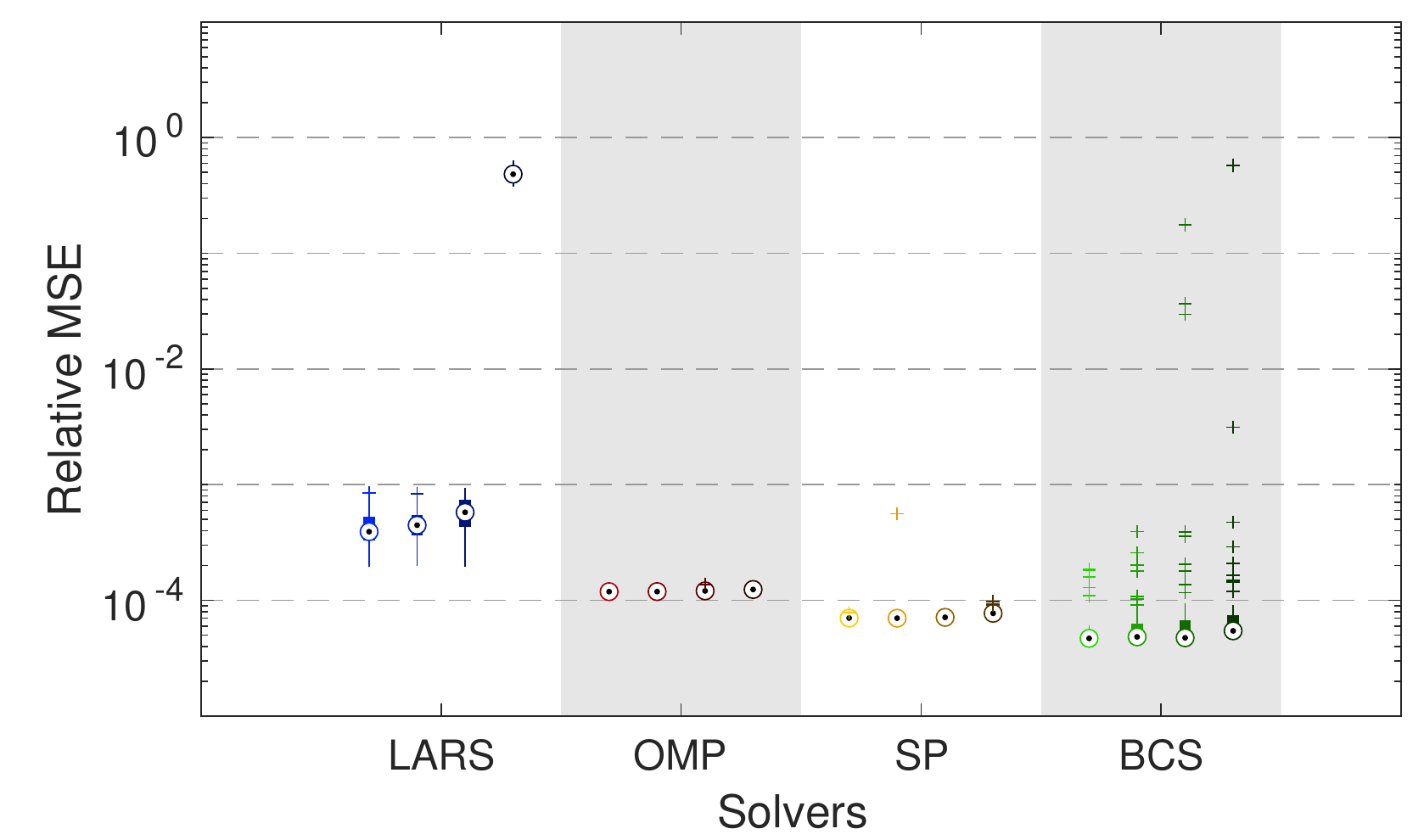}}
	\hfill
	\caption{\changed{Boxplots of relative MSE from the benchmark of four solvers and four sampling schemes for the \changed{100D} function ($d = 100, p = 4, q = 0.5$). Solvers are coded by colors. Sampling schemes are shown in varying shades and line styles. In (e) and (f), we show the relative MSE of each of the solvers combined with each sampling scheme in the order MC--LHS--coh-opt--D-opt(coh-opt). }}
	\label{fig:results_highdimfct_additional}	
\end{figure}

\changed{In Section~\ref{sec:results_sampling}, Figure~\ref{fig:results_sampling_aggregated} we showed aggregated results for the benchmark of solvers and sampling schemes.
	To give a more tangible impression of the data, in Figures~\ref{fig:results_ishigami_additional}--\ref{fig:results_highdimfct_additional} we display the boxplots of relative MSE against ED size for the four models Ishigami, borehole, two-dimensional diffusion, and \changed{100D} function. We show all combinations of solvers and sampling schemes, resulting in 16--20 combinations. Solvers are denoted by different colors. Sampling schemes are shown in varying shades and line styles. 
	We also show the same results sliced at small and large ED sizes to compare the performance between solvers.
}


\subsection{\changed{Comparison of sampling schemes together with solvers, using a smaller candidate basis}}
\begin{figure}[htbp]
	\centering
	\subfloat[][Ishigami model, LARS]{\includegraphics[width=.49\textwidth, height=0.24\textheight, keepaspectratio]{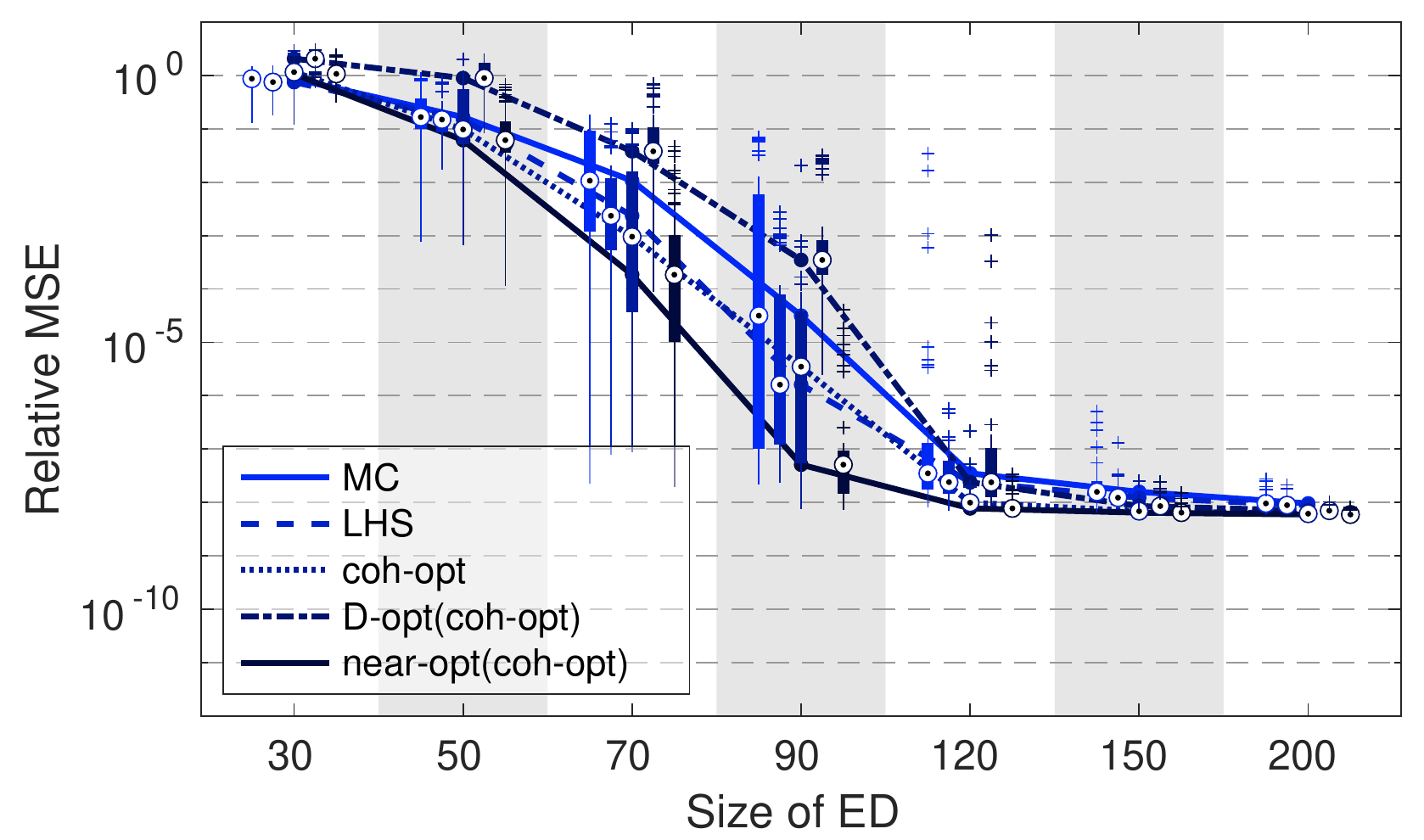}}
	\hfill
	\subfloat[][Borehole model, LARS]{\includegraphics[width=.49\textwidth, height=0.24\textheight, keepaspectratio]{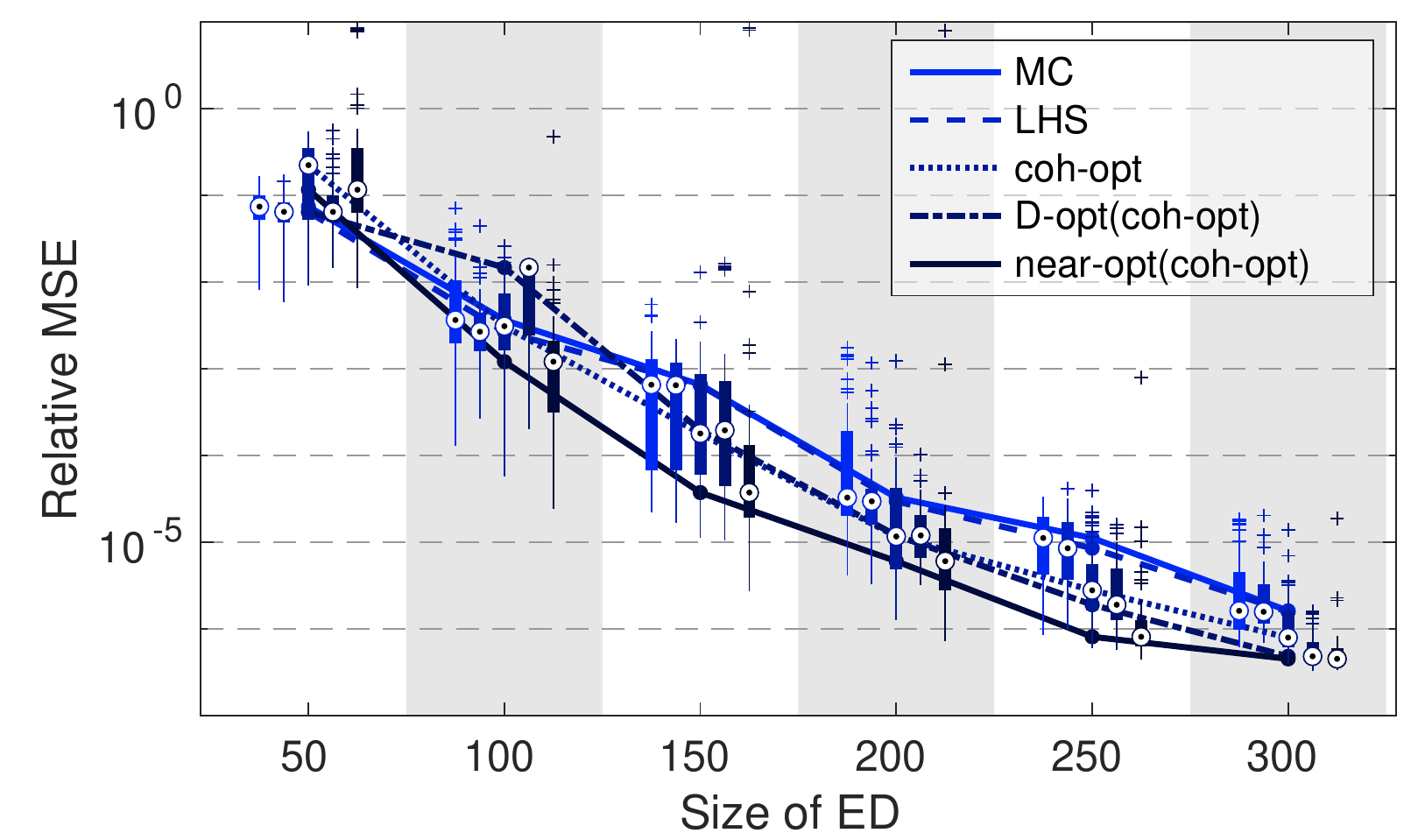}}
	\\
	\subfloat[][Ishigami model, SP]{\includegraphics[width=.49\textwidth, height=0.24\textheight, keepaspectratio]{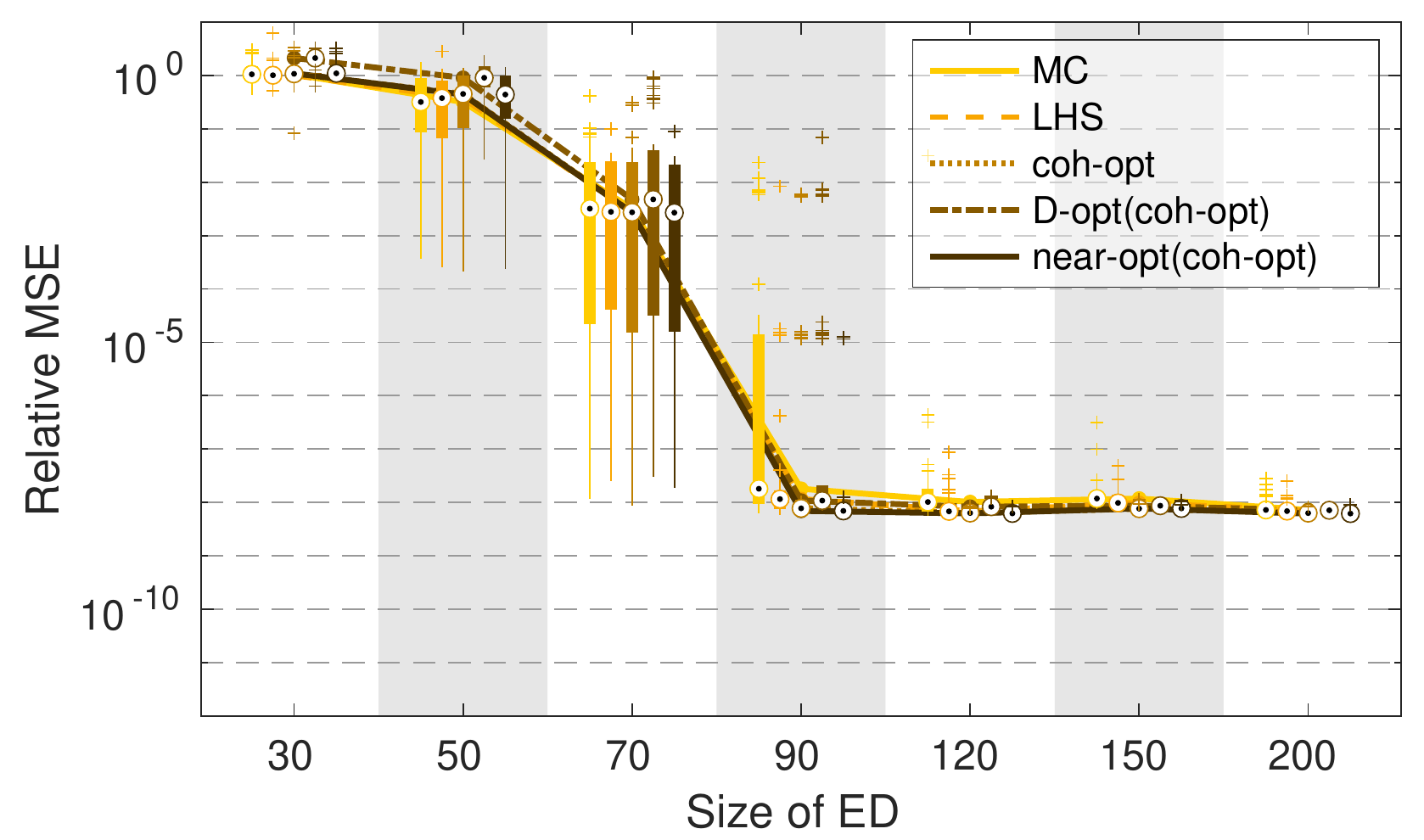}}
	\hfill
	\subfloat[][Borehole model, SP]{\includegraphics[width=.49\textwidth, height=0.24\textheight, keepaspectratio]{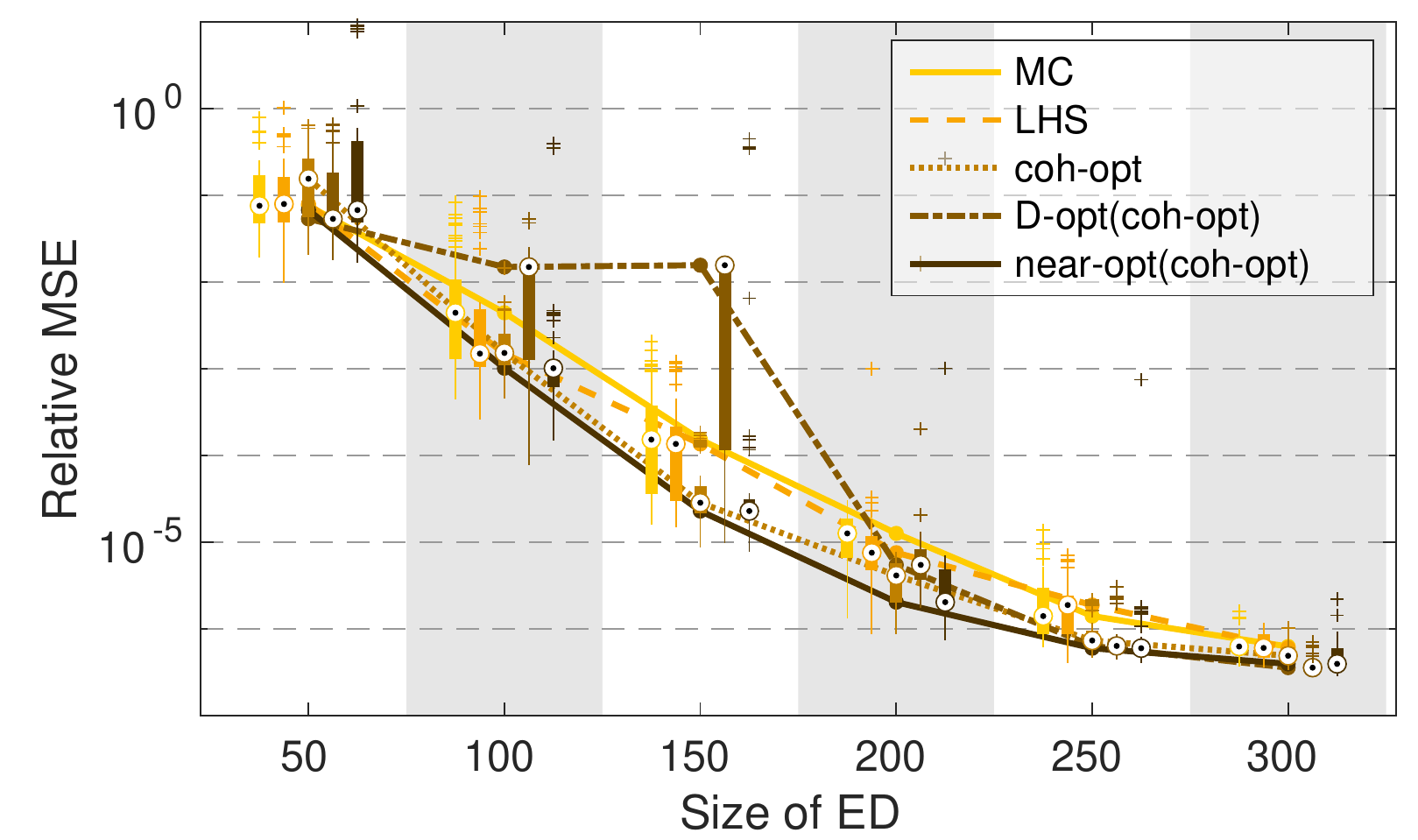}}
	\\
	\subfloat[][Ishigami model, BCS]{\includegraphics[width=.49\textwidth, height=0.24\textheight, keepaspectratio]{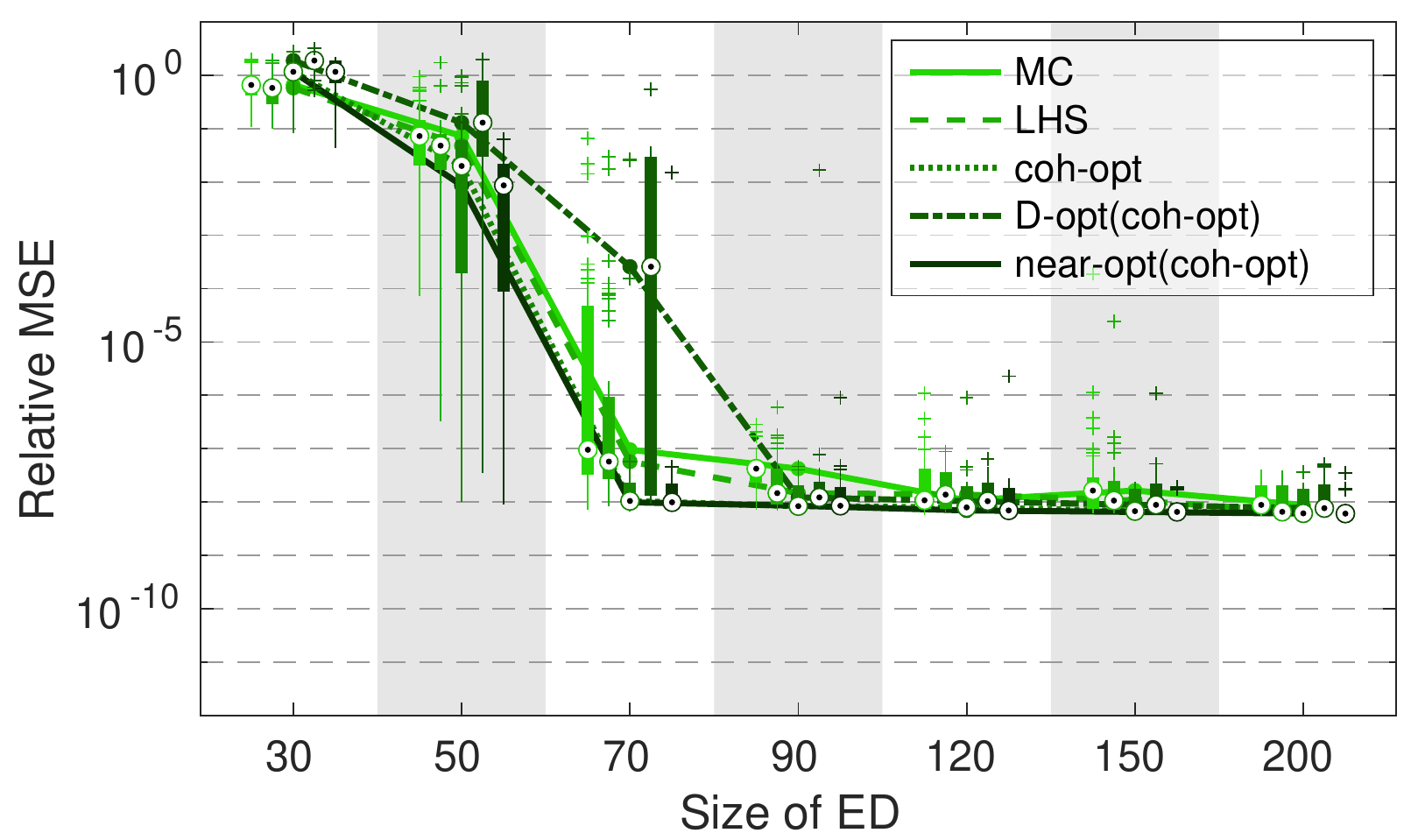}}
	\hfill
	\subfloat[][Borehole model, BCS]{\includegraphics[width=.49\textwidth, height=0.24\textheight, keepaspectratio]{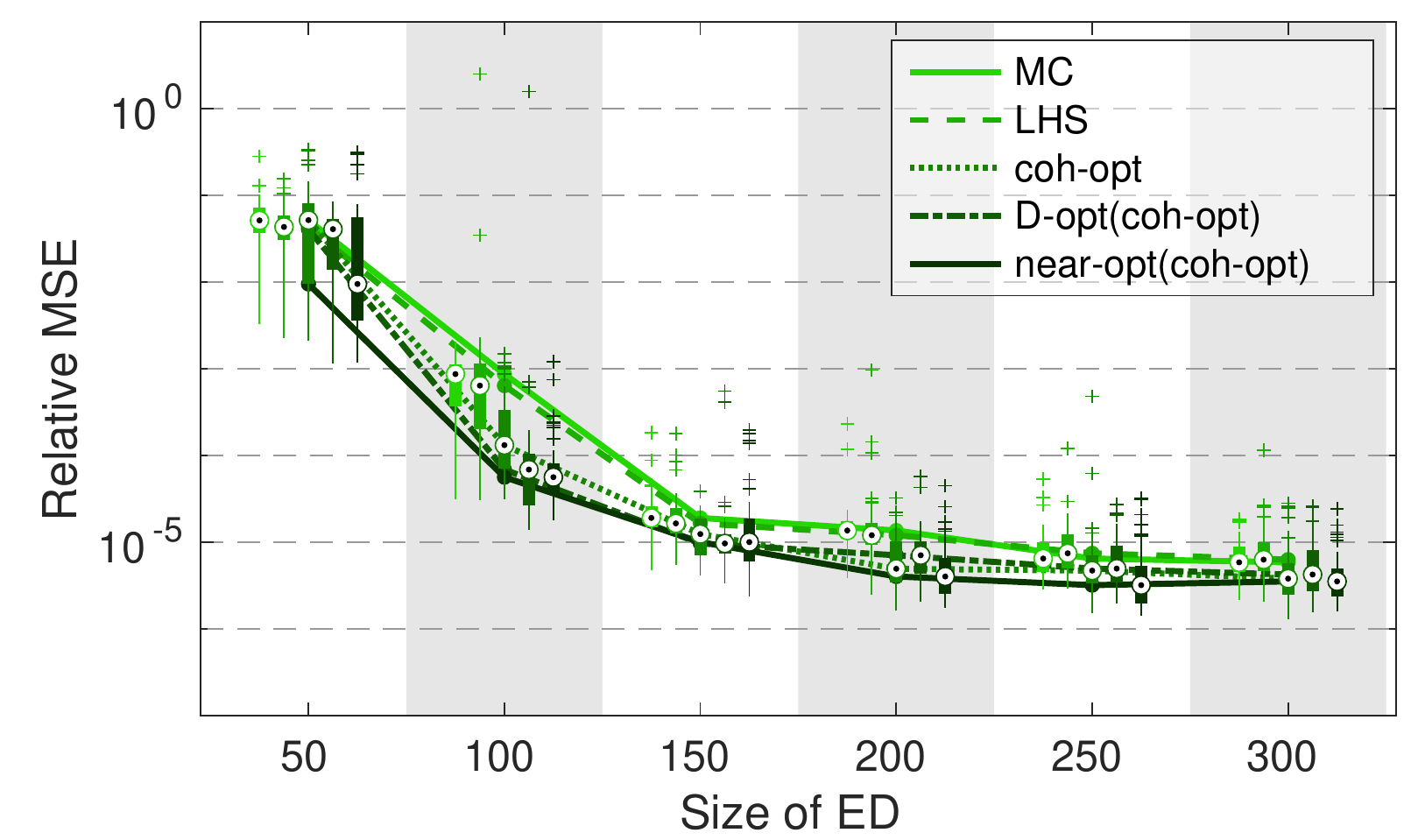}}
	\caption{\changed{Left column: results for the Ishigami model with a smaller basis ($d = 3, p = 12, q=1$), complementing the plots in Figure~\ref{fig:results_ishigami_smallbasis}. Right column: results for the borehole model with a smaller basis ($d = 8, p = 4, q = 1$), complementing the plots in Figure~\ref{fig:results_borehole_smallbasis}. Results for three sparse solvers and five experimental design schemes. Fifty replications.}}
	\label{fig:results_smallbasis_additional}	
\end{figure}

\changed{Due to space limitations, in Section~\ref{sec:results_nearopt}
	(Figures~\ref{fig:results_ishigami_smallbasis} and \ref{fig:results_borehole_smallbasis})
	we only showed results for two of the five solvers (OMP and \SPloo{}). In Figure~\ref{fig:results_smallbasis_additional}, we show boxplots of relative MSE against ED size for the three remaining solvers LARS, SP, and BCS. 
}

\clearpage

\section{Benchmark studies}
\label{app:benchmarkstudies}%
An overview of articles and benchmark studies comparing sparse PCE methods, including their main results, is given in Table~\ref{tab:benchmarkstudies}.

\begin{table}[tbhp]
	{\footnotesize
		\caption{\footnotesize Overview of some articles and benchmark studies comparing sparse PCE methods. The ``best method'' listed in the last column is the one delivering the smallest target error, as reported in the respective publications. Target quantities can be moments, Sobol' indices or the generalization error of the PCE surrogate. Abbreviations: `$\succ$' stands for `better than'; for the acronyms of solvers and sampling schemes see Sections \ref{sec:ED} and \ref{sec:solvers} or the respective publications; for the sampling schemes, the method given in parentheses indicates how the corresponding candidate set is created (e.g., D-opt(coh-opt) stands for D-optimal sampling based on a coherence-optimal candidate set).
			If the cited paper proposed a new method, this method is marked by a star (\markOwn).}
		\label{tab:benchmarkstudies}
		\centering
		\begingroup
		\renewcommand{\arraystretch}{1.3}
		\begin{tabular}{p{.17\textwidth} p{0.08\textwidth} p{.35\textwidth}  p{.29\textwidth}}
			\hline
			\textbf{Ref.} & \textbf{Type} & \textbf{Methods compared} & \textbf{Result} \\
			\hline
			\citep{Hampton2015, Hampton2015b} & sampling & MC, asymptotic\markOwn, and coh-opt\markOwn & coh-opt best
			\\ 
			\citep{FajraouiMarelli2017} & sampling & LHS, Sobol, D-opt(LHS), S-opt(LHS)\markOwn; \newline sequential sampling\markOwn & sequential S-opt(LHS) best; D-opt(LHS) worst
			\\ 
			\citep{Hadigol2018} & sampling (OLS) & MC, LHS, coh-opt, A-opt(coh-opt)\markOwn, D-opt(coh-opt)\markOwn, E-opt(coh-opt)\markOwn & D-opt(coh-opt) best for $p > d$; \newline MC, LHS best for $p < d$ 
			\\ 
			\citep{Jakeman2017} & sampling & MC, asymptotic, CSA\markOwn & for $p > d$: CSA better than MC and asymptotic
			\\ 
			\citep{Alemazkoor2018} & sampling & MC, coh-opt, near-opt(coh-opt)\markOwn & near-opt(coh-opt) $\succ$ coh-opt $\succ$ MC
			\\ 
			\citep{Diaz2018} & sampling & coh-opt, D-opt(coh-opt), sequential D-opt(coh-opt)\markOwn & seq. D-opt(coh-opt) $\succ$ \newline D-opt(coh-opt) $\succ$ coh-opt
			\\ 
			% 		\citep{szepietowska2018sensitivity} & sampling  & Sobol, Halton, D-opt(MC), D-opt(LHS), and others; with/without preconditioning \TODO{Omit?! Or ask her?} &  best: weighted D-opt(LHS); random subset of Hermite roots \mycomment{??}
			% 		\\ \hline 
			\citep{Dutta2020} & sampling & MC, LHS, Sobol, Importance Sampling & LHS best
			\\ 
			\citep{Hu2017} & solvers & OMP, SPGL1, BCS\citep{Babacan2010} & BCS $\succ$ OMP $\succ$ SPGL1
			\\ 
			\citep{Huan2018} & solvers & l1\_ls, SpaRSA, CGIST, FPC\_AS, ADMM with default parameters & all showed similar performance; ADMM slightly advantageous
			\\ 
			\citep{LiuWiart2020a} & solvers & OMP, LARS, rPCE\markOwn & rPCE $\succ$ LARS $\succ$ OMP
			\\ 
			\citep{Baptista2019} & solvers  & OMP, SPGL1, and two variants\markOwn{} of OMP (modified regressor selection, randomization) & best: OMP, and OMP with modified regressor selection
			\\ 
			\citep{Tarakanov2019} & solvers & OMP, LARS, Rank-PCE\markOwn & best: Rank-PCE
			\\ 
			\citep{Zhou2019c} & solvers & LARS, BCS \citep{Wipf2004b}, BCS \citep{Ji2008}, D-MORPH-reweighted \citep{Cheng2018b}, stepwise regression based on Bayesian ideas\markOwn & best: stepwise regression (based on 1 Sobol' design) 
			\\ \hline
		\end{tabular}
		\endgroup
	}
\end{table}

\end{document}